\ifx\shlhetal\undefinedcontrolsequence\let\shlhetal\relax\fi

\documentclass{amsart}%{tran-l}

\usepackage{amsmath} 
\usepackage{amssymb}

\newtheorem{theorem}{Theorem}[section] 
\newtheorem{fact}[theorem]{Fact}
\newtheorem{lemma}[theorem]{Lemma} 
\newtheorem{proposition}[theorem]{Proposition} 

\theoremstyle{definition}
\newtheorem{definition}[theorem]{Definition}
\newtheorem{problem}[theorem]{Problem} 
\newtheorem{question}[theorem]{Question} 
\newtheorem{conjecture}[theorem]{Conjecture} 
\newtheorem{thesis}[theorem]{Thesis}

\theoremstyle{remark}
\newtheorem{remark}[theorem]{Remark}

\newtheorem{discussion}[theorem]{Discussion}
\newtheorem{observation}[theorem]{Observation}

\newcommand{\pp}{{\rm pp}}
\newcommand{\pcf}{{\rm pcf}}
\newcommand{\cf}{{\rm cf}}
\newcommand{\tcf}{{\rm tcf}}
\newcommand{\Reg}{{\rm Reg}}
\newcommand{\bd}{{\rm bd}}
\newcommand{\cov}{{\rm cov}}
\newcommand{\CH}{{\rm CH}}
\newcommand{\Dom}{{\rm Dom}}
\newcommand{\Rang}{{\rm Rang}}
\newcommand{\bP}{{\mathbb P}}
\newcommand{\bQ}{{\mathbb Q}}
\newcommand{\bR}{{\mathbb R}}
\newcommand{\bS}{{\mathbb S}}
\newcommand{\cU}{{\mathcal U}}
\newcommand{\lesdot}{\mathrel{\mathord{<}\!\!\raise 
0.8 pt\hbox{$\scriptstyle\circ$}}} 
\newcommand{\bV}{{\bf V}}
\newcommand{\bL}{{\bf L}}
\newcommand{\otp}{{\rm otp}}
\newcommand{\conc}{{}^\frown\!}
\newcommand{\inv}{{\bf inv}}

\newcount\skewfactor
\def\mathunderaccent#1#2 {\let\theaccent#1\skewfactor#2
\mathpalette\putaccentunder}
\def\putaccentunder#1#2{\oalign{$#1#2$\crcr\hidewidth
\vbox to.2ex{\hbox{$#1\skew\skewfactor\theaccent{}$}\vss}\hidewidth}}
\def\name{\mathunderaccent\tilde-3 }

\numberwithin{equation}{section}

\title[On what I do not understand]{On what I do not understand (and have
something to say): Part I} 
\author{Saharon Shelah}
\address{Institute of Mathematics\\
 The Hebrew University of Jerusalem\\
 91904 Jerusalem, Israel\\
 and  Department of Mathematics\\
 Rutgers University\\
 New Brunswick, NJ 08854, USA}
\email{shelah@math.huji.ac.il}
\urladdr{http://www.math.rutgers.edu/$\sim$shelah}
\thanks{I would like to thank Alice Leonhardt for the beautiful
 typing. Publication 666\\
 Based on lectures in the Rutgers Seminar, Fall 1997 are: \S1, \S2, \S5,
 \S7, \S8}

\begin{document}

\begin{abstract} 
This is a non-standard paper, containing some problems in set theory I have
in various degrees been interested in.  Sometimes with a discussion on what
I have to say; sometimes, of what makes them interesting to me, sometimes
the problems are presented with a discussion of how I have tried to solve
them, and sometimes with failed tries, anecdote and opinion.  So the
discussion is quite personal, in other words, egocentric and somewhat
accidental. As we discuss many problems, history and side references are
erratic, usually kept at a minimum (``see $\ldots$" means: see the
references there and possibly the paper itself).

The base were lectures in Rutgers Fall '97 and reflect my knowledge then.
The other half, concentrating on model theory, will subsequently appear.  I
thank Andreas Blass and Andrzej Ros{\l}anowski for many helpful comments.
\end{abstract}

\maketitle

\tableofcontents
\eject

\section{Cardinal problems and pcf}
Here, we deal with cardinal arithmetic as I understand it (see
\cite{Sh:400a} or \cite{Sh:g}), maybe better called cofinality arithmetic
(see definitions below). What should be our questions? Wrong questions
usually have no interesting answers or none at all.

\noindent Probably the most popular\footnote{that is, most people who are
aware of this direction, will mention it, and probably many have tried it to
some extent} question is:

\begin{problem}
\label{1.1} 
Is $\pp(\aleph_\omega)<\aleph_{\omega_1}$? 
\end{problem}

Recall:
\begin{definition}
\label{1.2} 
Let ${\mathfrak a}$ be a set of regular cardinals (usually $|{\mathfrak a}|
< \min({\mathfrak a})$). We define  
\begin{enumerate}
\item $\pcf({\mathfrak a})=\{\cf(\prod{\mathfrak a}/D):D$ is an ultrafilter
on ${\mathfrak a}\;\}$.
\item $\cf(\prod{\mathfrak a})=\min\{|F|:F\subseteq\prod{\mathfrak a}$ and
$(\forall g\in\prod {\mathfrak a})(\exists f\in F)(g\le f)\}$. 
\item For a filter $D$ on ${\mathfrak a}$, $\tcf(\prod{\mathfrak a}/D)=
\lambda$ means that in $\prod{\mathfrak a}/D$ there is an increasing cofinal
sequence of length $\lambda$. 
\item For a singular cardinal $\mu$ and a cardinal $\theta$ such that
$\cf(\mu)\leq\theta<\mu$ let 
\[\begin{array}{ll}
\pp_\theta(\mu)=\sup\big\{\tcf(\prod{\mathfrak a}/I):&{\mathfrak a}\subseteq
\Reg\cap\mu,\ |{\mathfrak a}|<\min({\mathfrak a}),\ \sup({\mathfrak a})=
\mu,\\ 
&I \mbox{ an ideal on ${\mathfrak a}$ such that } J^{\bd}_{\mathfrak a}
\subseteq I,\mbox{ and }|{\mathfrak a}|\leq\theta\},
\end{array}\]
were for a set $A$ of ordinals with no last element, $J^{\bd}_A$ is the
ideal of bounded subsets of $A$. 
\item Let $\pp(\mu)=\pp_{\cf(\mu)}(\mu)$.
\item We define similarly $\pp_\Gamma(\mu)$ for a family (equivalently: a
property) $\Gamma$ of ideals; e.g., $\Gamma(\theta,\tau)=$ the family of
$(<\tau)$-complete ideals on a cardinal $<\theta$, $\Gamma(\theta)=\Gamma(
\theta^+,\theta)$. 
\end{enumerate}
\end{definition}

\begin{definition}
\label{1.3}
\begin{enumerate}
\item For a partial order $P$, 
\[\cf(P)=\min\{|Q|:Q\subseteq P\mbox{ and }(\forall p\in P)(\exists q\in
Q)(p \le q)\}.\] 
\item For cardinals $\lambda,\mu,\theta,\sigma$,
\[\begin{array}{lr}
\cov(\lambda,\mu,\theta,\sigma)=\min\{|{\mathcal A}|:&{\mathcal A}
\subseteq [\lambda]^{<\mu}\mbox{ and any }a\in [\lambda]^{<\theta}\mbox{ is 
included }\ \ \\
&\mbox{in the union of $<\sigma$ members of }{\mathcal A}\ \}.
  \end{array}\]
\end{enumerate}
\end{definition}

Problem \ref{1.1} is for me the right form of 

\begin{question}
\label{1.4} 
\begin{enumerate}
\item Assume $\aleph_\omega$ is strong limit. Is $2^{\aleph_\omega}<
\aleph_{\omega_1}$?  
\item Assume $2^{\aleph_0}<\aleph_\omega$. Is $(\aleph_\omega)^{\aleph_0}
<\aleph_{\omega_1}$? 
\end{enumerate}
\end{question}
Why do I think \ref{1.1} is a better form?  Because we know that:
\begin{enumerate}
\item[$(*)_1$] If $\aleph_\omega$ is strong limit, then $2^{\aleph_\omega}=
(\aleph_\omega)^{\aleph_0}$ (classical cardinal arithmetic). 
\item[$(*)_2$]  $\pp(\aleph_\omega)=\cf([\aleph_\omega]^{\aleph_0},
\subseteq)$ (see \cite{Sh:E12}),
\item[$(*)_3$]  $\aleph_\omega^{\aleph_0}=2^{\aleph_0}+\cf([
\aleph_\omega]^{\aleph_0},\subseteq)$ (trivial). 
\end{enumerate}
So the three versions are equivalent and say the same thing when they say
something at all, but Problem \ref{1.1} is always meaningful. 

To present what I think are central problems, we can start from what I
called the solution of the ``Hilbert's first problem", see \cite{Sh:460}
(though without being seconded). 

\begin{theorem}
\label{1.5} 
For $\lambda\ge\beth_\omega$, there are $\kappa<\beth_\omega$ and ${\mathcal
P} \subseteq [\lambda]^{<\beth_\omega}$, $|{\mathcal P}|=\lambda$ such that
every $A\in [\lambda]^{<\beth_\omega}$ is equal to the union of $<\kappa$
members of ${\mathcal P}$.  
\end{theorem}
So ${\mathcal P}$ is ``very dense''. E.g., if $c:[\lambda]^n\longrightarrow
\beth_n$ {\em then\/} for some $B_m\in {\mathcal P}$ (for $m<\omega$), the
restrictions $c\restriction [B_m]^n$ are constant and $|B_m|=\beth_m$. We
can replace $\beth_\omega$ by any strong limit cardinal $>\aleph_0$. 

In \cite[\S8]{Sh:575} the following application of \ref{1.5} to the theory
of Boolean Algebras is proved: 

\begin{theorem}
\label{1.6}
If $B$ is a c.c.c.~Boolean algebra and $\mu=\mu^{\beth_\omega}\le |B|\le
2^\mu$, {\em then} $B$ is $\mu$-linked, i.e., $B\backslash\{0\}$ is the
union of $\mu$ sets of pairwise compatible elements.
\end{theorem}
\noindent (See also \cite{Sh:92}, \cite{Sh:126} and Hajnal, Juh\'asz and
Szentmiklossy \cite{HaJuSz}.) 

We also have the following application: 

\begin{theorem}
[See \cite{Sh:454a}]  
\label{1.7}
If $X$ is a topological space (not necessarily $T_2$) with $\lambda$ points,
$\mu\le\lambda<2^\mu$ and $>\lambda$ open sets and $\mu$ is strong limit of
cofinality $\aleph_0$, {\em then} $X$ has $\ge 2^\mu$ open sets.
\end{theorem}
Another connection to the general topology is the following 

\begin{definition}
\label{1.8}
For topological spaces $X,Y$ and a cardinal $\theta$, write $X\rightarrow
(Y)^1_\theta$ iff for every partition $\langle X_i:i<\theta\rangle$ of $X$
into $\theta$ parts, $X$ has a closed subspace $Y'$ homeomorphic to $Y$
which is included in one part of the partition. 
\end{definition}
Arhangel'skii asked whether for every compact Hausdorff space $X$, 
\[X\nrightarrow (\mbox{Cantor discontinuum})^1_2.\] 
Arhangel'skii's problem $+ \neg\CH$ is sandwiched between two pcf statements
of which we really do not know whether they are true. If, for simplicity
$2^{\aleph_0}\ge\aleph_3$, then e.g.:

\begin{enumerate}
\item[$(*)_1$] if for no ${\mathfrak a}$, ${\mathfrak a}\subseteq\Reg$,
where $\Reg$ is the class of regular cardinals, $|{\mathfrak
a}|\ge\aleph_2$, $\prod{\mathfrak a}/[{\mathfrak a}]^{\le\aleph_0}$ is
$\sup({\mathfrak a})$--directed, {\em then\/} the answer is: for every
Hausdorff space $X$, we have $X\nrightarrow$ (Cantor discontinuum)$^1_2$ and
more. 
\item[$(*)_2$] If for some ${\mathfrak a}\subseteq\Reg\setminus
2^{<\kappa}$, $|{\mathfrak a}|=2^{\aleph_0}\le\kappa$ and $\prod{\mathfrak
a}/[{\mathfrak a}]^{\leq\aleph_0}$ is $\sup({\mathfrak a})$--directed, {\em
then\/} in some forcing extension there exists a zero-dimensional Hausdorff
space $X$ such that $X \rightarrow$ (Cantor discontinuum)$^1_2$. 

\noindent The Stone-\u Cech compactification of this space gives a negative
answer to Arhangel'skii's question.  
\end{enumerate}
(On the problem, see \cite{Sh:460} and more in \cite{Sh:668}.) 

However, we can start from inside pcf theory.

\begin{problem}
\label{1.10} 
Is $\pcf({\mathfrak a})$ countable for each countable set of cardinals? 
\end{problem}
This seems to me more basic than \ref{1.1}, but yet \ref{1.1} is weaker. I
think it is better to look at the battlefield between independence by
forcing from large cardinals and proofs in ZFC (I would tend to  say
between the armies of Satan and God but the armies are not disjoint).

The advances in pcf theory show us ZFC is more powerful than expected
before. I will try to give a line of statements on which both known methods
fail -- so far. 

\begin{conjecture}
\label{1.11}
If ${\mathfrak a}$ is a set of regular cardinals $>|{\mathfrak a}|$, {\em
then\/} for no inaccessible $\lambda$ the intersection $\lambda\cap
\pcf({\mathfrak a})$ is unbounded in $\lambda$. 
\end{conjecture}

\begin{conjecture}
\label{1.12}
For every $\mu\ge\aleph_\omega$, for every $\aleph_n<\aleph_\omega$ large
enough there is no $\lambda<\mu$ of cofinality $\aleph_n$ such that
$\pp_{\Gamma(\aleph_n)}(\lambda)>\mu$ (or replace $\aleph_n<\aleph_\omega$
by $\aleph_\alpha<\aleph_{\omega^2}$ or even $\aleph_\alpha<
\aleph_{\omega_1}$, or whatever).  
\end{conjecture}

\begin{conjecture}
\label{1.13}
\begin{enumerate}
\item[(A)] It is consistent, for any uncountable $\theta$ (e.g.,
$\aleph_1$), that for some $\lambda$  
\[\theta\leq |\{\mu<\lambda:\cf(\mu)=\aleph_0,\ \pp(\mu)>\lambda\}|.\]
\item[(B)] It is consistent that for some $\lambda$, the set
\[\{\mu<\lambda:\cf(\mu)>\aleph_0,\ \pp_{\aleph_1\text{--complete}}(\mu)>
\lambda\}\]
is infinite. 
\end{enumerate}
\end{conjecture}
Those three conjectures seem to be fundamental. Note that having
ZFC-provable answer in \ref{1.11}, \ref{1.12}, but independent answer for
\ref{1.13} are conscious choices. For all of those problems, present methods
of independence fail, and in addition they are known to require higher
consistency strength. Of course, we can concentrate on other variants;
e.g., in \ref{1.13}(B) use $\theta$ instead of $\aleph_0$.  

Other problems tend to be sandwiched between those, or at least those more
basic problems are embedded into them.  E.g., \ref{1.12} implies that in
\ref{1.5} we can replace $\beth_\omega$ by $\aleph_\omega$ if we replace
equal by included (or demand $\lambda\ge\sum\limits_{n<\omega}
2^{\aleph_n}$) and this implies $|{\mathfrak a}|\le\aleph_0\ \Rightarrow\ 
|\pcf({\mathfrak a})|\le\aleph_\omega$, while e.g.~$|{\mathfrak a}|\le
\aleph_{\omega n}\ \Rightarrow\ |\pcf({\mathfrak a})|\le\aleph_{\omega
n+\omega}$ implies the analog of \ref{1.5} for $\aleph_{\omega^2}$, see
\cite{Sh:460}, \cite{Sh:513}. See \cite{Sh:513} for more on the ZFC side; it
is very helpful in preventing futile attempts to force. 

Note that $\pp(\aleph_\omega)>\aleph_{\omega_1}$, {\em implies\/} that for
some countable ${\mathfrak a}$, $\pcf({\mathfrak a})$ is uncountable,
which {\em implies\/} that clause (A) from Conjecture \ref{1.13} holds. Also
$\pp(\aleph_\omega)>\aleph_{\omega_2}$ {\em implies\/} that for some
countable ${\mathfrak a},|\pcf({\mathfrak a})|\ge \aleph_2$ which {\em
implies} that clause (B) of Conjecture \ref{1.13} fails). 

So there is no point to try to prove CON($\pp(\aleph_\omega)>\aleph_{
\omega_1}$) before having the consistency of \ref{1.13}(A) and, thus,
CON($\pp(\aleph_\omega)>\aleph_{\omega_1}$) is a more specialized case.
(Also if we look at the earlier history of consistency proofs -- clearly
there is no point to start with Problem \ref{1.1}). 

In Conjecture \ref{1.11} the situation (which we say is impossible) may look
bizarre, as $\pcf({\mathfrak a})$ is extremely large. Of course, much better
is $|\pcf({\mathfrak a})|<$ ``first inaccessible $>|{\mathfrak a}|$'' and
even $|\pcf({\mathfrak a})|\le |{\mathfrak a}|^{+\omega}$, which follows
from Conjecture \ref{1.12}. Of course, replacing in \ref{1.11}, ``$\lambda$
inaccessible'' by ``$\lambda$-Mahlo'' is still a very important conjecture
while getting $\pcf({\mathfrak a})<$ ``the first fixed point $>|{\mathfrak
a}|$'' is much better, so why from all variants of \ref{1.11}, those we have
just mentioned and others, ``the accumulation inaccessible'' was chosen? The
point is that it implies 
\begin{enumerate}
\item[$(*)$]  $\cf(\prod\pcf({\mathfrak a}))=\cf(\prod{\mathfrak a})$,

if ${\mathfrak a}$ is a set of regular cardinals $>|{\mathfrak a}|$
\end{enumerate}
(see \cite[Ch.VIII,\S3]{Sh:g}, \cite{Sh:E11}; note that in the notation of
\cite{Sh:E11}, conjecture \ref{1.11} says that $\pcf({\mathfrak a})\in   
J^*_{\mathfrak a}$.) If there is a failure of Conjecture \ref{1.11} then
consistently $(*)$ fails. We can force by $(<\lambda)$--complete forcing
iterating adding $f\in\prod(\pcf({\mathfrak a})\cap\lambda)$ dominating the
old product (or for any $\mu$, just adding $\mu$ many $\lambda$--Cohen
functions, i.e., forcing with 
\[\{f:f\mbox{ is a partial function from $\mu$ to }\lambda,\
\Dom(f)|<\lambda\}).\]
So \ref{1.11} denotes the significant dividing line between chaos and
order. 

Concerning the last conjecture \ref{1.13}, maybe the proofs in Gitik and
Shelah \cite{GiSh:344} are relevant. There we force for hypermeasurable
cardinals $\kappa_0<\kappa_1<\ldots<\kappa_n$ with a forcing which makes
each $\kappa_i$ hypermeasurable indestructible under reasonable forcing
notions, including those which may add new Prikry sequences of ordinals
$>\kappa_\ell$ of length $<\kappa_\ell$. (So in this case supercompact
cannot serve, unlike in many proofs which do with hypermeasurable cardinals
what is relatively easy to do with supercompact cardinals.) Let $\lambda=
\lambda^{<\lambda}>\kappa_n$, $\theta_\ell<\lambda$. Then we blow up
$2^{\kappa_n}$ to $\lambda$, change $\cf(\kappa_n)$ to $\theta_n$; blow up
$2^{\kappa_{n-1}}$ to $\lambda$, change $\cf(\kappa_{n-1})$ to
$\sigma_{n-1}$, etc.  

The point is that when we arrive to $\kappa_i$ the forcing so far is fairly
``$\kappa_i$-complete for pure extensions'', etc, so does not destroy
``$\kappa_i$ is $\lambda$-hypermeasurable''. So for Conjecture \ref{1.13} we
(fix the desired cofinality $\theta$ and we) need to do it not $n$ times but
$\theta^*$ times ($\theta^*=$ inverse order of $\theta$) so we need
``anti-well-founded iteration''. In other words, we have $\langle\kappa_i:i
<\theta\rangle$ increasing; $\kappa_i$ is $\lambda$-hypermeasurable 
indestructible (necessarily in a strong way), and $\lambda>\sum\limits_{j<
\theta}\kappa_j$. 

\noindent{\sc First try}:\\
We may try to define by induction on $i<\theta$, $\lesdot$--decreasing
sequences $\bar{\bP}^i=\langle \bP^i_j:j\le i\rangle$ of forcing notions
such that $|\bP^i_0|=\lambda_i$, 
\[\Vdash_{\bP^i_0}\mbox{`` }\cf(\kappa_j)=\aleph_0\ \mbox{ for }j<i\mbox{
''}\]
(or whatever fixed value, but $\aleph_0$ is surely easier), $\bP^i_j$ is
$\kappa^+_j$-c.c., purely $\kappa_j$-complete, $\bP^i_0$ makes
$\pp(\kappa_j) = \lambda$ for every $j<i$. 

In successor stage - no problem:\quad $i=j+1$ and
\[
\begin{array}{lcl}
\bP^{i+1}_{i+1}&=&(\mbox{blowing up } 2^\kappa\mbox{ to }\lambda\mbox{
changing cf}(\kappa_j)\mbox{ to }\aleph_0)\\ 
\bP^{i+1}_j&=&\bP^{i+1}_{i+1} * \bP^i_j.
\end{array}\]
Not good enough: $\bP^{i+1}_i$ changes the definition of: ``blowing up
$2^{\kappa_j}$ to $\lambda$'' as there are more $\omega$--sequences. So we
should correct ourselves to $|\bP^i_j|=\sum\limits_{\zeta<i}\kappa_\zeta$: 
\[\bP^i_j\mbox{ blows up } 2^{\kappa_j}\mbox{ to essentially }\sum_{\zeta<i}
\kappa_\zeta. \]
So we have to prove the forcing notions extend as they should. If
$\bar{\bP}^i$ is defined, there is no problem to choose an appropriate
$\bP^{i+1}_{i+1}$. Now for each $j\le i$ separately we would like to choose
$\bP^{i+1}_j$ to be a $\lesdot$-extension of $\bP^i_j$ and of
$\bP^{i+1}_{i+1}$, but we have to do it for all $j\le i$ together. The limit
case seems harder. 

\[* \qquad * \qquad *\]

Why, in \ref{1.13}(A), do we have $\theta\ge\aleph_1$? Moti Gitik shows
consistency for $\theta=\aleph_0$ by known methods. 

\noindent{\sc Audience Question:}\quad How dare you conjecture ZFC can show
\ref{1.11}, \ref{1.12}? 

For Conjecture \ref{1.13} I have a scenario for an independence proof
(outlined above). For \ref{1.11} and \ref{1.12} the statements imply there
is quite a complicated pcf structure you necessarily drag with you. So it is
reasonable to assume that if we shall know enough theorems on the pcf
structure we shall get a contradiction.  Of course, those arguments are not
decisive. 

\[* \qquad * \qquad *\]

Traditionally, remnants of GCH have strongly influenced the research on
cardinal arithmetic, so e.g.~people concentrate on the strong limit case,
see \cite[AG]{Sh:g}, \cite{Sh:E12}; probably also it was clear what to do
and easier. On the other hand, \cite{Sh:g} aims to get ``exponentiation-free
theorems'', so we put forward:

\begin{thesis}
\label{1.14}
``Everything'' is expressibly by cases of $\pp_J$ (and $2^\kappa$ for
$\kappa$ regular). 
\end{thesis}
E.g. in \cite[\S2]{Sh:589} this is done to the tree power of $\lambda$,
\[\lambda^{\langle\kappa,{\rm tr}\rangle}=\sup\{|{\lim}_\kappa(T)|:T\mbox{ a
tree with } \le\lambda \mbox{ nodes and }\kappa \mbox{ levels}\},\] 
where $\lim_\kappa(T)$ is the set of $\kappa$--branches of $T$ (well, using
$\kappa^{\langle\kappa,{\rm tr}\rangle}$ for regular $\kappa$, which is
malleable by forcing, a relative of $2^\kappa$ for $\kappa$ regular).

But maybe there are also forcing proofs by which we can get interesting
situations say below the continuum, whose strong limit counterparts are
false, or have bigger consistency strength, or at least are harder to prove.
The known forcing proofs may be open to such variations, e.g., when we add
many Prikry sequences to one $\kappa$ we may have the order between them
such that every condition decides little about it. The following problem may
be relevant to \ref{1.14}, and anyhow is a central one.  

\begin{problem}
\label{1.15}
For a singular cardinal $\mu>\theta=\cf(\mu)$, is 
\[\cov(\mu,\mu,\theta^+,\theta)=\pp_{\Gamma(\theta^+,\theta)}(\mu)\ ?\]
\end{problem}
Note that other cases of cov can be reduced to those above. Now, this is
almost proved:  it holds when $\theta=\cf(\mu)>\aleph_0$. Furthermore, if
$\mu$ is strong limit, $\aleph_0=\cf(\mu)$ and the two expressions in
\ref{1.15} are not equal, both are quite large above $\mu$ as in Conjecture
\ref{1.13}. Also, e.g., for a club of $\delta<\omega_1$
\[\cov(\beth_\delta,\beth_\delta,\aleph_1,\aleph_0)=\pp(\beth_\delta),\] 
(see \cite{Sh:E12}, the ``$\beth_{\omega_1}$'' can be weakened to strong
limit in cov sense). But

\begin{question}
\label{1.16}
Can we force that there is $\mu<2^{\aleph_0}$ such that $\cf(\mu)=\aleph_0$
and $\cov(\mu,\mu,\aleph_1,\aleph_0)>\pp(\mu)$?
\end{question}
\noindent [Why $<2^{\aleph_0}$? As blowing up the continuum does not change
the situation, proving the consistency for $\mu<2^{\aleph_0}$ can be only
easier. But for $\mu<2^{\aleph_0}$ maybe it is even consistent that 
\[\cov(\mu,\mu,\aleph_1,\aleph_0)>\mu^+=\pp(\mu),\]
that is, by our present ignorance, it is even possible that the behaviour
below the continuum is different than above it.]

Note that all cases of $\lambda^\kappa$ can be reduced to cases of
$2^\theta$, $\theta$ regular, and $\cf([\mu]^{\le\theta},\subseteq)$ where 
$\mu>\cf(\theta)=\theta\ge\cf(\mu)$.\\
Why? If $\kappa$ is regular, $\lambda\le 2^\kappa$ then $\lambda^\kappa=
2^\kappa$. If $\kappa$ is regular and $\lambda>2^\kappa$ then
$\lambda^\kappa= \cf([\lambda]^{\le\kappa},\subseteq)$. So assume 
$\kappa$ is singular and let $\sigma=\cf(\kappa)$ and $\kappa=
\sum\limits_{i<\sigma} \kappa_i$, where each $\kappa_i$ is regular and
$\sigma<\kappa_i<\kappa$, so $\lambda^\kappa=\lambda^{\sum\{\kappa_i:i<
\sigma\}}=\prod\limits_{i<\sigma}\lambda^{\kappa_i}$. Thus, if $\lambda\le
2^\kappa$ then 
\[\lambda^\kappa=2^\kappa=\prod\limits_{i<\sigma}2^{\kappa_i}=(
\sum\limits_{i<\sigma}2^{\kappa_i})^\sigma=\cf([\sum\limits_{i<\sigma}
2^{\kappa_i}]^\sigma,\subseteq).\] 
Lastly, if $\lambda>2^\kappa$, 
\[\lambda^\kappa=\prod\limits_{i<\sigma}\lambda^{\kappa_i}=(\sum\limits_{i< 
\sigma}\lambda^{\kappa_i})^\sigma=(\max_{i<\sigma}\;\lambda^{\kappa_i}
)^\sigma=\max_{i<\sigma}\;\lambda^{\kappa_i}=\max_{i<\sigma}\;\cf([\lambda
]^{\kappa_i},\subseteq)\]
(on the third equality see Hajnal and Hamburger \cite{HH},
\cite[2.11(4),p. 164]{Sh:233}). 

If the answer to \ref{1.15} is yes, then we can reduce all cases of
$\lambda^\kappa$ and of $\cov$ to statements on cases of $\pp$.

\begin{problem}
\label{1.16a}
If $\cf(\mu)=\aleph_0$, is $\pp(\mu)$ equal to $\pp^{\text{cr}}_{J^{
\bd}_\omega}(\mu)$, where
\[\begin{array}{ll}
\pp^{\text{cr}}_{J^{\bd}_\omega}(\mu)=\sup\{\lambda:&\mbox{for some
increasing sequence }\langle\lambda_n: n<\omega\rangle\mbox{ of regular}\\ 
&\mbox{cardinals converging to }\mu\mbox{ we have }\lambda=\tcf(
\prod\limits_{n<\omega}\lambda_n/J^{\bd}_\omega)\}\ ?
  \end{array}\]
A variant is: except when $\pp^{\text{cr}}_{J^{\bd}_\omega}(\mu)$ has
cofinality $\aleph_0$ and $\pp(\mu)$ is its successor.
\end{problem}
By pcf calculus, if $\pp(\mu)<\mu^{+\omega_1}$ then this is true. Similarly,
if $\theta<\mu_0<\mu$ and
\[(\forall\mu')([\cf(\mu')\le\theta\ \&\ \mu'\in(\mu_0,\mu)\quad\Rightarrow
\quad \pp(\mu')<\mu^{+\theta^+}]\]
then $\pp(\mu)=\pp_\theta(\mu)$ and see \cite[6.5]{Sh:430}. Also, by
\cite[Part C]{Sh:E12}, e.g., for a club of $\delta<\omega_1$, $\mu=
\beth_\delta$ satisfies the conclusion. 

\[* \qquad * \qquad *\]

On pcf for set theories with weak versions of Choice (say ${\rm DC}_\kappa$,
the dependent choice of length $\kappa$) see \cite{Sh:497}.

\begin{problem}  
\label{1.16A}
Develop combinatorial set theory generally and, in particular, pcf theory
using only little choice (say ${\rm DC}_\kappa$).
\end{problem}
Inner model theory and descriptive set theory are not hampered by lack of
choice, and much was done on variants of the axiom of choice. \cite{Sh:497}
may be a beginning of combinatorial set theory, and pcf in particular;
i.e., it is enough to show that there are interesting theorems. In
particular 

\begin{question}
\label{1.16Aa}
\begin{enumerate}
\item[(a)]  Does ${\rm DC}_\kappa$ for $\kappa$ large enough imply the
existence of a proper class of regular cardinals? 
\item[(b)]  Does ${\rm DC}_\kappa$ for $\kappa$ large enough imply that for
a class of $\lambda$, ${\mathcal P}(\lambda)$ is not the union of $<\lambda$
sets, each of cardinality $\le\lambda$?
\end{enumerate}
\end{question}
See more in \cite{Sh:497}.  Gitik \cite{Gi80} had proved 
\[{\rm CON}\big((\forall\delta)(\cf(\delta)\le\aleph_0)\big)\]
relative to suitable large cardinals. Woodin asked if 
\[{\rm CON}\left({\rm DC}_{\aleph_0}+(\forall\delta)(\cf(\delta)\le
\aleph_1)\right).\]
Specker asked if, consistently, for every $\lambda$, for some $\langle
A_n:n<\omega\rangle$ we have ${\mathcal P}(\lambda)=\bigcup\limits_n A_n$,
$|A_n|\le\lambda$.

\[* \qquad * \qquad *\]

On how the problem of the existence of universal objects is connected to pcf
see Kojman and Shelah \cite{KjSh:409}, and \cite{Sh:552}, \cite{Sh:622}. The 
following conjecture will simplify the answers: 

\begin{conjecture}
\label{1.16B}
For every limit of limit cardinals $\mu$, for arbitrarily large regular
$\lambda<\mu$, we have  
\[(\forall^* \mu_1<\mu)[\cf(\mu_1)=\lambda\ \ \Rightarrow\ \ \pp_{\Gamma(
\lambda)}(\mu_1)<\mu],\]
where $\forall^*$ means ``for every large enough''.
\end{conjecture}
After we learned that, on the one hand, $2^{\aleph_\omega}$ ($\aleph_\omega$
strong limit) has a bound (in fact, every $2^{\aleph_\delta}$, if
$\aleph_\delta$ is strong limit $>|\delta|$, in \cite[Ch.XIII]{Sh:b}), and
on the other hand there are bounds for $2^\mu$, $\mu$ the $\omega_1$-th fix
point (when $\mu$ is strong limit or less), it becomes natural to ask: 

\begin{conjecture}
\label{1.17}
If $\aleph_\delta$ is the first fix point (i.e., the first such that
$\aleph_\delta=\delta$, so it has cofinality $\aleph_0$), {\em then}
$\pp(\aleph_\delta)<(2^{\aleph_0})^+$--th fix point.\\
(Even assuming GCH below $\aleph_\delta$ and proving just
$\pp(\aleph_\delta)<$ ``the first inaccessible'' is good, but
``$<\omega_4$-th fix point'' is better, and ``$<\omega_1$-th fix point'' is
best, but seems pointless to ask as long as \ref{1.1} is open).
\end{conjecture}
Note that we almost know: if $\aleph_\delta$ is the $\omega_1$-th fix point
(strong limit), then $\pp(\aleph_\delta)<\omega_4$--th fix point, we know it
if the answer to \ref{1.11} is yes, see \cite{Sh:420} and see
\cite[Ch.V]{Sh:g}. 

\[* \qquad * \qquad *\]

Traditionally we have asked: ``can we find all the laws of cardinal
arithmetic?'' This had been accomplished for regular cardinals, and we
prefer  
\begin{problem}
\label{1.18} 
Find all the rules of the pcf calculus or at least find more (or show that
the set of rules is inherently too complicated). 
\end{problem}
Note: if for simplicity $|\pcf({\mathfrak a})|<\min({\mathfrak a})$, then on
$\pcf({\mathfrak a})$ the pcf structure is naturally a compact topology:
${\mathfrak b}$ is closed iff ${\mathfrak b}=\pcf({\mathfrak b})$, and the
theorem on existence of generators $\langle {\mathfrak b}_\lambda:\lambda
\in\pcf({\mathfrak a})\rangle$ says that the topology is a particularly nice
one. If \ref{1.11} holds this is true whenever $|{\mathfrak a}|<\min(
{\mathfrak a})$ (see \cite{Sh:E11}). 

There may well be some ``global phenomena''. Also there may be special
behaviour near 
\[\min\{\lambda:\mbox{for some } A\subseteq\lambda,\mbox{ there is no
indiscernible class for }K[A]\},\]
as above it the covering theorem (Dodd and Jensen \cite{DJ2}) shows that
cardinal arithmetic is trivial. On the other hand, on the behaviour below
it, see \cite[Ch.V]{Sh:g}.  

An extreme case of our non-understanding concerning global behavior is:
\begin{question}
\label{1.18A}
Is it possible that:
\begin{quotation}
if ${\mathfrak a}$ is a set of odd [even] regular cardinals $>|{\mathfrak
a}|$,\\ 
{\em then} every $\theta\in\pcf({\mathfrak a})$ is odd [even]?
\end{quotation}
(where $\aleph_{2 \alpha}$ is even and $\aleph_{2 \alpha +1}$ is odd).
\end{question}

Instead of looking more on $\pp(\aleph_\omega)$ we may ask if the best
result was derived from the known laws of cardinal arithmetic. 

\begin{question}
\label{1.19}
Let $\ell<4$. Can there be $\delta\in [\omega_\ell,\omega_{\ell+1})$ and a
closure operation $c\ell$ on ${\mathcal P}(\delta+1)$ such that all the
rules used in the proof of $\pp(\aleph_\omega)<\aleph_{\omega_4}$ hold?
(see Jech and Shelah \cite{JeSh:476}). 
\end{question}

\begin{question}
\label{1.19A} 
\begin{enumerate}
\item Characterize the possible sequences 
\[\langle J_{<\theta}[\{\aleph_n:n\in [1,\omega)\}]:\ \theta\in\pcf\{
\aleph_n:\ n\in [1,\omega)\}\rangle.\] 
\item For every ordinal $\gamma$ characterize the possible $\langle J_{
<\theta}[{\mathfrak a}]:\ \theta\in\pcf({\mathfrak a})\rangle$ up to
isomorphism when $\otp({\mathfrak a})=\gamma$. 
\end{enumerate}
\end{question}
\noindent [For ${\mathfrak a}'$, ${\mathfrak a}''$ we have an isomorphism if
there there is a one--to--one order preserving $f:{\mathfrak a}'
\longrightarrow {\mathfrak a}''$ such that $\{J_{<\theta}[{\mathfrak a}'']:
\theta\in \pcf({\mathfrak a}'')\}=\big\{\{f[{\mathfrak b}]:{\mathfrak 
b}\in J_{<\theta}[{\mathfrak a}']\}:\theta\in\pcf({\mathfrak a}')\big\}$.]   

\[* \qquad * \qquad *\]

I feel that

\begin{thesis}
\label{1.20}
Proving a theorem from ZFC + ``cardinal arithmetic assumptions'' is a ``semi
ZFC result". 
\end{thesis}
This view makes proofs from cases of the failure of the SCH related to the
thesis below more interesting.  

\begin{thesis}
\label{1.21}
Assumptions on the failure of GCH (and even more so, of SCH) are good
assumptions, practical ones, in the sense that from them you can deduce
theorems. 
\end{thesis}
Traditionally this is how instances of GCH were treated (with large
supporting evidence). Clearly \ref{1.21} may be supported by positive
evidence (though hard to refute), whereas \ref{1.20} remains a matter of
taste. So Magidor would stress looking at ``existence of a large cardinal''
as semi-ZFC axioms (unlike some randomly chosen consistent theorems), which
seems to mean in our terminology that we will look at consequences of it as 
semi-ZFC theorems. Jensen stresses that showing $\psi$ holds in a universe
with structure is much better than mere consistency (so the fine structure
in $\bL$ was the only one we know of at one time, but e.g. $K$ is no less
good than $\bL$; the statement in \cite{Sh:E16} was inaccurate).

I agree with both, just to a lesser degree. Kojman criticized \ref{1.20}
saying cases of failure of SCH are large cardinal assumption in disguise;
and I agree that $2^\lambda > \aleph_{\lambda^{+4}}$ is a weaker assumption
than $2^{\beth_\omega}>\beth^+_\omega$, but I still stick to \ref{1.20}. We
may hope to really resolve problems by partitioning to cases according to
what the cardinal arithmetic is.

\[* \qquad * \qquad *\]

\noindent {\sc Discussion:}\quad The following should be obvious, but I have
found that mentioning them explicitly is helpful. Assume, e.g., that
$\cf(\mu)=\aleph_0$, $\pp(\mu)>\mu^{+\omega_n}$, $n>0$ ($\omega_n$ for 
simplicity) and let 
\begin{enumerate}
\item[$(*)_{\mu,n}$] for stationarily many $\delta<\omega_n$ of cofinality
$\aleph_0$, $\pp_{\aleph_n}(\mu^{+\delta})<\mu^{+\omega_n}$\\
(a ``soft'' assumption, see \cite[IX,\S4]{Sh:g}).
\end{enumerate}
{\em Then\/} we can find pairwise disjoint countable ${\mathfrak a}_i
\subseteq\Reg\cap\mu$ unbounded in $\mu$ and $\alpha_i<\omega_n$
successor, strictly increasing and such that 
\[\mu^{+\alpha_i}=\max\;\pcf({\mathfrak a}_i),\quad \mu^{+\alpha_i}\notin
\pcf(\bigcup_{j\ne i}{\mathfrak a}_j),\]
moreover $\mu'<\mu\ \Rightarrow\ \mu^{\alpha_i}=\max\pcf({\mathfrak a}_i
\setminus\mu')$. 

\noindent [Why? We can find, by the assumption and Fodor Lemma,
$\alpha^*<\omega_n$ such that 
\[\alpha\in [\alpha^*,\omega_n)\ \ \Rightarrow\ \ \max\;\pcf\{\mu^{+
\beta+1}:\beta\in (\alpha_0,\alpha)\}<\mu^{+ \omega_n}.\]
By the assumption $\pp(\mu)\ge\mu^{+\omega_n}$, there is ${\mathfrak a}
\subseteq \mu \backslash\omega_n$, $|{\mathfrak a}|=\aleph_n$ such that
$\alpha<\omega_n\ \Rightarrow\ \mu^{+\alpha+1}\in\pcf({\mathfrak a})$. First
assume $2^{\aleph_n}<\mu$, so without loss of generality $\min({\mathfrak
a})> 2^{\aleph_n}$, and we have a smooth closed generating sequence $\langle
{\mathfrak b}_\lambda:\lambda\in\pcf({\mathfrak a})\rangle$ for
$\pcf({\mathfrak a})$ (so ${\mathfrak b}_\lambda\subseteq\pcf({\mathfrak
a})$, etc.). Now choose by induction on $i<\omega_n$ pairs $(\alpha_i,
{\mathfrak b}'_i)$ as follows. If $\langle\alpha_j:j<i\rangle$ has been
defined, we know that 
\[\max\;\pcf\{\mu^{+\beta}:\beta\mbox{ successor, }\alpha^*\le\beta\le
(\alpha^* +2)\cup\bigcup_{j<i}\alpha_j\}<\mu^{+ \omega_n},\]
and hence we can find $m_i<\omega$ and successor ordinals $\gamma^i_\ell\in 
[\alpha^*,(\alpha^* +2)\cup\bigcup\limits_{j<i}\alpha_j]$ (for $\ell<m_i$)
such that 
\[\{\mu^{+\beta}:\beta\mbox{ a successor, }\alpha^*\le\beta\le(\alpha^*+1)
\cup \bigcup_{j<i}\alpha_j\}\subseteq\bigcup_{\ell<m_i} {\mathfrak
b}_{\gamma^i_\ell}.\]
Let $\alpha_i<\omega_n$ be the minimal successor such that 
\[\mu^{+\alpha_i}>\max\;\pcf\{\mu^{+\beta}:\beta\mbox{ a successor, }
\alpha^* \le\beta\le(\alpha^*+1)\cup\bigcup_{j<i}\alpha_j\},\]
and let ${\mathfrak a}_i={\mathfrak b}_{\alpha_i}\setminus\bigcup\limits_{
\ell<m_i}{\mathfrak b}_{\gamma_\ell}$. If $\neg(2^{\aleph_\omega}<\mu)$ use
the end of \cite[\S6]{Sh:430}.]\\
If, weakening $(*)_{\mu,n}$, we assume that for some $\alpha^*<\omega_n$ we
have 
\[\delta>\alpha^*\ \&\ \delta<\omega_n\mbox{ is limit}\quad \Rightarrow\quad
\pp(\mu^{+\delta})<\mu^{+\omega_n},\]
then we can get the same conclusion. Of course, omitting $(*)_{\mu,n}$ if
$2^{\aleph_0}<\omega_n$, by the $\Delta$--system lemma, we can get $\langle
({\mathfrak a}_i,\alpha_i):i<\omega_n\rangle$ as above but demanding only
$i\neq j\ \Rightarrow\ \mu^{+\alpha_i}\notin\pcf({\mathfrak a}_j)$. Of
course, we cannot let $\alpha_i=i+1$, as e.g.~for some infinite
$A\subseteq\omega$, $\mu^{+\omega+1}=\tcf(\prod\limits_{n\in A}\mu^{+n}/
J^{\bd}_A)$, and hence $\mu^{+\omega+1}\in\pcf(\bigcup\limits_{n\in A}
{\mathfrak a}_n)$. 
\medskip

\noindent {\sc Another remark:}\quad Even if $\pcf({\mathfrak a})$ is large
and ${\mathfrak a}$ is countable, we can find a c.c.c.~forcing notion $\bQ$
such that in $\bV^\bQ$ we can find $\langle {\mathfrak b}_\lambda:\lambda\in
\pcf({\mathfrak a})\setminus {\mathfrak a}\rangle$ satisfying:\quad
${\mathfrak b}_\lambda\subseteq{\mathfrak a}$ has order type $\omega$ and
$\prod {\mathfrak b}_\lambda/J^{\bd}_{{\mathfrak b}_\lambda}$ has true
cofinality $\lambda$. [Why?  If $\langle {\mathfrak b}_\lambda:\lambda\in
\pcf({\mathfrak a})\rangle$ is a generating sequence, let $\bQ$ force for
each $\lambda$ an $\omega$--sequence $\subseteq {\mathfrak b}_\lambda$,
almost disjoint to ${\mathfrak b}_{\lambda_1}$ for $\lambda_1<\lambda$.]
Such forcing does not change the pcf structure (in fact, if $\langle
{\mathfrak b}_\lambda:\lambda\in\pcf({\mathfrak a})\rangle$ is a generating
sequence for ${\mathfrak a}$, $\bQ$ is a $\min({\mathfrak a})$--c.c.~forcing
notion, {\em then\/} $\langle {\mathfrak b}_\lambda:\lambda\in\pcf(
{\mathfrak a})\rangle$ is still a generating sequence for ${\mathfrak a}$, 
witnessed by the same $\langle f^*_\alpha:\alpha<\lambda\rangle$).

\begin{question}
\label{1.22}
For a regular cardinal $\theta$, can we find an increasing sequence $\langle
\lambda_i:i<\theta\rangle$ of regular cardinals such that for some successor
$\lambda$ and $f_\alpha\in \prod\limits_{i<\theta}\lambda_i$ for $\alpha<
\lambda$ we have:    
\begin{enumerate}
\item[$(*)$]  if $C_i$ is a club of $\lambda_i$ for $i<\theta$, then for
every large enough $\alpha<\lambda$ for every large enough $i<\theta$ we
have $f_\alpha(i)\in C_i$. 
\end{enumerate}
\end{question}
By \cite[\S6]{Sh:620} an approximation to this holds: if $\mu$ is a strong
limit singular cardinal, $\pp(\mu)=^+ 2^\mu$ and $\lambda=2^\mu=\cf(2^\mu)$
then the answer is yes, i.e.~$(*)$ holds true, but $2^\mu$ may be a limit
cardinal (if $2^\mu$ is singular, a related statement holds). 

\begin{question}
\label{1.23}
Assume $\kappa=\cf(\kappa)$, $\langle\mu_i:i\le\kappa\rangle$ is an
increasing continuous sequence of strong limit cardinals, for nonlimit $i$,
$\cf(\mu_i)=\aleph_0$ and $\prod\limits_{i<\kappa}\mu^{+n}_i/J^{\bd}_\kappa$
has true cofinality $\mu^{+n}$. Can we find an interesting colouring theorem
on $\mu^{+n}$?  (The point is that for $n\ge 2$, we can have both a
colouring as $\mu^{+n}$ is a successor of regulars (as in \cite{Sh:280},
\cite{Sh:572}) and using a witness to $\tcf(\mu^{+n}_i/J^{\bd}_\kappa)=
\mu^{+n}$ as in \cite{Sh:575}, \cite{Sh:620}.) The question is whether
combining we shall get something  startling.
\end{question}

\begin{question}
\label{1.24}
\begin{enumerate}
\item Are there non-metrizable first countable Hausdorff topological spaces
which are $\aleph_2$-metrizable (i.e., the induced topology on any $\le
\aleph_1$ points is metrizable)?
\item Are there non-collectionwise Hausdorff, first countable Hausdorff
topological spaces which are $\aleph_1$--collectionwise Hausdorff?
\end{enumerate}
\end{question}
See \cite{Sh:E9}. Concerning hopes to answer yes note that if SCH fails (or
just $\cf(\mu)=\aleph_0$, $\pp(\mu)>\mu^+$) then there are examples
(see \cite[\S1]{Sh:E9}), so we are allowed to assume $2^{\beth_\omega}=
\beth^+_\omega$, etc.  

\begin{question}
\label{1.25}
Let $D$ be an ultrafilter on $\kappa$ and ${\rm Spc}(D)=\{\prod\limits_{i<
\kappa} \lambda_i/D,\; \lambda_i \ge 2^\kappa$ for $i<\kappa\}$. Is ${\rm
Spc}(D)$ equal to $\{\mu:2^\kappa\le\mu=\mu^{<\text{reg}(D)}\}$?\\
(Where ${\rm reg}(D)=\sup\{\theta:$ for some $A_i\in D$, $i<\theta$, for
every $\alpha<\kappa$, the number of $i<\theta$ such that $\alpha\in A_i$ is
finite $\}$.)
\end{question}
See on this \cite{Sh:589} where some information is gained.

\begin{question}
\label{1.26}
For which $\lambda\geq\mu$ we can find an almost disjoint family ${\mathcal
A}\subseteq [\lambda]^{\aleph_0}$ such that 
\[(\forall X\in [\lambda]^\mu)(\exists A\in {\mathcal A})(A \subseteq^* X)\
?\]
At least when $\lambda\ge\mu=\beth_\omega$?  (See \cite{Sh:460},
\cite{Sh:668}). 
\end{question}

\begin{question}
\label{1.27}
Is it consistent that for some strong limit singular cardinal $\mu$, for no
regular $\lambda\in [\mu,2^\mu]$ do we have a c.c.c.~Boolean Algebra which
is not $\lambda$--Knaster? 
\end{question}
On related ZFC constructions see \cite{Sh:575}, \cite{Sh:620}; see also \S6
here. 

\begin{question}
\label{1.33}
Are all the assumptions in the result of \cite{Sh:580} (see below)
necessary? In particular, are assuptions (a), (b), (c) below sufficient? 
\end{question}

\begin{theorem}
[See \cite{Sh:580}]
\label{1.34}
Assume that
\begin{enumerate}
\item[(a)] $\bV$ is our universe of sets, ${\bf W}$ is another model of ZFC
(i.e., a transitive class of $\bV$ containing all the ordinals),
\item[(b)] $\kappa$ is a regular cardinal in $\bV$,
\item[(c)] $({\bf W},\bV)$ has $\kappa$--convering (that is, every set of
$<\kappa$ ordinals from $\bV$ is included in a set of $<\kappa$ ordinals
from ${\bf W}$),
\item[(d)] the successor of $\kappa$ in $\bV$ is the same as its succesor in
${\bf W}$, call it $\kappa^+$,
\item[(e)]  $({\bf W},\bV)$ has $\kappa^+$--convering.
\end{enumerate}
{\em Then\/} $({\bf W},\bV)$ has the strong $\kappa$--covering (that is,
for every structure $M$ with universe an ordinal $\alpha$ and a countable
vocabulary, and a set $X$ from $\bV$ of cardinality $<\kappa$, there is a
set $Y$ from ${\bf W}$ of cardinality $<\kappa$ including $X$ which is the
universe of an elementary submodel of $M$).
\end{theorem}

\section{The quest for the test: on the theory of Iterated Forcing for the
continuum}
On the subject see \cite{Sh:f}, and recent papers, too, but this section is
hampered by some works in progress. 

The issue is:

\begin{problem}
\label{2.0}
\begin{enumerate}
\item[(a)] Assuming we know something about each iterand $\name{\bQ}_i$,
what can we say about $\bP_\alpha$, where $\langle \bP_i,\name{\bQ}_j:i \le
\alpha,\; j<\alpha\rangle$ is an iteration (which may be FS (finite
support), or CS (countable support) or RSC (revised countable support) and
more) ? 
\item[(b)]  Find more useful ways to iterate (say, new ``supports'').
\end{enumerate}
\end{problem}
So ``c.c.c.~is preserved by FS iteration'', ``properness is preserved by 
CS iteration'' can be seen as prototypes. But also many times: ``adding no 
Cohen real over $\bV$'', ``adding no dominating real over $\bV$'', etc.,
and, very natural, ``adding no new real''.

Note that this is not the same as having forcing axioms, e.g., having (the
very important) MA does not discard the interest in FS iterations of
c.c.c.~forcing. The point is that in many questions you want to add reals
for some purpose (which appear as generic sets for some forcing notions),
but not another (e.g., a well-ordering of $\omega$ of order type
$\omega_1$). Also considering an axiom speaking on forcing notions with some
property, when considering a candidate, a forcing notion $\bP$, during an
iteration we may force that it will not satisfy the property, discard it
instead ``honestly'' forcing with it.

What we get by iterations as above can be phrased as having some axioms, but
we have many combinations of adding reals of kinds A, B, and C while
preserving properties ${\rm Pr}_1$, ${\rm Pr}_2$, in other words
practically one preservation theorem may be used in many such contexts. 

In fact, some of the most intriguing problems are fine distinctions: adding
solution to one kind, but not to a close variant, e.g., the old problem: 

\begin{question}
\label{2.1}
${\rm CON}({\mathfrak p}<{\mathfrak t})$ ?\\
(Note that if ${\mathfrak p}<{\mathfrak t}$ then $2^{\aleph_0}\ge\aleph_3$.
See \ref{3.7}). 
\end{question}
With FS iteration, all values of the continuum were similar, except
$\aleph_1$ (well, also there is a distinction between regular and
singular). 

In fact, the advances in proper forcing make us ``rich in forcing'' for
$2^{\aleph_1}=\aleph_2$, making the higher values more mysterious. (So in
\cite[Ch.VII,VIII]{Sh:f} we separate according to the size of the
$\name{\bQ}_i$'s and whether we add reals, but we concentrate on the length
$\omega_2$). So, because we know much more how to force to get $2^{\aleph_0}
= \aleph_2$, the independence results on the problems of the interrelation
of cardinal invariants of the continuum have, mostly dealt with
relationships of two cardinals, as their values are $\in\{\aleph_1,
2^{\aleph_0}\}$. Thus, having only two possible values $\{\aleph_1,
\aleph_2\}$ among any three, two are equal; the Pigeonhole Principle acts
against us. As we are rich in our knowledge to force for $2^{\aleph_0}=
\aleph_2$, naturally we are quite poor concerning ZFC results. If we try for
cardinal invariants ${\bf c}_1,{\bf c}_2$ to prove they consistently are
$\aleph_1,\aleph_2$, respectively, much of our way exists (quoting existing
preservation theorems) and we can look at the peculiarities of those
invariants which may be still intractable. We are not poor concerning
forcing for $2^{\aleph_0}=\aleph_1$ (and are rich in ZFC). But for
$2^{\aleph_0}\ge \aleph_3$ we are totally lost: very poor in both
directions. We would like to have iteration theory for length
$\ge\omega_3$. I tend to think good test problems will be important in
developing such iterations.   

In some senses, most suitable is 

\begin{problem}
\label{2.2}
Investigate cardinal invariants of the continuum showing $\ge 3$ may have
prescribed order. 
\end{problem}
Of course, the lack of forcing ability does not stop you from proving
hopeful ZFC theorems about them, if true. Now I think there are some, but: 

\begin{thesis}
\label{2.3}
They are camouflaged by the independent statements.
\end{thesis}

\noindent [Yes, I really believe there are interesting restrictions.]
However, once we prove 90 percent of the problems are independent we will
know where to look (as in hindsight occurs in cardinal arithmetic). So
cardinal invariants from this perspective are excellent excuses to find
iteration theorems. Mainly for $2^{\aleph_0}\ge\aleph_3$, but, of course,
there is more to be said on $2^{\aleph_0}=\aleph_1$ (though not for
\ref{2.2}), and even $2^{\aleph_0} =\aleph_2$. 

Without good test problems you are in danger of imitating the king who
painted the target after shooting the arrow.  Let us consider some additional
well known problems:

\begin{question}
[See Just, Mathias, Prikry and Simon \cite{JMPS}]
\label{2.4}
Is there a filter $D$ on $\omega$ such that:
\begin{enumerate}
\item[(a)] every co-finite subset of $\omega$ belongs to $D$,
\item[(b)] $D$ is a $P$--filter (i.e., if $A_n\in D$ for $n<\omega$, then
for some $A\in D$, $n<\omega\ \Rightarrow\ A\subseteq^* A_n$),
\item[(c)] $D$ is not feeble, i.e., if $0=n_0<n_1<\ldots$, then for some
$A\in D$ for infinitely many $i<\omega$ we have $[n_i,n_{i+1})\cap A=
\emptyset$. 
\end{enumerate}
\end{question}

\begin{question}
[See Garcia--Ferreira and Just \cite{GFJ97}]
\label{2.5}
Is there an almost disjoint family ${\mathcal A}\subseteq
[\omega]^{\aleph_0}$ (i.e., $(\forall A \ne B\in {\mathcal A})[|A \cap B|<
\aleph_0]$) of cardinality ${\mathfrak b}$ satisfying the following
condition:  
\begin{quotation}
if $A_n\in {\mathcal A}$ are pairwise distinct and $h:\omega
\longrightarrow\omega$\\  
{\em then\/} for some $B\in {\mathcal A}$ we have $(\exists^\infty n)(A_n
\cap B\nsubseteq h(n))$ ?
\end{quotation}
\end{question}
If not, then $2^{\aleph_0}>\aleph_\omega$; on both questions see the
discussions after \ref{2.12}. 

\begin{question}
[See van Mill {\cite[Problem 4, p.563]{vM84}}, Miller {\cite[Problem
9.1]{Mi91}}] 
\label{2.6}
\[{\rm CON}(\mbox{no $P$-point and no $Q$-point})\ ?\]
\end{question}
If so, $2^{\aleph_0}\ge\aleph_3$. [Why? Mathias \cite{Mt78} showed that if
${\mathfrak d}$ (the minimal size of a dominating family) is $\aleph_1$,
then there is a $Q$--point. Ketonen \cite{Kt76} showed that ${\mathfrak
d}=2^{\aleph_0}$ implies the existence of $P$--points.] 

\begin{question}
\label{2.7} 
{\rm CON}$({}^\omega(\omega +1)$ with box product topology is not
paracompact)? 
\end{question}
If so, $2^{\aleph_0}\ge\aleph_3$.  See on this Williams \cite{Wl84}.

\begin{question}
[See Miller {\cite[Problem 16.3]{Mi91}}] 
\label{2.8} 
\[{\rm CON}(\mbox{Borel Conjecture and Dual Borel Conjecture})\ ?\]
(See \ref{3.3}).
\end{question}

\begin{question}
\label{2.9}
{\rm CON}($\cf(\cov({\rm meagre}))<{\rm additivity}({\rm meagre})$)?\\
(See before \ref{2.15}).
\end{question}

\begin{problem}
\label{2.10}
\begin{enumerate}
\item {\rm CON}(every function $f:{}^\omega 2\longrightarrow {}^\omega 2$ is
continuous when restricted to some non-null set)?\\
\relax [Here ``null'' means of Lebesgue measure zero.]
\item Similarly for other natural ideals. This in particular means if $\bQ$
is a nicely defined forcing notion (see \S5 below, e.g., Souslin c.c.c.),
$\name{\eta}$ a $\bQ$-name of a real, $A\subseteq {}^\omega 2$ is called
$(\bQ,\name{\eta})$--positive if for every countable $N\prec ({\mathcal
H}(\chi),\in,<^*_\chi)$ to which $\bQ$, $\name{\eta}$ belong, some $\eta\in
A$ is $\name{\eta}[G]^N$ for some $G\subseteq \bQ^N$ generic over $N$; so
the question for such $\bQ$ is ``{\rm CON}(every $f:{}^\omega
2\longrightarrow {}^\omega 2$ has a continuous restriction to some
$(\bQ,\name{\eta})$--positive set $A$)?''  
\item Is the following consistent:\\
if $A\subseteq {}^\omega 2$ is non-null, $f:A\longrightarrow {}^\omega 2$
{\em then\/} for some positive $B\subseteq A$, $f\restriction B$ is
continuous.  Similarly for general ideals as in part (2).
\item If $A\subseteq {}^\omega 2\times {}^\omega 2$ is not equivalent to a
Borel set modulo one ideal $I_1$ (as described in part (2) above), {\em
then} for some continuous $f:{}^\omega 2\longrightarrow {}^\omega 2$, the
set $\{\eta\in {}^\omega 2:(\eta,f(\eta))\in A\}$ is not equivalent to a
Borel set modulo another ideal $I_0$ for suitable pairs $(I_0,I_1)$. 
\end{enumerate}
\noindent See Fremlin \cite{Fe94}, Ciesielski \cite[Theorem 3.13, Problem
5]{Ci97}; \cite{Sh:473} shows ``yes'' for (2) for non meagre, Ciesielski and
Shelah \cite{CiSh:695} prove ``yes'' for (4) for non meagre, on work in
progress see Ros{\l}anowski and Shelah \cite[\S 2]{RoSh:670}. With Juris
Stepr\=ans we have had some discussions on trying to use the oracle cc to
the case of non-megre ideal in (2). See \ref{3.8}.
\end{problem}

\noindent{\sc Note:}\quad Mathematicians who are not set theorists
generally  consider ``null'' as senior to ``meagre'', that is as a more
important case; set theorists inversely, as set-theoretically Cohen reals
are much more manageable than random reals and have generalizations,
relatives, etc. Particularly, in FS iterations, we get Cohen reals ``for
free'' (in the limit), which kills our chances for many things and until now 
we have nothing parallel for random reals (but see \cite{RoSh:670}). 

Judah suggests: 

\begin{question}
[$\bV=\bL$]  
\label{2.11} 
Find a forcing making ${\mathfrak d}=\aleph_3$ but not adding Cohen reals.
\end{question}
I am skeptical whether this is a good test question, as you may make
${\mathfrak d}=\aleph_3={\mathfrak b}$ by c.c.c.~forcing, then add
$\aleph_1$ random reals $\langle\nu_i:i<\omega_1\rangle$ by a measure
algebra; so over $\bL[\langle\nu_i:i<\omega_1\rangle]$ we have such a
forcing. But certainly ``not adding Cohen'' is important, as many problems
are resolved if $\cov({\rm meagre})=2^{\aleph_0}$.

There is a basic question for us:

\begin{problem}
\label{2.12}
Is there an iteration theorem solving all the problems described above or at
least for all cases involving large continuum not adding Cohen reals? 
\end{problem}
I suspect not, and the answers will be ramified.

Let us review some problems. Now, Problems \ref{2.4}, \ref{2.5} are for
$2^{\aleph_0}>\aleph_\omega$, as: in \ref{2.4}, if $\cf([{\mathfrak d}]^{
\aleph_0},\subseteq)={\mathfrak d}$ then there is such a filter (see
\cite{JMPS}), and also in \ref{2.5}, if $\cf([{\mathfrak b}]^{\aleph_0},
\subseteq)={\mathfrak b}$ there is a solution (see Just, Mathias, Prikry and
Simon \cite{JMPS}).

It may well be that the solution will look like: let $\mu$ be a strong limit
singular cardinal with $2^\mu\ge\mu^{++}>\mu^+$ and we use FS iteration of
length $\mu^{++}$. This will be great, but probably does not increase our
knowledge of iterations. If on the other hand along the way we will add new
$\omega$-sequences say to $\mu$ (say $\cf(\mu)=\aleph_0$) and necessarily we
use more complicated iteration, then it will involve better understanding of
iterations, probably new ones.

We can (\cite[Ch.XIV]{Sh:f}) iterate up to ``large'' $\kappa$, and for many
$\alpha<\kappa$, $\alpha$ strongly inaccessible, we have $\bQ_\alpha$ change
its cofinality to $\aleph_0$.  Sounds nice, but no target yet.

\[* \qquad * \qquad *\]

We may note that ``FS iterations of c.c.c.~forcing notions'' is not
dead. Concerning \ref{2.2} and \ref{2.9}, there are recent indications
that FS iteration of c.c.c.~still can be exploited even in cases for which
for a long time we thought new supports are needed. We can iterate with FS,
$\langle \bP_i, \name{\bQ}_i:i<\alpha\rangle$, where $\name{\bQ}_i$ is
(partially random or is Cohen) adding a generic real $r_i$, $\name{\bQ}_i$
is Cohen forcing or random forcing over $\bV[\langle r_j:j\in A_i\rangle]$,
where $A_i\subseteq i$ and each $\name{\bQ}_i$ is reasonably understood; but
we do not require $j\in A_i\ \Rightarrow\ A_j\subseteq A_i$ (so called
transitive memory). It is not so immediate to understand this sort of
iterations, e.g., can the iteration add a dominating real?   

It appears that if the $A_i$'s are sufficiently closed, it will not, see
\cite{Sh:592}, more generally look at \cite{Sh:619}. There we prove: 
\[{\rm CON}(\exists\mbox{ non-null $A$ such that the null ideal restricted
to $A$ is $\aleph_1$-saturated}).\]
Clearly we should use a measurable cardinal $\kappa$, a normal ultrafilter
${\mathcal D}$ on $\kappa$ in $\bV$ and we add $\kappa$ random reals
$\langle r_\zeta:\zeta<\kappa\rangle$, but how do we make 
\[A\in {\mathcal D}^\bV\quad \Rightarrow\quad \{r_\zeta:\zeta\in \kappa
\setminus A\}\mbox{ is null ?}\]
Cohen forcing does the job, but unfortunately too strongly. (In the case
with non meagre, by Komjath \cite{Ko}, there are no problems in this
respect). The solution is that we use FS iteration, but first we add
$2^\kappa$ Cohens, that is $\bQ_i$ is Cohen forcing for $i<2^\kappa$; only,
then we add the (somewhat) random reals: 
\[\name{\bQ}_{\lambda+\zeta}\mbox{ for }\zeta<\kappa\mbox{ is random forcing
over }\langle r_{\lambda+\xi}:\xi<\zeta\rangle,\ \langle r_i:i\in A_\zeta
\subseteq\lambda\rangle.\]
We use $\{r_i:i<\lambda\}$ such that: $r_i$ makes $\{r_{\lambda+\zeta}:
\zeta \mbox{ satisfies } i\in A_\zeta\}$ null for $i<2^{\aleph_0}$. So we
need: for $A\in {\mathcal D}^\bV$ for some $i<\lambda$ for every
$\xi<\kappa$, $\xi\in\kappa\setminus A\ \Leftrightarrow\ i\notin A_\xi$.

This works for specially chosen $A_\zeta$'s.

\begin{problem}
\label{2.15}
\begin{enumerate}
\item Can you make this into a general method?
\item Can you deal with $n$ or even $\kappa$ kinds of ``reals'' (getting 
interesting results)?
\end{enumerate}
\end{problem}
\noindent What does this mean? This means that we use FS iteration $\langle
\bP_i,\name{\bQ}_j:i\le\delta,\ j<\delta\rangle$ and $h:\alpha
\longrightarrow\beta$, for $\zeta<\beta$, $\bR_\zeta$ is a nep
c.c.c.~forcing notion in $\bV$, (on nep see \S5 and more see \cite{Sh:630};
e.g., $\bR_\zeta$ is Cohen, random, or as in \cite{RoSh:628},
\cite{RoSh:672} or whatever), and $\name{s}_\zeta\in {}^\omega 2$ is a
generic real for $\bR_\zeta$, and $\name{\bQ}_i$ is $\bR_{h(\zeta)}$ as
interpreted in $\bV[\langle \name{r}_j:j\in A_i\rangle]$ and $\name{r}_j$ is
$\name{s}_{h(j)}$ there, and $A_i\subseteq i$. So the idea is that $0=
\delta_0 <\delta_i<\ldots<\delta_n=\delta$ and $j\in [\delta_\ell,
\delta_{\ell +1})\ \Rightarrow\ h(j)=\ell$. 

\[* \qquad * \qquad *\]

In \cite{Sh:538} we use $\aleph_\varepsilon$--support. This is less than
$(<\aleph_1)$--support (i.e., countable support). This looks quite special,
but 
\begin{problem}
\label{2.16}
Can we make a general (interesting) theorem?  
\end{problem}
We can note that long FS iterations not only add Cohen reals, they also add,
e.g., $\aleph_2$--Cohens, i.e.~generics for $\{f:f$ a finite function from
$\omega_2$ to $\{0,1\}\;\}$. So we may like to iterate, allowing to add
Cohen reals but not $\aleph_2$--Cohens in the sense above. This is done in
\cite{Sh:538}, but the family of allowable iterands can be probably widened.

\[* \qquad * \qquad *\]

If we agree that preservation theorems are worthwhile, then after not
collapsing $\aleph_1$, probably the most natural case is adding no reals.
Now, whereas properness seems to me both naturally clear and covers
considerable ground for not collapsing $\aleph_1$ and there are reasonable
preservation theorems for ``proper $+ X$'' for many natural properties $X$
(e.g., adding no dominating reals, see \cite[Ch.VI,\S1,\S2,\S3]{Sh:f},
\cite[Ch.XVIII,\S2]{Sh:f}), the situation with NNR (no new real) is
inherently more complicated. In the early seventies when I heard on Jensen's
${\rm CON}({\rm GCH} + {\rm SH})$, I thought it would be easy to derive an
axiom; some years later this materialized as Abraham, Devlin and Shelah
\cite{ADSh:81}, but reality is not as nice as dreams. One obstacle is the
weak diamond, see Devlin and Shelah \cite{DvSh:65}, more in
\cite[Ch.XIV,\S1]{Sh:b}, \cite[AP,\S1]{Sh:f}, \cite{Sh:638}. For a time the
iteration theorem in \cite[Ch.V,\S5,\S7,Ch.VIII,\S4]{Sh:b} seemed
satisfactory to me. [There we use two demands. The first was ${\mathbb
D}$--completeness (this is a ``medicine'' against the weak diamond, and
${\mathbb D}$ is a completeness system, $\aleph_1$--complete in 
\cite[Ch V, \S 5,7]{Sh:b}, that is any countably many demands are compatible,
and 2--complete in \cite[Ch VIII, \S 4]{Sh:b}, that is any two demands are
compatible). The second demand was $\alpha$--properness for each countable
ordinal $\alpha$ (or relativized version, see \cite[Ch VIII, \S 4]{Sh:b},
\cite[Ch VIII, \S 4]{Sh:f}).] But \cite[\S1]{Sh:177} (better \cite[Ch.XVIII,
\S1]{Sh:f}) gives on the one hand very nice and easy forcing notions not
adding reals (running away from club guessing sequences) which are not
covered as they fail $(<\omega_1)$--properness and on the other hand, shows
by a not so nice example that generally you cannot just omit the
$(<\omega_1)$--properness demand and promise an iteration theorem covering
them. The problem concerning that forcing was resolved (promised in
\cite{Sh:177}, carried out in a different way in lectures in MSRI '89 =
\cite[Ch.XVIII,\S2]{Sh:f}), but resulted in a dichotomy: we can get by
forcing  ${\rm CON}({\rm ZFC}+{\rm CH}+{\rm SH})$ and we can get by forcing
${\rm CON}({\rm ZFC}+{\rm CH}+\mbox{ no club guessing})$, but can we have
both? More generally, can we have two other such contradictory statements
(more generally for such results see Shelah and Zapletal \cite{ShZa:610}).

\begin{question}
\label{2.17}
Can we have two statements of the form\footnote{${\mathcal H}(\lambda)$ is
the family of sets with transitive closure of cardinality $<\lambda$}
\[(\forall x\in {\mathcal H}(\aleph_2))(\exists y\in {\mathcal H}(\aleph_2))
\varphi\]
each consistent with 
\[\begin{array}{ll}
{\rm CH} + [{\rm Axiom}(\bQ\mbox{ is}&(<\omega_1)\mbox{--proper and
${\mathbb D}$--complete}\\
&\mbox{for some simple $2$-completeness system ${\mathbb D}$})]
  \end{array}\]
but not simultaneously?\\  
(We may change the axiom used, we may speak directly on the iteration; we
may deal with CS and proper or with RCS and semi-proper, etc.) 
\end{question}
Note: possible failure of iteration does not prove a ZFC consequence, we may
have freedom in the iteration only in some stages (like c.c.c.~productive
under MA). 

This leaves me in bad shape: the iteration theorems seem not good enough,
but the test problem (of getting both) does not seem so good. Now,
\cite{Sh:656} deals with NNR solving the specific dichotomy (and really
satisfies the \cite{Sh:177} promise circumvented in \cite[Ch.XVIII,
\S2]{Sh:f}) but left \ref{2.16} open.

Eisworth suggested to me (motived by Abraham and Todor\v{c}evi\'{c}
\cite{AbTo97}) 

\begin{question}
\label{2.17A}
Is the following consistent with ZFC+CH
\begin{enumerate}
\item[$(*)$]  if $A_\alpha\in [\omega_1]^{\aleph_0}$ and $\alpha<\beta\
\Rightarrow\ A_\alpha\subseteq A_\beta$ mod finite, and for every stationary
$S\subseteq\omega_1$ the set $\bigcup\{[A_\alpha]^{<\aleph_0}:\alpha\in S\}$
contains $[E]^{< \aleph_0}$ for some club $E$ of $\omega_1$,\\
{\em then\/} for some club $C$ of $\omega_1$ we have
\[(\forall\alpha<\omega_1)(\exists\beta<\omega_1)(C\cap\alpha\subseteq
A_\beta).\]
\end{enumerate}
\end{question}

For long, an exciting problem for me has been

\begin{problem}
\label{2.18}
\begin{enumerate}
\item Can we find a consequence of ZFC + CH which ``stands behind'' the
``club objection to NNR", e.g.~it implies the failure of CH + Axiom($\bQ$
proper $\mathbb D$--complete for some single 2-completeness system) ?
\item Similarly for other limitations on iteration theorems? 
\end{enumerate}
\end{problem}

\begin{question}
\label{2.18A}
Is ``${\rm CH}+{\mathcal D}_{\omega_1}$ is $\aleph_2$--saturated''
consistent, where ${\mathcal D}_{\omega_1}$ is the club filter on
$\omega_1$?\\ 
Recall that a filter $D$ on a set $A$ is $\lambda$--saturated if there are
no $A_i\in D^+$ for $i<\lambda$ such that $i<j\ \Rightarrow\ A_i\cap A_j= 
\emptyset$ mod $D$. 
\end{question}
\noindent See \cite[Ch.XVI]{Sh:f}. Woodin proved that if there is a
measurable cardinal then no. So we may look at $L[A]$, $A\subseteq\kappa$
codes ${\mathcal H}(\chi)$, $\kappa$ large and try to collapse it to
$\omega_2$. 

Note that by \cite{Sh:638}, if ${\rm CH} +{\mathcal D}_{\omega_1}$ is
$\aleph_2$--saturated, then essentially we have the weak diamond for three
colours (or any finite number). 

\[* \qquad * \qquad *\]

Baumgartner \cite{B4} asked 

\begin{question}
\label{2.19}
Is it consistent that $2^{\aleph_0}>\aleph_2$ and any two $\aleph_2$--dense
subsets of $\bR$ of cardinality $\aleph_2$ (that is, any interval has
$\aleph_2$ points) are isomorphic (as linear orders). 
\end{question}

I think it is more reasonable to try 

\begin{question}
\label{2.20}
Is it consistent that:\quad $2^{\aleph_0}>\lambda\ge\aleph_2$ and there are
no two far subsets $A\in [\bR]^\lambda$, where
\end{question}

\begin{definition}
\label{2.21}
The (linear orders) $I,J$ are $\theta$--far if there is no linear order of
cardinality $\theta$ embedded into both.  If $\theta$ is omitted, we mean
$\min\{|I|,|J|\}$.  
\end{definition}

On OCA$'$ (i.e., OCA$_{\aleph_1,\omega}$, see the definition below) see
Abraham, Rubin and Shelah \cite{ARSh:153}, continued for OCA$''$ in
Todor\v{c}evi\'{c} \cite{To89}, Veli\v ckovi\'c \cite{Ve92}; on a parallel
for subsets of the plain which follows from MA, see Stepr\=ans and Watson
\cite {SW87}. 

\begin{question}
\label{2.22}
\begin{enumerate}
\item Is OCA$'_{\aleph_2}$ consistent?  Is OCA$''_{\aleph_2}$ consistent?
\item The parallel problems for ${}^\kappa 2$ and $\lambda$, even for
$\lambda=\kappa^+$, $\kappa>\aleph_0$, where
\end{enumerate}
\end{question}

\begin{definition}
\label{2.23}
\begin{enumerate}
\item OCA$'_{\lambda,\kappa}$ means $\lambda\le 2^\kappa$ and: for any
$A\in [{}^\kappa 2]^\lambda$ and an open symmetric set ${\mathcal
U}\subseteq {}^\kappa 2\times {}^\kappa 2$ there is $B\subseteq A$ of
cardinality $\lambda$ such that $\{(a,b):a\ne b\mbox{ are from } B\}$ is
included in ${\mathcal U}$ or is disjoint to  ${\mathcal U}$ (we use the
space ${}^\kappa 2$ for simplicity). 
\item OCA$''_{\lambda,\kappa}$ is defined similarly only we have $B_i
\subseteq A$ for $i<\kappa$, $A=\bigcup\limits_{i<\kappa} B_i$, each $B_i$
as in part (1). 
\item If we omit $\lambda$ we mean $\lambda=\kappa^+$, if in addition we
omit $\kappa$, we mean $\kappa = \aleph_0$. 
\end{enumerate}
\end{definition}

\section{Case studies for iterated forcing for the reals}
The following was suggested during the lecture on \S2 by Juh\'{a}sz who was in
the audience: 

\begin{question}
\label{3.1}
Does CH imply that there is an $S$--space of cardinality $\aleph_2$, where
$S$-space is defined as being regular, hereditarily separable, not
Lindel\"of?

Eisworth prefers the variant:\\ 
Does CH imply the existence of a locally compact $S$--space? 
\end{question}
This problem looks important, but it is not clear to me if it is relevant to
developing iteration theorems, though an existence proof may be related to
the weak diamond, consistency to NNR iterations.

The same goes for  the well known: 

\begin{question}
\label{3.2}
CON$({\mathfrak d}<{\mathfrak a})$ ?
\end{question}
This definitely seems not to be connected to the iteration problem. It seems
to me that a good test problem for our purpose in \S2 should have one step
clear but the iteration problematic, whereas for those two problems the
situation is the inverse.
 
Note: by existing iteration theorems to get the consistency of ${\mathfrak
d}<{\mathfrak a} + 2^{\aleph_0} = \aleph_2$ it is enough to show 
\begin{enumerate}
\item[$(*)$] for any MAD family $\{A_i:i<i^*\}\subseteq[\omega]^{\aleph_0}$,
there is an ${}^\omega\omega$--bounding proper forcing notion $\bQ$ of
cardinality $\aleph_1$ adding $\name{A}\in [\omega]^{\aleph_0}$ almost
disjoint to each $A_i$.  
\end{enumerate}
You are allowed to assume CH (start with $\bV\models 2^{\aleph_0}=\aleph_1 +
2^{\aleph_1} = \aleph_2$ and use CS iteration of such forcing notion); even
$\diamondsuit_{\aleph_1}$ (if $\bV\models\diamondsuit_{\aleph_1}$). We can
weaken $|\bQ|=\aleph_1$ to ``$\bQ$ satisfies $\aleph_2$--pic'' (this is a
strong form of $\aleph_2$--c.c.~good for iterating proper forcing, see
\cite[Ch.VIII,2.1,p.409]{Sh:f}). If you agree to use large cardinals, it is
okay to assume in $(*)$ that an appropriate forcing axiom holds and not
restrict $|\bQ|$, and as we can first collapse $2^{\aleph_0}$ to $\aleph_1$,
we can get $\diamondsuit_{\aleph_1}$ for ``free''. I idly thought to use
free forcing for the problem, (\cite[Ch.IX]{Sh:f}), but no illumination
resulted.   

We can try in another way: start with a universe with a forcing axiom (say
MA) and force by some $\bP$, which makes ${\mathfrak d}=\aleph_1$, but $\bP$
is understood well enough and we can show that ${\mathfrak a}$ is still
large (just as adding a Cohen real to a model of MA preserves some
consequences of MA (see Roitman \cite{Rt79}, Judah and Shelah
\cite{JdSh:335}]). So clearly FS iteration will not do.

I think that a more interesting way is to consider, assuming CH,
\[\begin{array}{ll}
K_{\omega_1} =\{(\bar{\bP},\name{\bar{r}}):&\bar{\bP}=\langle\bP_i:i<
\omega_1\rangle\mbox{ is $\lesdot$--increasing},\ |\bP_i|\le\aleph_1,\\
&\name{\bar{r}}=\langle\name{r}_i:i<\omega_1\rangle,\ \name{r}_i\mbox{
is a $\bP_{i+1}$--name},\\ 
&\Vdash_{\bP_{i+1}}\mbox{`` }\name{r}_i\in {}^\omega\omega\mbox{ dominates }
({}^\omega \omega)^{\bV^{\bP_i}]}\mbox{ ''}\}
\end{array}\]
ordered naturally, and for a generic enough $\omega_2$--limit $\langle
(\bar{\bP}^\zeta,\name{\bar{r}}^\zeta):\zeta<\omega_2\rangle$ we may use
$\bigcup\limits_{\scriptstyle i<\omega_1,\atop\scriptstyle\zeta<\omega_2}
\bP^\zeta_i$. Another way to try is the non-Cohen Oracle \cite{Sh:669}. The
difference is small. Also the ``$\omega_2+\omega_1$--length mix
finite/countable pure support iteration'' seems similar. 

I have just heard about CON$({\mathfrak u}<{\mathfrak a})$ being an old
problem, clearly related to CON$({\mathfrak d}<{\mathfrak a})$. I do not see
much difference at present.

Another direction is to develop the historic $\aleph_\varepsilon$--support
iteration from \cite{Sh:538}.
\medskip

\begin{discussion}
\label{3.3}
Concerning \ref{2.8}, I had not really considered it (except when Judah
spoke to me about it) but just before the lecture, Bartoszy\'nski reminded
me of it (see \cite{BaJu95}). Now, ``the'' proof of CON(Borel conjecture) is
by CS iteration of Laver forcing (see Laver \cite{L1}), whereas the
consistency proof of the dual is adding many Cohen reals (see Carlson
\cite{Ca93}). So in a (hopeful) iteration proving consistency we have two
kinds of assignments. We are given, say in stage $\alpha$, in $\bV^{
\bP_\alpha}$ a set $A=\{\eta_i:i<\omega_1\}\subseteq {}^\omega 2$, and we
should make it not of strong measure zero, so we should add an increasing
sequence $\bar{n}=\langle n_\ell:\ell<\omega\rangle$ of natural numbers such
that for no $\bar{\nu}=\langle\nu_\ell:\ell<\omega\rangle\in\prod\limits_{
\ell<\omega}{}^{(n_\ell)} 2$ do we have $(\forall i<\omega_1)(\exists^\infty
\ell)(\nu_\ell\vartriangleleft\eta_i)$. Now, even if we define $\bQ_\alpha$
to add such $\bar{n}$, we have to preserve it later, so it is easier to
preserve, for some family $F\subseteq\prod\limits_{\ell<\omega} 2^{n_\ell}$,
the demand  
\[(\forall f\in F)\neg(\exists\langle\nu_{\ell,k}:\ell<\omega,k<f(\ell)
\rangle)(\forall i<\omega_1)(\exists^\infty \ell)(\exists k<f(\ell))(
\nu_{\ell,k} \vartriangleleft \eta_i).\]
The second kind of assignment which we have in stage $\alpha$ is the
following. In $\bV^{\bP_\alpha}$, we are given $A=\{\eta_i:i<\omega_1\}
\subseteq {}^\omega 2$ and we should make it non-strongly meagre, so we
should add, by $\bQ_\alpha$, a subtree $T_\alpha\subseteq {}^{\omega>}2$
(i.e., $\langle\rangle\in T_\alpha$, $\eta\in T_\alpha\ \&\ \nu
\vartriangleleft\eta\ \Rightarrow\ \nu\in T_\alpha$, $\eta\in T_\alpha\
\Rightarrow\ (\exists \ell<2)(\eta\conc\langle\ell\rangle \in T_\alpha))$ of
positive measure (i.e., $0<\inf\{|T_\alpha\cap {}^n 2|/2^n:n<\omega\}$) such
that $(\forall\eta\in {}^\omega 2)(\exists i<\omega_1)[\eta\oplus\eta_i
\notin\bigcup\limits_n(T_\alpha)^{[n]}]$, 
where
\[\begin{array}{ll}
(T_\alpha)^{[n]}=\{\nu:&\mbox{ for some }\rho\in T_\alpha\mbox{ we have }
  \ell g(\nu)=\ell g(\rho) \\ 
&\mbox{ and } (\forall\ell)(n\le \ell<\ell g(\nu)\ \Rightarrow\ \rho(\ell)=
  \nu(\ell))\}. 
  \end{array}\]
Again we have to preserve this.

A way to deal with such preservation problems is to generalize ``oracle
c.c.c.'' (see \cite[Ch.IV]{Sh:f}) replacing Cohen by other things. To
explain this, it seems reasonable to look at the ``oracle for random'' (or
even sequence of c.c.c.~Souslin forcing, from \cite{Sh:669}). This evolves
to: for iterations of length $\le\omega_2$ of forcing notions of cardinality
$\aleph_1$, prove that we can preserve the following condition on $\bP=
\bP_\alpha$ for some $\langle M_\delta,M^+_\delta,r_\delta:\delta\in S
\rangle$, $S\subseteq\omega_1$ stationary such that $\langle M_\delta:\delta
\in S\rangle$ is an oracle, i.e., a $\diamondsuit^*$--sequence and $M_\delta
\models\delta=\omega_1$, $M_\delta\models {\rm ZFC}^-_*$,
$M^+_\delta\models$ ``${\rm ZFC}^-_* + M_\delta$ is countable'' and
$r_\delta$ is random over $M^+_\delta$. Now  without loss of generality,
$\bP \subseteq\omega_1$ and 
\[\begin{array}{lrl}
\{\delta\in S:&\bP\restriction\delta\in M_\delta,\mbox{ and for every }p\in
\bP\cap\delta,\mbox{ for some }q\mbox{ we have }\ &\\
&p\le q\in\bP\mbox{ and }q\Vdash\mbox{`` }r_\delta\mbox{ is random over
}M^+_\delta[\name{G}_{\bP}\cap\delta]\mbox{ '' }\}&\in {\mathcal
D}_{\omega_1}\restriction S
  \end{array}\]
(so this is like the oracle c.c.c.~(\cite[Ch.IV]{Sh:f}), but the support is
not countable so on other stationary $S_1\subseteq\omega_1\setminus S$ we
may have different behaviour). Of course, we use ``small'' $S$ so that we
have ``space'' for more demands, see \cite{Sh:669}. But trying to explain it
(to Ros{\l}anowski) it seemed the proof is too simple, so we can go back to
good old CS and just preserving an appropriate property, a watered-down
relative in the nep family (\cite{Sh:630}).
\end{discussion}

We mainly try to combine the two iterations (of Cohen and of Laver forcing
notions): 

\begin{definition}
\label{3.4}
A forcing notion $\bQ$ is 1--e.l.c.~if the following condition is satisfied:
\begin{quotation}
whenever $\chi$ is large enough, $M_0\prec M_1\prec({\mathcal H}(\chi),
\in)$, $\bQ\in M_0$, $M_0\in M_1$ and $M_0,M_1$ are countable and $p\in\bQ
\cap M_0$,\\
{\em then\/} for some condition $q\in\bQ$ stronger than $p$ we have
\[q \Vdash\mbox{``}\mbox{for every }{\mathcal I}\in M_1\mbox{ such that }
{\mathcal I}\cap M_0\mbox{ is predense in }\bQ^{M_0}\mbox{ we have }
\name{G}_\bQ \cap {\mathcal I}\ne\emptyset\mbox{''.}\]
\end{quotation}
(Note that $q\Vdash$`` $M_0[\name{G}_\bQ\cap M_0]$ is a generic extension of
$M_0$ for a forcing notion which $M_1$ thinks is countable ''.)
\end{definition}
Note: e.l.c.~stands for {\em elementary locally Cohen}. This is, of course,
close to Cohen, or more accurately is another way to present strongly
proper. But we also seem to need Laver forcing (or a close relative of it),
but it is far from being strongly proper.  Still it satisfies the parallel
if we demand ``${\mathcal I} \subseteq Q^{M_1}$ is predense under pure
extensions'', i.e., with the same trunk. This approach seems to me promising
but it is not clear what it delivers. 

We may consider a more general definition (and natural preservation):

\begin{definition}
\label{3.5}
Let ${\rm Pr}$ be a property. A forcing notion $\bQ$ with generic $\name{X}
\subseteq\alpha_\bQ$ (i.e.~$\bV[\name{G}_\bQ]=\bV[\name{X}[\name{G}_\bQ]]$,
$\alpha_\bQ$ an ordinal) is called e.l.--${\rm Pr}$ forcing if:
\begin{quotation}
for $\chi$ large enough, if $\bQ,\name{X}\in M_1\prec M_2\prec({\mathcal
H}(\chi),\in)$, $M_1,M_2$ countable, $M_1\in M_2$, $p\in\bQ\cap M_1$,\\
{\em then\/} we can find $q,\bQ'$ such that 
\begin{enumerate}
\item[(a)]  $p\le q\in\bQ$,
\item[(b)]  $\bQ'\in M_2$ is a forcing notion with $\name{X}'\subseteq
\alpha_\bQ$ generic, 
\item[(c)]  $M_2\models {\rm Pr}(\bQ',M_1,p)$,
\item[(d)]  $q\Vdash$`` $\name{X}\restriction M_2$ is a $\bQ'$-generic over
$p$ and for some set $G'\subseteq(\bQ')^{M_2}$, generic over $M_2$ we have
$\name{X}\restriction M_2=\name{X}'[G']$ ''.
\end{enumerate}
\end{quotation}
\end{definition}
This seems to me interesting but though Laver forcing satisfies some
relatives of those properties it does not seem to be enough.

Note: this definition tells us that generically for many countable models $M
\prec ({\mathcal H}(\chi),\in)$, we have some $q\in\name{G}_\bQ$ which is
almost $(M,\bQ)$--generic, but not quite. The ``almost'' is because this
holds for another forcing $\bQ'$. So when the whole universe is extended
generically for $\bQ$, $M$ ``fakes'' and is instead extended generically for
$\bQ'$. So for preservation in iteration it is not natural to demand $M_2
\prec ({\mathcal H}(\chi),\in)$, but rather to proceed as in \cite{Sh:630},
this will be n.e.l.--${\rm Pr}$.

We may wonder (considering \ref{2.8}) whether we can replace Laver forcing
in the proof of the consistency of the Borel conjecture, by a forcing notion
not adding a dominating real.  So a related question to \ref{2.8} is 

\begin{question}
\label{3.6}
\[{\rm CON}({\mathfrak b}= \aleph_1 +\mbox{ Borel Conjecture})\ ?\]
\end{question}
It is most natural to iterate, one basic step will be $\bQ$, adding an
increasing sequence $\langle \name{n}_i:i<\omega\rangle$ such that on the
one hand: 
\begin{enumerate}
\item[(a)]  no old non-dominated family $\subseteq{}^\omega\omega$ is
dominated (or at least some particular old family remains undominated), 
\end{enumerate}
while on the other hand
\begin{enumerate}
\item[(b)]  for any uncountable $A\subseteq {}^\omega 2$, from $\bV$, we
have: 
\[\Vdash_\bQ\mbox{`` for no $\name{\eta}_i\in {}^{n_i}2$, ($i<\omega$) do
we have $(\forall\nu\in A)(\exists^\infty i)(\name{\eta}_i\vartriangleleft
\nu)$ '',}\] 
\end{enumerate}
or at least
\begin{enumerate}
\item[(b)$'$] like (b) for one $A$ given by bookkeeping.
\end{enumerate}
(To preserve we need to strengthen the statement, replacing $\langle
\name{\eta}_i:i<\omega\rangle$ by a thin enough tree.) The $\name{\eta}_i$
should ``grow'' fast enough, so naturally we think of forcing notions as in
Ros{\l}anowski and Shelah \cite{RoSh:470}, \cite{RoSh:670}, which proved
easily checked sufficient conditions for what we desire (so in the
``neighborhood'' of Blass and Shelah \cite{BsSh:242}). But what should be
the norm?  

\[* \qquad * \qquad *\]

\begin{discussion}
\label{3.7}
Concerning ${\mathfrak p}<{\mathfrak t}$, I have made quite a few failed
tries.  Some try to use long iterations $(\ge\aleph_{\omega+1})$ or a new
support. But also I thought that Blass and Shelah \cite{BsSh:257} would be a
reasonable starting point, the point is how to extend $\aleph_1$-generated
filters to a good enough $P$-point. 

That is, trying to force ${\mathfrak p}=\aleph_2$, ${\mathfrak t}=\aleph_3 
={\mathfrak c}$ start, say, with $\bV=\bL$ and use a FS iteration $\langle
\bP_i,\name{\bQ}_j:i\le\omega_3,j<\omega_3\rangle$, where $\name{\bQ}_i$ is
a Cohen forcing adding $\name{r}_i\in {}^\omega 2$ for some $i$'s, and
$\name{\bQ}_i$ is shooting an $\omega$--sequence through a $P$-point filter
(or ultrafilter) on $\omega$ for some $i\ge\omega_2$. The point is that when
we have to find a $\le^*$--lower bound to the downward directed ${\mathcal
A} \in [{\mathcal P}(\omega)]^{\aleph_1}$, we extend it to a $P$-point,
possibly also for the $\omega_2$--towers we have to do this. It is natural
to try to preserve, for $\alpha \in [\omega_2,\omega_3)$, the statement: 
\begin{quotation}
in $\bV^{\bP_\alpha}$, noting that ${\mathcal H}(\aleph_1)$ has cardinality
$\aleph_2$, if ${\mathcal H}(\aleph_1)=\bigcup\limits_{\alpha<\omega_2}
M_\alpha$, $M_\alpha$ increasing continuous, $\|M_\alpha\|<\aleph_2$,\\
{\em then\/} the following set is $=\emptyset\mod {\mathcal D}_{\omega_2}+
S^2_1$:  
\[\begin{array}{lr}
\big\{\delta:&\mbox{if some } a\in M_\delta\cap [\omega]^{\aleph_0}\mbox{ is
almost included in }r^{-1}_i(\{1\})\mbox{ for many } i<\delta,\quad\\
&\mbox{then } a \mbox{ is almost disjoint to } r^{-1}_\delta(\{1\})\ \big\}.
  \end{array}\]
\end{quotation}
\end{discussion}

\[* \qquad * \qquad *\]

\begin{discussion}
\label{3.8}
Concerning \ref{2.10} consider the problem ``every $f:\bR\longrightarrow\bR$
is continuous on a non-null set''. 

We can try to use a forcing notion which looks locally random (like the
forcing for ``non meager set'' of \cite{Sh:473} looked locally like the
Cohen forcing notion) or a mixture of random and quite bounding ones. Such
forcing notions are considered in \cite{RoSh:670}, do they help for ``every
function $f:{}^\omega 2\longrightarrow {}^\omega 2$ is continuous on a
non-null set''?  

How can we try to prove the consistency of ``for every non-meagre
$A\subseteq {}^\omega 2$ and $f:{}^\omega 2\longrightarrow {}^\omega 2$ for
some non-meagre $B\subseteq A$, $f\restriction B$ is continuous''?  We may
use CS or even FS iteration of length $\omega_2$, (with $\bV\models{\rm GCH}
+ \diamondsuit_{\{\delta<\aleph_2:\cf(\delta)=\aleph_1\}}+S_\alpha\subseteq
\omega_1$ ($\alpha<\omega_2$) increasing $\mod{\mathcal D}_{\omega_1}$).

In Stage $\alpha$ we have $\bar{r}^\alpha=\langle\name{r}^\alpha_i:i\in
S_\alpha\rangle$ such that $\bar{r}^\alpha\in N\prec ({\mathcal
H}(\chi),\in)\ \Rightarrow\ \name{r}^\alpha_{N\cap\omega_1}$ is forced to be
Cohen over $N$ and $\beta<\alpha\ \Rightarrow\ \{i\in
S_\beta:\name{r}^\beta_i \ne \name{r}^\alpha_i\}$ is not stationary. 

Sometimes in Stage $\alpha$, bookkeeping gives us a $\bP_\alpha$-name
$\name{A}_\alpha$ of a non-meagre subset of ${}^\omega 2$ and we choose
$\bar{r}^{\alpha+1}$ such that $\bar{r}^{\alpha+1}\restriction S_\alpha=
\bar{r}^\alpha$ and $\{r^{\alpha+1}_i:i\in S_{\alpha+1}\setminus S_\alpha\}$
is a non-meagre subset of $\name{A}_\alpha$.

Sometimes in Stage $\alpha$, bookkeeping gives us a stationary subset
$S^\alpha$ of $S_\alpha$ (from $\bV$) and a $\bP_\alpha$--name
$\name{f}_\alpha$ of a function from $\{\name{r}^\alpha_i:i\in S^\alpha\}$
to ${}^\omega 2$ and we try to choose $\bar{r}^{\alpha+1}$ such that:
$\bar{r}^{\alpha+1}\restriction S_\alpha=\bar{r}^\alpha$, $\{\name{r}^{
\alpha +1}_i:i\in S_{\alpha+1}\setminus S_\alpha\}\subseteq\{
\name{r}^\alpha_i:i\in S^\alpha\}$ and $\name{f}_\alpha\restriction
\{\name{r}^{\alpha+1}_i:i\in S_{\alpha+1}\setminus S_\alpha\}$ is
continuous. So the aim is that in $\bV^{\bP_{\omega_2}}$, every non-meagre
$A\subseteq {}^\omega 2$ contains a subset of the form $\{\name{r}^\alpha_i
:i\in S_{\alpha+1}\setminus S_\alpha\}$ and 
\[S'\subseteq S_\alpha\ \&\ \alpha<\omega_1\ \&\ S'\mbox{ stationary }\
\Rightarrow\ \{r^\alpha_i:i\in S'\}\mbox{ is non-meagre.}\]
We may try to define iterations for forcing related to measure: we can use
CS, or try to imitate the measure algebra, there are various ways to
interpret it. If each $\bQ_i$ is as in \cite{RoSh:470}, so each condition
has possible pair: a norm $\in\omega$ and real $r\in (0,1)$ and using those
we define what is a condition in the iteration.  See \cite{RoSh:670}.
\end{discussion}

More on preservation (for not necessarily c.c.c.~ones) commutativity,
associativity, generic sets and countable for pure/finite for a pure support
iterations see \cite{Sh:669}. Remember that an automorphism $F$ of
${\mathcal P}(\omega)$/finite is called trivial if it is induced by a
permutation $f$ of $\mathbb Z$ (the integers, where $\omega=\{n\in{\mathbb
Z}:n \ge 0\}$) such that $\{n\in{\mathbb Z}:n<0\ \Leftrightarrow\ f(n)\ge
0\}$ is finite. 

\begin{question}
\label{3.9}
What can ${\rm AUT}({\mathcal P}(\omega)/{\rm finite})$ be? Is it consistent
that 
\begin{quotation}
${\rm AUT}({\mathcal P}(\omega)/{\rm finite})$ is not the group of trivial
automorphisms of ${\mathcal P}(\omega)/{\rm finite}$, but is of cardinality
continuum (or even is generated by adding one automorphism to the subgroup
of the trivial automorphisms of ${\mathcal P}(\omega)/{\rm finite}$) ?
\end{quotation}
\end{question}

It is reasonable try to combine Shelah and Stepr\=ans \cite{ShSr:427} and the
later part of the proof in \cite[Ch.IV,\S6]{Sh:f} (from being ``locally
trivial'' to being trivial).

\[* \qquad * \qquad *\]

\begin{discussion}
\label{3.10}
Next we deal with the variants of OCA and isomorphisms or farness of sets of
$\aleph_2$ reals, i.e.~\ref{2.20}, \ref{2.23}. 

Concerning \ref{2.20} (on `` no far $A,B\in [{}^\omega 2]^\lambda$ ''),
assume that for $\ell=1,2$ we have $A_\ell=\{\eta^\ell_\alpha:\alpha < 
\lambda\} \subseteq {}^\omega 2$ with no repetitions.  Considering
Baumgartner \cite{B} and Abraham, Rubin and Shelah \cite{ARSh:153}, it is
natural to try to find $\bar{f}=\langle f_\alpha:\alpha<\lambda\rangle$ such
that: 
\begin{enumerate}
\item[(a)]  $f_\alpha$ is a partial, countable, non-empty function from
$\omega_2$ to $\omega_2$, (for the present aim, $\Dom(f_\alpha)$ a singleton
in (a) and $\gamma=1$ in (b) are fine, so we assume so),
\item[(b)]  for some $\gamma=\gamma^*\le\omega$, the sequence 
\[\langle\bigcup\limits_{n<\gamma}(\Dom(f_{\gamma\alpha+n})\cup\Rang(f_{
\gamma\alpha+n})):\alpha<\lambda\rangle\]
is a sequence of pairwise disjoint sets.
\end{enumerate}
We let $\hat{f}_\alpha=\{(\eta^1_i,\eta^2_{f_\alpha(i)}):i\in\Dom(f_\alpha)
\}$, so $\hat{f}_\alpha(\eta)$ is well defined iff $i\in\Dom(f_\alpha)$,
$\eta=\eta^1_i$, and then $\hat{f}_\alpha(\eta^1_i)=\eta^2_{f_\alpha(i)}$. 

It is natural to try the following forcing notion
\[\begin{array}{lr}
\bQ_{\bar{f}}=\big\{g:&g\mbox{ is a finite 1-to-1 order preserving function
  from } A_1 \mbox{ to } A_2,\quad \\ 
&\mbox{which has the form }\bigcup\limits^n_{\ell=1}f_{\alpha_\ell},\
  \alpha_\ell<\omega_2,\mbox{ and such that}\quad\\ 
&\ell_1\ne\ell_2\quad\Rightarrow\quad\neg(\exists\alpha)(\{\alpha_{\ell_1},
\alpha_{\ell_2}\}\subseteq [\gamma^*\alpha,\gamma^*\alpha+\gamma^*])\big\}.
  \end{array}\]
The order is the inclusion.

It is enough to have ``$\bQ_{\bar{f}}$ satisfies the c.c.c.'' (for some
$\bar{f}$ as above): clearly $\bigcup\{g:g\in\name{G}_{\bQ_{\bar{f}}}\}$ is
an order preserving function from some $A'_1\subseteq A_1$ into $A_2$; but
does it have cardinality $\lambda$? Essentially yes, as e.g.~if
$\cf(\lambda)>\aleph_0$ then some $p\in\bQ_{\bar{f}}$ forces this. For this
it is enough to have: if $u\in [\lambda]^{\aleph_1}$ then $\bQ_{\bar{f}
\restriction u}$ satisfies the c.c.c.~and $\lambda=\aleph_2$, it is enough
to check for $u=\alpha\in [\omega_1,\omega_2)$. So as in \cite{ARSh:153}, it
is enough that
\begin{enumerate}
\item[(c)]  if $n<\omega$ and $C\subseteq {}^{2n}({}^\omega 2)$ is closed
 and
\[(\eta_0,\eta_1,\ldots,\eta_{2n-1})\in C\quad \Rightarrow\quad (\forall\ell
<m<n)(\eta_{2\ell}<_{\ell ex}\eta_{2m}\ \Leftrightarrow\ \eta_{2\ell+1}
<_{\ell ex}\eta_{2m+1})\]
and there are $p_\zeta=\{(\eta_{\zeta,\ell},\nu_{\zeta,\ell}):\ell\le n\}
\in \bQ_{\bar{f}}$ for $\zeta<\omega_1$, $\ell<n$ with 
\[(\zeta_1,\ell_1)\ne (\zeta_2,\ell_2)\quad\Rightarrow\quad \eta_{\zeta_1,
\ell_1}\ne\eta_{\zeta_2,\ell_2}\mbox{ and }\nu_{\zeta_1,\ell_1}\ne\nu_{
\zeta_2,\ell_2},\]
hence 
\[(\forall\zeta<\omega_1)(\forall\ell<m<n)(\eta_{\zeta,\ell}<_{\ell ex}
\eta_{\zeta,m}\ \Leftrightarrow\ \nu_{\zeta,\ell}<_{\ell ex}\nu_{\zeta,m})\]
(if $\gamma^*=1$, and each $\Dom(f_\alpha)$ is a singleton, $p_\zeta\in
\bQ_{\bar{f}}$ means, in addition only $\{(\eta_{\zeta,\ell},\nu_{\zeta,
\ell})\}=f_{\alpha_{\zeta,\ell}}$ for some $\alpha_{\zeta,\ell}$), 

\noindent {\em then\/} there are $(\eta'_0,\eta'_1,\ldots,\eta'_{2n-1}),
(\eta''_0,\eta''_1,\ldots,\eta''_{2n-1})\in C$ such that $\eta'_\ell\ne
\eta''_\ell$ and $\eta'_{2\ell}<_{\ell ex} \eta''_{2 \ell}\ \Leftrightarrow\
\eta'_{2\ell+1}<_{\ell ex} \eta''_{2\ell}$.
\end{enumerate}
Any counterexample to clause (c) induces a continuous partial function for
which we get dependencies. If $\bV\models {\rm CH}$, and $\bP$ is adding
Cohen reals, then in $\bV^\bP$ this holds, so it is natural to try to retain
during the iteration similarly to this case. For \ref{2.19} we can use
similar, but somewhat more involved forcing notion as in \cite{ARSh:153}.  
\end{discussion}

There are similar considerations on OCA. We may consider trying to get
negative ZFC results, so Kojman and Shelah \cite{KjSh:409} seems to me a
reasonable starting point (of course, the problem there is different). 

\[* \qquad * \qquad *\]

Baumgartner \cite{B6} defines

\begin{definition}
\label{upi.0} 
\begin{enumerate}
\item For a (non-principal) ultrafilter $D$ on $\omega$, and a countable
ordinal $\delta$, we say $D$ is a $\delta$--ultrafilter if: 
\begin{enumerate}
\item[$(*)$]  for every function $f$ from $\omega$ to $\omega_1$, for some
$A\in D$ we have $\otp(f(A))<\delta$.
\end{enumerate}
\item We say that $D$ is a weak $\delta$--ultrafilter if:
\begin{enumerate}
\item[$(*)^-$]  for every function $f$ from $\omega$ into $\delta$, for
some $A\in D$ we have $\otp(f(A))<\delta$.
\end{enumerate}
\item We say that $D$ is a NWD if for every function $f$ from $\omega$ to
$\bR$, for some $A\in D$, $f(A)$ is a nowhere dense subset of $\bR$.
\item For an ideal $I$, $D$ is an $I$--filter if for any $f:\Dom(D)
\longrightarrow\Dom(I)$ and $A\in D^+$ there is $B\subseteq A$, $B\in D^+$
such that $f(B)\in I$.
\end{enumerate}
\end{definition}
Then he asked whether such ultrafilters exist (if CH yes, so):

\begin{question}
\label{upi.0A}
Prove the consistency of ``there is no $\delta$--ultrafilter on $\omega$''.
\end{question}
It seemed the solution of the related ``CON(there is no NWD-ultrafilter)''
in \cite{Sh:594} should give this, but it did not, and it is not clear if
the question \ref{upi.0A} is harder (the NWD eluded me several times, but
when solved, the solution seems a straightforward generalization of CON(no
$P$-point), which was also a priori the natural starting point).

A nice feature of $P$-points is that ``$D$ generates a $P$-point ultrafilter
on $\omega$'' is preserved in limit for CS iterations, so $P$-points
generated by $\aleph_1<2^{\aleph_0}$ sets are gotten naturally. Are they the
only ones? Of course, by c.c.c.~forcing $P$ you may have ultrafilters on
$\omega$ generated by $<2^{\aleph_0}$ sets, and forcing by a subforcing
$\bQ\lesdot\bP$, in $\bV^\bQ$ we get an ultrafilter preserved (see more in
Brendle and Shelah \cite{BnSh:642}); but we have no understanding, though
the suggestion in \ref{3.2} may help. 

\begin{question}
\label{upi.0B}
Are ultrafilters $D$ as defined below in \ref{upi.8} preserved in limit
stages of CS iterations?  This means that: 
\begin{quotation}
{\em if\/} $\bar{\bQ}=\langle\bP_i,\name{\bQ}_j:i\le\delta,j<\delta\rangle$
is a CS iteration of proper forcing notions and $D$ is an ultrafilter as
above in $\bV$,\\
{\em then\/}\quad $(\forall i<\delta)(\boxtimes_i)\ \ \Rightarrow\ \ (
\boxtimes_\delta)$, where
\begin{enumerate}
\item[$(\boxtimes_\beta)$]  $\Vdash_{\bP_\beta}$ ``in $\bV^{\bP_\beta}$ the
filter on $\omega$ that the family $D$ generates in $\bV^{\bP_\beta}$ is an
ultrafilter". 
\end{enumerate}
\end{quotation}
\end{question}
Less nice, but still good, is to prove the preservation of ``$D$ generates
an ultrafilter $+{\rm Pr}$'' (where ${\rm Pr}$ is some additional property
like: $({}^\omega \omega)^\bV$ is dominating).

The following \ref{upi.2} --- \ref{upi.9} suggest an approach to question
\ref{upi.0B}.

\begin{definition}
\label{upi.2}
\begin{enumerate}
\item Let
\[\begin{array}{ll}
{\mathcal T}=\big\{t:&t\subseteq {}^{\omega>}\omega,\ t\mbox{ has a }
  \vartriangleleft\mbox{--minimal element }{\rm rt}(t),\\
  &t\mbox{ is closed under initial segments of length }\ge\ell g({\rm
  rt}(t)), \\ 
  &\mbox{for }\eta\in t,\ {\rm Suc}_t(\eta)=\{\eta\conc\langle\ell\rangle:
  \eta\conc\langle\ell\rangle\in t\}\mbox{ is empty or infinite}\big\}. 
  \end{array}\]
For $t\in {\mathcal T}$ let $h_t:t\longrightarrow\omega_1\cup\{\infty\}$ be
defined by
\[h_t(\eta)=\bigcup\big\{h_t(\nu)+1:\nu\in{\rm Suc}_t(\eta)\big\}.\]
(So $h_t(\eta)=\infty$ iff there is an $\omega$-branch through $\eta$).

\noindent We say that $t$ is standard if 
\[\eta\in t\ \&\ \beta<h_t(\eta)\quad\Rightarrow\quad (\forall^\infty \nu\in
{\rm Suc}_t(\eta))(\beta\le h_t(\nu)).\]
If not said otherwise, every $t$ is standard.  Note: 
\[s\in{\rm sub}(t)\ \&\ t\mbox{ is standard}\quad\Rightarrow\quad s\mbox{ is
standard},\]
where on ${\rm sub}(t)$ see part (8) below. 

\item For an ordinal $\alpha < \omega_1$, let
\[\begin{array}{lcl}
{\mathcal T}_\alpha&=&\big\{t\in {\mathcal T}:\Rang(h_t)\subseteq\omega_1
\mbox{ and } h_t({\rm rt}(t))=\alpha\big\},\\
{\mathcal T}_{<\alpha}&=&\bigcup\limits_{\beta<\alpha}{\mathcal T}_\beta.
  \end{array}\]
\item For $t \in {\mathcal T}$ let
\[\begin{array}{lr}
{\mathcal A}_t=\big\{\bar{A}:&\bar{A}=\langle A_\eta:\eta\in t\rangle,\
  A_\eta \in [\omega]^{\aleph_0}, \mbox{ and if }\eta\in t\mbox{ is not
  maximal, }\ \\
&\mbox{ then }\langle A_\nu:\nu\in{\rm Suc}_T(\eta)\rangle\mbox{ is a
sequence of pairwise disjoint }\ \\
&\mbox{subsets of } A_\eta\big\},
  \end{array}\]
\[\begin{array}{lcl}
{\mathcal A}_\alpha&=&\bigcup\big\{{\mathcal A}_t:t\in {\mathcal T}_\alpha
  \big\},\\ 
{\mathcal A}_{<\alpha}&=&\bigcup\limits_{\beta<\alpha}{\mathcal A}_\beta.
 \end{array}\]
\item For $t\in {\mathcal T}$, $\eta\in t$ let $t^{[\eta]}=\{\nu:\eta
\trianglelefteq\nu\in t\}\in {\mathcal T}$. 
\item For $t\in {\mathcal T}$ let $\lim(t)=\{\eta\in {}^\omega\omega:
(\forall \ell<\omega)(\ell\ge \ell g({\rm rt}(t))\ \Rightarrow\ \eta
\restriction \ell\in t)\}$. 
\item $\max(t)=\{\eta\in t:{\rm Suc}_t(\eta)=\emptyset\}$.
\item We say that $y$ is a front of $t\in {\mathcal T}$ {\em if\/}:\ $y
\subseteq t$, $(\forall\eta,\nu\in y)(\neg\eta\vartriangleleft\nu)$ and
\[(\forall\eta)[(\eta\in\lim(t)\vee\eta\in\max(t))\ \Rightarrow\ (\exists
\ell \le\ell g(\eta))(\eta\restriction \ell\in y)].\]
We let ${\rm fr}(t)=\{y:y\mbox{ is a front of } t\}$.
\item For $t\in {\mathcal T}$ let ${\rm sub}(t)=\{s\in {\mathcal T}:
s\subseteq t \mbox{ and }\max(s)=\max(t)\cap s\}$.

\noindent Clearly for a standard $t$ and $s\in{\rm sub}(t)$ we have $h_s=
h_t\restriction s$. Let 
\[{\rm sub}^-(t)=\{s\in {\mathcal T}:s\subseteq t,\ s\in{\rm sub}(t)\mbox{ 
and }{\rm rt}(s)={\rm rt}(t)\}.\]
\item For $s,t\in {\mathcal T}$ and a set ${\mathcal Y}$ of fronts of $t$
let ${\rm FT}^0_{\mathcal Y}(t,s)$ be the set of embeddings $f:s
\longrightarrow t$ (i.e., $f$ is one--to--one, $\Dom(f)=s$, $\Rang(f)
\subseteq t$, $f({\rm rt}(s))={\rm rt}(t)$ and $(\forall\eta,\nu\in t)(\eta
\vartriangleleft\nu\ \Rightarrow\ f(\eta)\vartriangleleft f(\nu))$) which
respect each $y\in {\mathcal Y}$, i.e., $\{\eta:f(\eta)\in y\}$ is a front
of $s$ for every $y\in {\mathcal Y}$.  
\item If ${\mathcal Y}$ is the set of all fronts of $t$ we may omit it.
\item ${\rm FT}^1(t,s)$ is the set of $f\in {\rm FT}^0(t,s)$ such that 
\[\eta \in{\rm Suc}_s({\rm rt}(s))\quad \Rightarrow\quad h_s(\eta)=h_t(f
(\eta)).\] 
\item ${\rm FT}^2(t,s)$ is the set of all $f:s\longrightarrow t$ such that:
\[\begin{array}{l}
\eta\vartriangleleft\nu\ \Leftrightarrow\ f(\eta)\vartriangleleft f(\nu),\\
\{f(\eta):\eta\in{\rm Suc}_s({\rm rt}(s))\}\mbox{ is a front of }t,\ \mbox{
and}\\
\mbox{for }\eta\in{\rm Suc}_s({\rm rt}(s)),\quad f\restriction s^{[\eta]}
\mbox{ is one-to-one onto }t^{[f(\eta)]}.
  \end{array}\]
\item For an ideal $I$ let
\[\begin{array}{lr}
{\mathcal S}{\mathcal T}={\mathcal S}{\mathcal T}_I=\big\{(t,g):&t\in
{\mathcal T}\mbox{ and } g:t\setminus\{{\rm rt}(t)\}\longrightarrow\Dom(I)
\mbox{ are such that }\ \\
&\eta\in t\setminus\max(t)\quad\Rightarrow\quad\{g(\nu):\nu\in{\rm Suc}_T
(\eta)\}\in I^+\big\}.
  \end{array}\]
We usually omit $I$ if it is clear from the context (here it is fixed).
\item For $(t^1,g^1),(t^2,g^2)\in {\mathcal S}{\mathcal T}$ for $\ell=0,1,2$
let ${\mathcal F}{\mathcal T}^\ell((t^1,g^1),(t^2,g^2))$ be the set of $f\in
{\rm FT}^\ell(t^1,t^2)$ such that: $g^2=g^1\circ f$.
\item Let
\[\begin{array}{lcl}
{\mathcal T}{\mathcal T}_t&=&\big\{(t,\bar{A},g):\bar{A}\in {\mathcal A}_t
\mbox{ and }(t,g)\in {\mathcal S}{\mathcal T}\big\},\\
{\mathcal T}{\mathcal T}&=&\bigcup\limits_t {\mathcal T}{\mathcal T}_t,\\
{\mathcal T}{\mathcal T}_\alpha&=&\big\{(t,\bar{A},g):t\in {\mathcal
T}_\alpha\mbox{ and }(t,\bar{A},g)\in {\mathcal T}{\mathcal T}_t\big\},\\
{\mathcal T}{\mathcal T}_{<\alpha}&=&\bigcup\limits_{\beta<\alpha}{\mathcal
T}{\mathcal T}_\beta. 
  \end{array}\]
\item Let ${\mathcal F}{\mathcal T}^\ell((t^1,\bar{A}^1,g^1),(t^2,\bar{A}^2,
g^2))$ be the set of $f\in{\mathcal F}{\mathcal T}^\ell((t^1,g^1),(t^2,
g^2))$  such that $\eta\in t^2\quad \Rightarrow\quad A^1_{f(\eta)}
\subseteq A^2_\eta$. 
\item If we omit $\ell$ (in ${\rm FT}^\ell,{\mathcal F}{\mathcal T}^\ell$)
we mean $\ell = 0$.
\end{enumerate}
\end{definition}

\begin{definition}
\label{upi.3}
We define a partial order on ${\mathcal T}{\mathcal T}_{< \alpha}$:
\[(t^1,\bar{A}^1,g^1)\le^\ell (t^2,\bar{A}^2,g^2)\quad \mbox{ if and only if
}\quad{\mathcal F}{\mathcal T}^\ell((t^1,\bar{A}^1,g^1),(t^2,\bar{A}^2,g^2))
\neq\emptyset.\]
\end{definition}

\begin{observation}
\label{upi.4}
$\leq^\ell$ really is a partially order of ${\mathcal T}{\mathcal T}$.
\end{observation}

\begin{definition}
\label{upi.5}
\begin{enumerate}
\item For $t\in {\mathcal T}$, let
\[\begin{array}{lr}
{\rm fsub}(t)=\{s\in{\mathcal T}:&\mbox{for some finite } w\subseteq{\rm
Suc}_t({\rm rt}(t))\mbox{ we have }\ \\
&s=\{\eta\in t:\neg(\exists\nu\in w)(\nu\trianglelefteq \eta)\}\}.
  \end{array}\]
\item For $(t^1,\bar{A}^1,g^1),(t^2,\bar{A}^2,g^2)\in {\mathcal T}{\mathcal
T}$ we define ${\mathcal F}{\mathcal T}^\ell_*((t^1,\bar{A}^1,g^1),(t^2,
\bar{A}^2,g^2))$ as the set of all $f$ such that for some $t^3\in {\rm
fsub}(t^2)$ we have  
\[f\in {\mathcal F}{\mathcal T}^\ell((t^1,\bar{A}^1,g^1),(t^3,\bar{A}^2
\restriction t^3,g^2\restriction t^3)).\]
\item $(t^1,\bar{A}^1,g^1)\le^\ell_* (t^2,\bar{A}^2,g^2)$ if and only if
${\mathcal F}{\mathcal T}^\ell_*((t^1,\bar{A}^1,g^1),(t^2,\bar{A}^2,g^2))$
is not empty. 
\end{enumerate}
\end{definition}

\begin{fact}
\label{upi.6}
\begin{enumerate}
\item $<^\ell_*$ is a partial order of ${\mathcal T}{\mathcal T}$ such that:
$<^\ell$ is a subset of $<^\ell_*$. 
\item Any $<^\ell_*$--increasing chain of length $\omega$ in ${\mathcal
T}{\mathcal T}_{\le \alpha}$ has an upper bound in ${\mathcal T}{\mathcal
T}_{\le \alpha}$. 
\item ${\mathcal T}{\mathcal T}_{\le \omega_1} = {\mathcal T}{\mathcal T}_{<
\omega_1}$. 
\end{enumerate}
\end{fact}

\begin{proposition}
\label{upi.6a}
\begin{enumerate}
\item If $(t,\bar{A},g)\in {\mathcal T}{\mathcal T}_{\le\alpha}$, and
$B\subseteq \omega$, 

\noindent {\em then\/} for some $(t',\bar{A}',g')\in {\mathcal T}{\mathcal
T}_{\le\alpha}$ we have:  
\begin{enumerate}
\item[$(\alpha)$] $(t,\bar{A},g)\le^1 (t',\bar{A}',g')$, in fact $t'\in{\rm
sub}(t)$, 
\item[$(\beta)$]  $\bigcup\{A_\eta:\eta\in t\setminus\{{\rm rt}(t)\}\}$ is a
subset of $B$ or is disjoint to $B$. 
\end{enumerate}
\item Similarly for ${\mathcal T}{\mathcal T}$.
\end{enumerate}
\end{proposition}

\begin{proof}  
(1)\quad By induction on $\alpha$.\\
(2)\quad Similarly.
\end{proof}

\begin{proposition}
\label{upi.7}
If $(t,\bar{A},g)\in {\mathcal T}{\mathcal T}_{\le \alpha}$ and $E$ is an
equivalence relation on $\omega$, then for some $(t,\bar{A}',g')$ and front
$y$ of $t'$ we have: 
\begin{enumerate}
\item[$(\alpha)$] $(t,\bar{A},g)\le^1 (t',\bar{A}',g')$,
\item[$(\beta)$]  for $\eta\in y$, $A'_\eta$ is included in one
$E$--equivalence class, 
\item[$(\gamma)$] for $\eta\ne\nu$ from $y$, the $E$--equivalence classes in
which $A'_\eta,A'_\nu$ are included, are distinct (hence disjoint), 
\item[$(\delta)$] on $A'_{{\rm rt}(t')}$ we have: $E$ is either trivial or
refines $\langle A'_\eta:\eta\in{\rm Suc}_{t'}({\rm rt}(t'))\rangle$.
\end{enumerate}
\end{proposition}

\begin{proof}
By induction on $\alpha$.
\end{proof}

\begin{proposition}
[CH] 
\label{upi.8}
Let $I$ be an ideal such that $|\Dom(I)|\le\aleph_1$ and
\[(\forall X\in I^+)(\exists Y\in I^+)(Y\subseteq X\ \&\ |Y|=\aleph_0).\]
There is a $\leq^1_*$--increasing sequence $(\langle t^\zeta,\bar{A}^1,
g^\zeta\rangle:\zeta<\omega_1)$ of members of ${\mathcal T}{\mathcal T}_{\le
\omega_1}$ such that $\{\bigcup\limits_{\eta\in t^\zeta} A^1_\eta:\zeta<
\omega_1\}$ generates a non-principal ultrafilter $D$ on $\omega$ which is a
$Q$-point, and for every equivalence relation $E$ on $\omega$, for some
$\zeta$, $\langle A^\zeta_\eta:\eta\in{\rm Suc}_{t^\zeta}(rt(t^\zeta))
\rangle$ refines $E\restriction A^\zeta_{{\rm rt}(t^\zeta)}$. 
\end{proposition}

\noindent{\sc Remark:}\quad  This construction gives an ultrafilter $D$ on
$\omega$, a $Q$--point such that 
\[D'\le_{\rm RK} D\ \Rightarrow\ D'\mbox{ is not NWD.}\]

\begin{proposition}
\label{upi.9}
In \ref{upi.8}, if in addition an ideal $I'$ satisfies
\begin{enumerate}
\item[$(*)$]  $|\Dom(I')|\le\aleph_1$ and if $(t,\bar{A},g)\in {\mathcal
T}{\mathcal T}_{<\omega_1}$ and $g':\omega\longrightarrow\Dom(I')$

\noindent {\em then\/} for some $(t',\bar{A}',g')\ge (t,\bar{A},g)$ we have
$g'[A'_{\langle\rangle}]\in I'$, 
\end{enumerate}
{\rm then\/} we can demand that $D$ is an $I'$-ultrafilter (see
\ref{upi.0}(4)). 
\end{proposition}

\[* \qquad * \qquad *\]

There are many problems on the $\sigma$-versions of cardinal invariants, and
I think for some the method of \cite{Sh:592}, \cite{Sh:619} is relevant, e.g.

\begin{question}
[See Brendle and Shelah \cite{BnSh:642}]
\label{3uf.1}
Does $\chi_\sigma(D)=\chi(D)$ for all ultrafilters $D$ on $\omega$? Recall
that 
\[\begin{array}{rl}
\chi(D)=\min\{|{\mathcal A}|:&{\mathcal A}\subseteq D\mbox{ and for every
 } A \in D\mbox{ for some } B \in {\mathcal A}\mbox{ we have } B \subseteq
 A\}\\ 
\chi_\sigma(D)=\min\{|{\mathcal A}|:&{\mathcal A}\subseteq {}^\omega D
\mbox{ is such that for every } \bar{A}\in {}^\omega D,\mbox{ for some }\ \\
&\quad\bar{B}\in {\mathcal A}\mbox{ we have } (\forall n<\omega)(\exists m< 
\omega)(B_m\subseteq^* A_n)\}.
\end{array}\]
\end{question}

So a reasonable scenario to prove the consistency of a negative answer runs
as follows: let, e.g., $\mu=\aleph_\omega$. We use FS iteration of
c.c.c.~forcing notions, $\langle\bP_i,\name{\bQ}_j: i\le\delta^*,j<\delta^*
\rangle$. We have $\bP_i$-names $\name{D}^i_u,\name{A}_{u,\gamma}$ (for
$\gamma<\gamma^i_u$ and $u\in[\mu]^{<\aleph_0}$) such that:
\begin{itemize}
\item $\name{D}^i_u$ is the filter on $\omega$ generated by $\{\name{A}_{u,
\gamma}:\gamma<\gamma^i_u\}$ and the co-bounded sets, 
\item $\gamma^i_u$ are increasing with $i$, $\name{D}^i_u\subseteq\name{D}^i_v$
if $u\subseteq v$, and $\name{D}^j_u\subseteq\name{D}^i_u$ for $j<i$.
\end{itemize}
To simplify we decide:
\begin{enumerate}
\item[$(*)$]  if $j<i$, $\Vdash_{\bP_j}$`` $\name{A}\subseteq\omega\mbox{
and } u\in [\mu]^{<\aleph_0}\mbox{ and }\name{A}\subseteq^* \name{A}_{u,
\gamma}\mbox{ for every }\gamma<\gamma^j_u$ ''

\noindent{\em then\/} $\Vdash_{\bP_i}$`` for every $v\in [\mu]^{<\aleph_0}$,
$\name{A}\notin\name{D}^i_v$ ''. 
\end{enumerate}
Also for the following it seems reasonable to try to be influenced by
\cite{Sh:592}, \cite{Sh:619}. 

\begin{question}
[See Brendle and Shelah \cite{BnSh:642}].  
\label{3uf.2}
Can $\pi\chi(D)$ be singular, where
\[\begin{array}{lr}
\pi\chi(D)=\min\{|{\mathcal A}|:&{\mathcal A}\subseteq[\omega]^{\aleph_0},
\mbox{ and for every } B\in D\quad\\
&\mbox{for some } A\in {\mathcal A}\mbox{ we have }  A \subseteq^* B\}. 
  \end{array}\]
\end{question}

\section{Nicely defined forcing notions}
Ros{\l}anowski and Shelah \cite{RoSh:470}, and \cite{Sh:630}, \cite{Sh:f}
relate as algebraic three dimensional varieties relate to manifolds in
$\bR^n$ and these, in turn, relate to general topology. In \cite{Sh:630}
(on nep and snep) and in Judah and Shelah \cite{JdSh:292} (on Souslin
forcings) we deal with forcing notions defined in an absolute enough way; in
\cite{RoSh:470} (more in \cite{RoSh:670}, \cite{RoSh:672}) with forcing
notions defined in an explicit way (say as tress and generally by
creatures), in \cite{Sh:f} we deal with forcing notions related to the
continuum. 

Our problem with speaking about \cite{RoSh:470}, \cite{RoSh:628},
\cite{RoSh:670} and \cite{RoSh:672} is that much work is in progress, still
orthogonal to it is the question whether in the main theorems of
\cite{RoSh:470}, all the assumptions are needed. That is, within the
framework of condition trees or $\omega$-sequences of creatures, are the
demands on the norms necessary? This is dealt with for the conditions for
properness in \cite{RoSh:470}, showing necessity but there are still gaps
remaining.  

\begin{question}
\label{4.1}
Are the sufficient conditions for properness in \cite[\S2]{RoSh:470}
necessary?  The test case (chosen in \cite{RoSh:470}) is 
\[\bQ=\{\langle w_n:n<\omega\rangle:w_n\subseteq 2^n,w_n\ne\emptyset\mbox{
and }\lim_{n\to\omega}|w_n|\}=\infty\]
ordered by $\bar{w}\le\bar{w}'\ \Leftrightarrow\ (\forall n\in\omega)(w'_n
\subseteq w_n)$. 
\end{question}
Though properness is the main thing and there we look for counterexamples
only for properness, it is interesting to know:

\begin{question}
\label{4.2}
Concerning other theorems of \cite{RoSh:470}, are they sharp?
\end{question}
There are more specialized problems, probably solvable in this context.

\begin{question}
\label{4.2A}
Is there an ${}^\omega\omega$--bounding forcing notion adding a perfect set
of random reals? 
\end{question}
It seems this should not be hard if true.

The following problems (raised by Komjath and Stepr\=ans respectively) seem to
me a matter of choosing the right variant of \cite{RoSh:470} or
\cite{RoSh:670} and having the right finite combinatorics.

\begin{question}
\label{4.2B}
\begin{enumerate}
\item Can each $A\in [{}^\omega 2]^{\aleph_1}$ be null while the union of
some $\aleph_1$ lines in $\bR\times\bR$ is not null? 
\item For reals $0<a_0<a_1\le 1$, is it consistent with ZFC that:\quad for
$\ell<2$,
\[\ell=0\ \mbox{ iff\ \ some }A\in [\bR]^{\aleph_1}\mbox{ has positive
Hausdorff capacity for }a_\ell\mbox{ ?}\]
\end{enumerate}
\end{question}

\[* \qquad * \qquad *\]

I suppose that the feeling that the Cohen forcing notion and the random real
forcing notion occupy a special place is old; probably more in the version
speaking on the ideal of null sets and the ideal of meagre sets. I feel the
former version is more interesting. For me this translates to

\begin{problem}
\label{4.3}
Among Souslin c.c.c.~forcing notions, are Cohen forcing and random forcing
special? 
\end{problem}
Some progress was made in \cite{Sh:480}.

\begin{theorem}
\label{4.4}
If a Souslin c.c.c.~forcing notion $\bQ$ adds $\name{\eta}\in{}^\omega\omega$
not dominated by any old $\nu\in {}^\omega\omega$, {\em then\/} forcing with
$\bQ$ adds a Cohen real.\\ 
(The ``Souslin'' is needed for enough absoluteness, so with the existence of
large cardinals we can allow a larger family).   
\end{theorem}
So the Cohen forcing notion is the minimal one among Souslin c.c.c.~forcing
notions adding an undominated real, so it is natural to conjecture: 

\begin{problem}
\label{4.5}
Show that any Souslin c.c.c.~forcing notion adding a real adds a Cohen real
or adds a random real. 
\end{problem}
This really will show that Cohen and random are special.

In a sense the realm of Souslin c.c.c.~forcing notion can be looked as being
divided between the ${}^\omega\omega$-bounding (with Random forcing as
prototype) and those forcing notions adding an undominated real (with Cohen 
forcing as prototype); we can further distinguish those adding a dominating 
real.

However, the situation is very unbalanced: among Souslin c.c.c.~forcing
notions adding an undominated real we have many examples and a $<^*$-minimal
one, Cohen, (see Definition \ref{4.6} below). 

On the other side we have no idea what occurs among the
${}^\omega\omega$--bounding ones: probably random real is the unique one,
but it is not out of the question that there is a plethora (adding one or
many randoms is an irrelevant distinction; we can even order $(\bQ,
\name{r})$, $\name{r}$ a $\bQ$--name of a real such that the order depends
only on the subforcing $\name{r}$ generates).

\begin{definition}
\label{4.6}
Let $\bQ_1,\bQ_2$ be definitions of forcing notions (absolute enough) say as
in \cite{Sh:630}, or Souslin. 
\begin{enumerate}
\item $\bQ_1\le^*_0 \bQ_2$ if forcing with $\bQ_2$ adds a generic for
$\bQ_1$ and we let $<^*$ mean $<^*_0$. 
\item $\bQ_1\le^*_1\bQ_2$ means: for some $n$, if we force by iteration $n$
times of $\bQ_2$, we add a generic for $\bQ_1$.
\item $\bQ_1\le^*_{2,{\rm fs}}\bQ_2$ is defined similarly using FS iteration
of length $<\omega_1$.  
\item $\bQ_1\le^*_{2,{\rm cs}}\bQ_2$ is defined similarly using CS of length
$<\omega_1$. 
\end{enumerate}
\end{definition}
Note: \ref{4.5} is on the interval between the {\em control measure
problem\/} (see Fremlin \cite{Fe94}) and von Neumann question which says: is
any complete c.c.c.~Boolean Algebra which as a forcing is
${}^\omega\omega$--bounding, a measure algebra. Another way to express the
thought that Cohen and random are special was Kunen's conjecture, see Kunen
\cite{Ku84}, Kechris and Solecki \cite{KeSo95}, Solecki \cite{So96},
\cite{So98} and Ros{\l}anowski and Shelah \cite{RoSh:628}. 

It is natural to investigate the partial orders from \ref{4.6}. So,

\begin{problem}
\label{4.7} 
Investigate the quasi order $\le^*$ (and its variants) for $\bQ$ which are
nep (see \cite{Sh:630}) or which are c.c.c.~$\aleph_0$--nep or which are
c.c.c.~$\aleph_0$--snep.\\
{[We may concentrate on those with a generic real (those are the main
interest for \ref{4.8}(1), (2) below).]}
\end{problem}

An example is (and probably not hard):

\begin{question}
\label{4.7a}
Prove that dominating real forcing (i.e., the Hechler forcing notion) is
$\le^*$--minimal among Souslin c.c.c.~forcing notions adding a dominating
real.\\ 
{[For $\le^*_1$ this is easy.]}
\end{question}
Looking more serious are 

\begin{question}
\label{4.7B}
Can you characterize the $<^*$--minimal $\bQ$, which add a Cohen real but
are not equivalent to the Cohen forcing (hopefully there is one or at least
there are only few).
\end{question}
 
\begin{question}
\label{4.7C} 
Can you characterize the $<^*$--minimal $\bQ$ among the non--minimal $\bQ$
which add a dominating real but are stronger than the Hechler forcing
notion?
\end{question}
A positive solution of \ref{4.5} would also show that the only symmetric
Souslin c.c.c.~forcing notions are the Cohen forcing and the random forcing
(by \cite[\S9]{Sh:630}).  

\[* \qquad * \qquad *\]

Why should we be interested in Souslin proper or in nep forcing or better
yet, why am I? The reason has been iteration theorems; when you are
interested in iterating some very special forcing notions, the proof of
their properness gives more, e.g., the existence of generic conditions over
models occurs also for countable models of versions of ZFC which are not
necessarily $\prec({\mathcal H}(\chi),\in)$. Moreover, some things are
preserved by iterations and this is helpful for specific problems which is
the point of Judah and Shelah \cite{JdSh:292}. [In \cite{JdSh:292} this was 
phrased using descriptive set theory, getting Souslin proper. However, this
does not cover the Sacks forcing notion, the Laver forcing notion, etc.,
which was accomplished by nep.]  

Needless to say, I think iteration theorems for forcing are important and
interesting (otherwise, normally I would not have written a book on the
subject - see \S2).

Another basic reason is that the family of nep forcing notions forms a
natural class. Now, while I feel that general sets are much more basic and
interesting then families of definable ones, and so prefer ${\mathcal
P}({\mathcal P}(\omega))$ to the family of projective sets, certainly they
are interesting and natural. 

Another reason is ``large'' ideals. Let $I$ be a $\kappa$-complete ideal on
$\lambda$. Gitik and Shelah \cite{GiSh:357} start by proving that ${\mathcal
P}(\kappa)/I$ cannot be (the Boolean algebra which up to isomorphism is
equivalent to) the Cohen forcing or random real forcing, an old question
which Fremlin promoted (see \cite{Fe93}), which comes from asking: can the 
classical result of Solovay \cite{So2} (saying that consistently
$2^{\aleph_0}$ is real valued measurable, now the Maharam type there was
large) be improved to get small Maharam type.

But then \cite{GiSh:357} turns to: 

\begin{problem}
\label{4.8}
\begin{enumerate}
\item Prove that ${\mathcal P}(\kappa)/I$ cannot be a Souslin c.c.c.~forcing
generated by the name of one real $\name{\eta}$ (where $I$ is a
$\kappa$-complete ideal on $\kappa$ or at least $\aleph_1$-complete). 
\item Similarly for Souslin proper (or weaken the definability demand -
natural as the existence of the ideal implies more absoluteness).
\item Even reasonable subclasses or cases are interesting.  
\end{enumerate}
\end{problem}

\begin{problem}
\label{4.8A}
Similarly, we can ask about a $\sigma^+$--complete ideal $I$ on $\kappa$
such that ${\mathcal P}(\kappa)/I$ has a dense subset isomorphic to a
partial order defined in $({\mathcal H}_{< \sigma}(\kappa'),\in)$ with
parameters.  
\end{problem}

In Gitik and Shelah \cite{GiSh:357}, \cite{GiSh:412}, \cite{GiSh:582}, 
in addition to information on adding not too many random or Cohen reals, and
(toward \ref{4.8}) to general criteria for impossibility, we consider more
specific cases (see then \cite{Sh:F257}). The problems lead us to
properties of definable forcing notions like symmetry. The theorems on
Cohen and random reals use the symmetry (i.e, the Fubini theorem), but other
properties pop up naturally, e.g., for Souslin c.c.c.~forcing $\bQ$ with a
dominating real $\name{\eta}$ as generic, to show impossibility it suffices
to show: $\Vdash_\bQ$`` ${\mathfrak b}=\aleph_1$ '' (by \cite{GiSh:357}). 
Maybe the work on the ideals is done and we just need to verify that always at
least one of the criteria applies (at least for large subclasses). Now
\cite{Sh:480}, \cite[\S8,\S9]{Sh:630} comes to my mind. 

Considerations like this lead to questions like

\begin{question}
\label{4.8Ax}
Find sufficient conditions on $\bQ$ for ``$\bQ*\name{\mbox{Random}}/
\mbox{Random}$ adds no random real''.\\  
(This question is chosen since it is also interesting because if the
condition is reasonable enough, it suffices for proving ${\rm
CON}(\cov^*({\rm null})<{\rm non}({\rm meagre}))$, see Bartoszy\'nski,
Ros{\l}anowski and Shelah \cite{BRSh:490}, \cite{BRSh:616}.)
\end{question}

\begin{question}
\label{4.8B} 
Investigate commuting pairs (see \cite{Sh:630}).
\end{question}
For such considerations I had felt that a peculiar property of Cohen forcing
and random forcing is their satisfying: ``being a maximal antichain is a
Borel property''; this leads to 

\begin{definition}
\label{4.10}
A forcing notion $\bQ$ is very Souslin c.c.c.~if it is Souslin c.c.c.~and
also the notion of ``$\{r_n:n<\omega\}$ is a maximal antichain'' is
$\Sigma^1_1$. 
\end{definition}

We hope this will turn out to be a good dividing line of the Souslin
c.c.c.~forcing (so helping to prove theorems). This is because I suspect the
answer to the following is yes. 

\begin{question}
\label{4.10A}
Prove: If $\bQ$ is a Souslin c.c.c.~forcing notion, say with generic real
$\name{r}$ and it is not very Souslin c.c.c.~above any $p\in\bQ$ {\em
then\/} $\Vdash_{\bQ}$`` ${\mathfrak b}=\aleph_1"$.  (This should help
\ref{4.8} by \cite[\S4]{GiSh:357}).  See on this \cite{Sh:669}. 
\end{question}
As in \cite{Sh:630} we can define (restricting $\kappa$ to be $\aleph_0$ for
simplicity)

\begin{definition}
\label{4.11}
\begin{enumerate}
\item A forcing notion $\bQ$ is $\omega$--nw--nep if there is a sequence
$\bar{\varphi}=\langle \varphi_0,\varphi_1,\varphi_2 \rangle$ of
$\Sigma^1_1$ definitions such that: 
\begin{enumerate}
\item[(a)] the set of members of $\bQ$ and $\le^\bQ$ are $\Sigma^1_1$ sets
(of reals) defined by $\varphi_0,\varphi_1$, respectively, 
\item[(b)] if $N$ is an $(\bar{\varphi},\omega)$--nw--candidate (that is, a
model of ZFC$^-_*$, suitable version of ZFC, not necessarily well-founded
but with standard $\omega$), and with the real parameters of the
$\varphi_\ell$'s, and $p\in\bQ^N=\{x:N\models\varphi_0(x)\}$, {\em then\/}
for some $q\in\bQ$, we have 
\begin{enumerate}
\item[$(\alpha)$]  $p\le q$,
\item[$(\beta)$]   $q$ is $(N,\bQ)$--generic, which means that for every
${\mathcal I}\in {\rm pd}(N,\bQ)=\{{\mathcal J}\in N:N\models ``{\mathcal
J}$ is predense in $\bQ$''$\}$, for some list $\langle r_n:n<\omega\rangle$
of $\{x:N\models$`` $x \in {\mathcal I}$ ''$\}$ we have $\{r_n:n<\omega\}$
is predense above $q$, 
\item[$(\gamma)$]  moreover $\varphi_2(q,\langle r_n:n<\omega\rangle)$
holds, 
\end{enumerate}
\item[(c)]if $\varphi_2(q',\langle r'_n:n<\omega\rangle)$ then the set
$\{r'_n:n<\omega\}$ is predense above $q'$. 
\end{enumerate}
\item Omitting the ``nw'' or writing just ``w'' means we allow only
well-founded candidates.
\item \relax[On ${\bf Ur}$ see below]. We say $\name{r}\subseteq {\bf Ur}$
(or $\name{r}\subseteq {\mathcal H}_{<\aleph_1}({\bf Ur})$ in the w-case) is
generic for $\bQ$ if: $\name{r}$ is a $\bQ$-name and ``$a\in\name{r}$'' is
determined by the truth value $\varphi_{\name{r},a}[\name{G}_\bQ]$ and the
sequence $\langle\varphi_{\name{r},a}:a\rangle$ is definable in ${\mathfrak
B}_\bQ$ (see \cite{Sh:669}). 
\end{enumerate}
\end{definition}
Now, there are more examples of $\omega$--nw--nep forcing notions in
addition to Cohen and random: say the Sacks forcing notion. But about the
Laver forcing notion we should beware; note that we can guarantee that for
every ${\mathcal I}\in{\rm pd}(\bQ,N)$, for some front $X$ of $q$ we have
$\eta \in X\ \Rightarrow\ q^{[\eta]}$ is above a member of ${\mathcal I}$, 
but being a front is not absolute from $(\bar{\varphi},
\omega)$--nw--candidates, as they are not necessarily well-founded. In fact,
we can easily craft a counterexample: assume $N \models$ ``$\alpha$ a
countable ordinal'', but from the outside not well ordered. There are $f,T
\in N$ (so actually $f=f_\alpha$, $t=T_\alpha$) such that 
\[\begin{array}{ll}
N\models\mbox{``}&T\subseteq {}^{\omega>}\omega\mbox{ is closed under
initial segments, }\quad\langle\rangle\in T,\\  
&f:T\longrightarrow\alpha +1,\qquad f(\langle\rangle)=\alpha,\\
&f(\eta)>0\quad \Rightarrow\quad (\forall n)(\eta\conc\langle n\rangle\in
T),\\
&f(\eta)=0\quad\Rightarrow\quad(\forall n)(\eta\conc\langle n\rangle\notin
T),\\ 
&f(\eta)=\beta+1\quad \Rightarrow\quad (\forall n<\omega)[\eta\conc\langle n
\rangle\in T\ \&\ f(\eta\conc\langle n\rangle)=\beta],\\
&\mbox{if }f(\eta)=\delta\mbox{ is limit}\\
&\mbox{then }\langle f(\eta\conc\langle n\rangle):n<\omega\rangle\mbox{ is
strictly increasing with limit }\delta\mbox{ ''.} 
  \end{array}\]
Let
\[\begin{array}{lcl}
{\mathcal I}_0&=&\{({}^{\omega>}\omega)^{[\eta]}:\eta\in T,\ f(\eta)=0\},\\ 
{\mathcal I}_{n+1}&=&\{({}^{\omega>}\omega)^{[\eta\conc\rho]}:\eta\in
{\mathcal I}_n,\ \rho \in T,\ f(\rho) = 0\}.
  \end{array}\]
Clearly, above no $q$ is every ${\mathcal I}_n$ predense. 

Still,

\begin{theorem}
\label{4.12}
The ${}^\omega \omega$-bounding and almost ${}^\omega\omega$--bounding
forcing notions covered by \cite{RoSh:470} and \cite{RoSh:670} are all
$\omega$--nw--nep. 
\end{theorem}

What about iterations?  Now, for unifying the treatment of finite support
and countable support we revise our definition to have two quasi-orders
$\le^\bQ,\le^\bQ_{\rm pr}$ such that $p\le_{\rm pr}q\ \Rightarrow\ p\le q$.
Hence in Definition \ref{4.11} we add $\varphi_{1,2}$ (absolute just like
$\varphi_1$) serving as a definition of $p\le_{\rm pr} q$. The support is
countable but finite for the apure cases, i.e., only for finitely many
$\alpha$ ,$\neg(\emptyset_\alpha\le_{\rm pr} p(\alpha))$. First assume the 
length is $\alpha^*<\omega_1$, so we can use a parameter coding a well
ordering on $\omega$ with this order type. We should repeat the proof in
\cite{Sh:630} in order to prove preservation in this case, but we better not
use the $L_{\omega_1,\omega}$--completeness as there, as we have problem
with well--foundedness. So we just demand: elements in $\bP_\alpha$ have
depth $<\omega\alpha$ (or so).

What about long iterations? It seems to me, at least now, better (and fit to
\cite{Sh:630}, too) to use a set of urelements ${\bf Ur}$; let ${\mathfrak
C}$ [and ${\mathfrak B}$] be models with universes ${\bf Ur}$ [or $\subseteq
{\bf Ur}$] and ${\mathcal S}\subseteq [{\bf Ur}]^{\le\aleph_0}$ be unbounded
(usually stationary, if $a\in {\mathcal S}$ then ${\mathfrak B}\restriction
a\subseteq {\mathfrak B}$ and ${\mathfrak C}\restriction a \subseteq
{\mathfrak C}$), and anyhow the family of nw--candidates should be
$<(\aleph_1)$--directed and ${\rm uord}\subseteq{\bf Ur}$ is a well ordered 
set (which will serve as the length of the iteration). Now, a candidate will
be a countable model $N \subseteq ({\mathcal H}(\chi),\in)$ of ZFC$^-_*$
(where ${\mathcal H}(\chi)$ includes the urelements), $N\cap{\bf Ur}\in S$,
where ${\bf Ur}$ and the relations of ${\mathfrak B}$ and ${\mathfrak C}$
are the considered relations. We define an nw--candidate similarly but now
$\in^N$ is a relation on $N$ and $N$ is not necessarily well-founded (but
the order type of the well ordered ordinals of $N$ is $>\otp(N\cap{\rm
uord})$.  In Definition \ref{4.11} we demand, of course, that $\varphi_\ell$
are upward absolute from the nw--candidates. Now we can use ${\rm uord}$ as
the index set for iteration (instead of the true ordinals) and there are no
problems. 

This set-up looks like a nice army with no enemy yet, but this seems to me a
natural dividing line among the nep forcing notions and therefore reasonable 
for our interest in those forcing notions per se. I hope it will help,
particularly with \ref{4.8} (and even more so by the c.c.c.~version see
\ref{4.13} below). 

More on preservation (for not necessarily c.c.c.~forcing notions),
commutativity, associativity of generic sets, and countable for pure/finite
for apure support iterations see \cite{Sh:669}.

A restricted version of the large continuum is
\begin{problem}
\label{4.12B}
Can we have long $(>\omega_2)$ iterations of ${}^\omega\omega$--bounding
forcing notions (or at least nw--nep ones) not collapsing cardinals and not
adding Cohen reals? 
\end{problem}
What about the c.c.c.~(nw-nep) ones (or even very Souslin c.c.c.) ones?

\begin{definition}
\label{4.13}
$\bQ$ is c.c.c.--nw--nep {\em if\/} it is a pair $\langle\varphi_0,\varphi_1
\rangle$ of formulas such that 
\begin{enumerate}
\item[$(\alpha)$] $\bQ=\bigcup\{\bQ^N:N$ is an nw-candidate$\}\subseteq
{}^\omega{\rm uord}$, similarly\\
$\le_Q=\bigcup\{\le^N_\bQ:N$ is an nw-candidate$\}$,\\  
recall that in $N$, $\bQ^N,\le^N_\bQ$ are defined by the appropriate
$\varphi_\ell$, 
\item[$(\beta)$] all $\varphi_\ell$ are upward absolute among nw-candidates,
\item[$(\gamma)$] if $N\models$`` ${\mathcal I}\subseteq\bQ$ is predense '',
{\em then\/} ${\mathcal I}^\bQ$ is really predense.
\end{enumerate}
\end{definition}

\begin{proposition}
\label{4.14}
\begin{enumerate}
\item The Cohen and random forcing notions are c.c.c.-{\rm nw}-nep.
\item The class of c.c.c.-{\rm nw}-nep is closed under FS iterations. 
\item This class is also closed under subforcings.
\end{enumerate}
\end{proposition}

And I am curious to know:
\begin{problem}
\label{4.15}
Does \ref{4.14} exhaust all c.c.c-nw-nep forcing notions (at least those
with a generic real)? 
\end{problem}

\[* \qquad * \qquad *\]

Being interested in classifying nep c.c.c.~forcing notions, we may consider
sweetness; the discussion below is in fact an introduction to
\cite{RoSh:672}. Sweetness phenomenons are when we can build homogeneous
forcing notions (as in \cite[\S7,\S8]{Sh:176}), sour phenomenons are strong
negations (as in \cite[\S6]{Sh:176}).

\begin{problem}
\label{4.15a}
\begin{enumerate}
\item For which $(\bQ,\name{\eta})$, nep c.c.c.~forcing notions, is it
consistent that:  
\begin{enumerate}
\item[(a)] there is a $(\bQ,\name{\eta})$--generic real over $\bL[A]$ for
every $A \subseteq\omega_1$, and 
\item[(b)] for every subset $B\in\bL[\bR]$ of ${}^\omega 2$ for some $A
\subseteq\omega$, for a dense set of $p\in (\bQ,\name{\eta})$:\quad for some
truth value ${\bf t}$, {\em if\/} $\eta$ is a $(\bQ^{\bL[A]}_{\ge p},
\name{\eta})$--generic real {\em then\/} $\eta\in B\ \Leftrightarrow\ {\bf
t}$. 
\end{enumerate}
\end{enumerate}
\end{problem}
So for this question, random reals are complicated (see \cite[\S6]{Sh:176}),
whereas Cohen real and universal-meagre one (and dominating = Hechler reals)
are low, see Judah and Shelah \cite{JdSh:446}. 

More generally,

\begin{problem}
\label{4.15b}
Let $(\bQ,\name{\eta})$ be a nep c.c.c.~forcing notion, and $\kappa$ be a
cardinal number. Let $I^\kappa_{(\bQ,\eta)}$ be the $\kappa$--complete ideal
generated by sets of the form 
\[A_N = \{\eta:\eta\mbox{ is not }(\bQ,\name{\eta})\mbox{--generic over
}N\}\]
for countable models $N\prec({\mathcal H}(\chi),\in)$ to which $(\bQ,
\name{\eta})$ belongs.\\ 
What is the consistency strength of `` every projective set, or even every
set from $\bL[\bR]$, is equal to a Borel set modulo $I^\kappa_{(\bQ,
\name{\eta})}$ '' ?
\end{problem}
We hope for a strong dichotomy phenomena, i.e., if the answer above is
negative for $(\bQ,\name{\eta})$ (so the sweetness fails), then a strong
negation holds, so we call such phenomena sourness. 

\begin{definition}
\label{4.15c}
Let $\bQ_1,\bQ_2$ be nep c.c.c.~forcing notions definable in $\bL$ (or
$\bL[r]$). We say that $\bQ_1,\bQ_2$ are explicitly sour over Cohen {\em if}
we can find $\bQ_\ell$-names $\name{\eta}_\ell$ of Cohen reals (for $\ell=
1,2$) such that 
\begin{quotation}
if $G_\ell\subseteq \bQ_\ell$ is generic over $\bL$ (or $\bL[r]$) for
$\ell=1,2$\\
then $\name{\eta}_1[G_1]\ne\name{\eta_2}[G_2]$.
\end{quotation}
\end{definition}
We should note that there may be homogeneity for wrong reasons, i.e., maybe
when we force, very few $(\bQ,\name{\eta})$--generic reals over $\bL$ are
added and then homogeneity holds for ``degenerated'' reasons; we may call
such cases sacharin. For more on this direction see \cite{RoSh:670}.

Speaking about the class of sweet forcing notions we should mention the
following problem.

\begin{problem}
\label{add1}
For any cardinal $\kappa$ and a large cardinal property (or consistency
strength) we may ask the following.
\begin{enumerate}
\item Is there a widest class ${\mathcal K}$ of absolute enough forcing
notions such that for some forcing notion $\bP$ we have
\begin{enumerate}
\item[(a)] $\Vdash_\bP\kappa=\aleph_1$,
\item[(b)] $\bP$ is homogeneous for complete subforcings from ${\mathcal
K}$, moreover
\item[(c)] if $\bP^*\lesdot\bP$, $\bP^*$ has a generic real then
$\bP/\name{G}_{\bP^*}$ satisfies (b) in $\bV[\name{G}_{\bP^*}]$.
\end{enumerate}
\item If not, at least give a wide enough such class. 
\item Are there two classes ${\mathcal K_1},{\mathcal K}_2$ as above, such
that there is no class with the respective property including ${\mathcal
K_1} \cup{\mathcal K}_2$ ? Or even that the consistency strength of the
(now) obvious conclusion is higher than the given one?
\end{enumerate}
\end{problem}
Now, the variants of sweetness try to deal with the case of
$\kappa=\aleph_1$ and the consistency strength ZFC (see Ros{\l}anowski and
Shelah \cite[\S 3]{RoSh:672}); the theory of determinacy is applicable to
the case $\kappa=\aleph_1$ and maximal consistency strength (see Woodin
\cite{Woxx}), and \cite{Sh:295} intends to deal with the case of ``ZFC + 
$\kappa$ is strongly inaccessible'' (and no further consistency strength
assumptions).   

\begin{question}
\label{add2}
Is there a sweet forcing notion (see \cite[\S 7]{Sh:176}), preferably
natural, such that it cannot be completely embedded into the forcing notion
constructed in \cite[\S 7]{Sh:176} (it was gotten by composing ${\mathcal
U}{\mathcal M}$, amalgamating and direct limits), or at least not
``below ${\mathcal U}{\mathcal M}$'' (in the sense of $\leq^*_1$, see
\ref{4.6}) ? 
\end{question}

\[* \qquad * \qquad *\]

Let me now mention free iterations: 

In nep (also above) we can replace CS by free limit as in
\cite[Ch.IX]{Sh:b}. (This was my third proof for preservation of properness
(the first was like \cite[\S2]{Sh:669}, the second was like the one in
\cite[Ch.III]{Sh:b}; this third proof looks very natural but no real reason
for replacing CS iteration by it has appeared).

In particular

\begin{proposition}
\label{4.16}
Definition of $L_{\omega_1,\omega}$--free iteration is absolute enough, so
we have our $\varphi_2$ (we ignore the {\rm nw}).
\end{proposition}
 
Preservation (of reasonable properties) by CS iterations of proper forcing
(or variants) seems to me a worthwhile subject. For ${}^\omega
\omega$--bounding the situation is nice, ``proper $+ {}^\omega
\omega$--bounding" is preserved, and analogous results hold for a large
family of properties even, e.g., ``$D$ generates a $P$-point ultrafilter on
$\omega$''. But some properties do not fit, though still are preserved in
limits (see \cite[Ch.XVIII,\S3]{Sh:f}). For example, we shall consider below
the case of ``$A\subseteq {}^\omega 2$ is non-null ''. The simple
preservation is: in the existence of generic conditions we can preserve `` a
given old $\eta\in {}^\omega 2$ is random over the model $N[\name{G}_{\bP}]$
''.  

\begin{question}
\label{4.17}
Assume that
\begin{enumerate}
\item[(a)]  $\bar{\bQ}$ is a CS iteration, or a $L_{\omega_1,\omega}$--free
iteration, 
\item[(b)]  each $\bQ_\alpha$ is proper, or a nep forcing notion,
\item[(c)]  each $\bQ_\alpha$ is ``non-null for
$(\bS,\name{r})$--preserving''\footnote{$\bS$ is a nep forcing notion,
$\name{r}$ is a hereditarily countable $\bS$--name, so $(\bS,\name{r})$
induces an ideal on ${}^\omega 2$, $\bS$ may be random and then the ideal is
the ideal of null sets.} 
(if each $\bQ_\alpha$ is nep and $\name{r}$ is, e.g., random, then by
\cite{Sh:630} this is equivalent to not making old ${}^\omega 2$ null).
\end{enumerate}
Does $\bP_\alpha$ have the property from clause (c)?  (so we have four
versions of the question, as for clause (a) and for clause (b) we can choose
the first or the second possibilities).
\end{question}

Assume:
\begin{enumerate}
\item[$(*)$]  each $\name{\bQ}_i$ has a generic $\name{r}_i\subseteq
\name{\alpha}_i$. 
\end{enumerate}
(Hence $\bP_\alpha$ have this property as it is preserved).

We assume knowledge of free iterations (see \cite[IX,\S1,\S2]{Sh:f}); in
short, if $\bP_n\lesdot\bP_{n+1}$ for $n<\omega$, let $\bP_\omega$ be
\[\begin{array}{ll}
\{\psi:&\psi\mbox{ is a sentence in the $L_{\omega_1,\omega}$ propositional
calculus with}\\
&\mbox{the set of propositional variables $\bigcup\{\bP_n:n<\omega\}$ such
that in some}\\
&\mbox{forcing extension of $\bV$ there is $G\subseteq
\bigcup\{\bP_n:n<\omega\}$ satisfying}\\
&\mbox{(a)\ for each $n<\omega$, $G\cap \bP_n$ is generic over $\bV$,}\\
&\mbox{(b)\ looking at $G$ as assigning truth values to members of
}\bigcup\{\bP_n:n<\omega\},\\
&\mbox{\qquad it assigns the value truth to }\psi\}.
  \end{array}\]
The order of $\bP_\omega$ is the natural one.

We deal with an $L_{\omega_1,\omega}$--free iteration $\bar{\bQ}$ such that 
\[\Vdash_{\bP_i}\mbox{`` }\name{\bQ}_i\mbox{ is nep and $({}^\omega
2)^{\bV^{\bP_i}}$ not null in }({}^\omega 2)^{\bV^{\bP_i*\name{\bQ}_i}}
\mbox{ ''.}\]
This is quite a wide case. What does it mean for a successor $\alpha= \ell
g(\bar{\bQ})$? Say, $\alpha=2$ so we know that in $\bV^{\bQ_0}$, $({}^\omega
2)^\bV$ is not null, and in $\bV^{\bQ_0*\name{\bQ}_1}$, $({}^\omega
2)^{\bV^{\bQ_0}}$ is not null. But for nep forcing notions preserving the
non-nullity of ${}^\omega 2$ implies preserving the non-nullity of any old
non-null set by \cite[\S7]{Sh:630}. There is no problem with successor
stages. So now assume $\delta=\ell g(\bar{\bQ})$ is a limit ordinal of
countable cofinality.

As $\bP_\delta$ is nep it is enough to assume
\begin{enumerate}
\item[$(**)$]  $p_0\Vdash_{P_\delta}$ ``$\name{T}\subseteq {}^{\omega >}2$
is a subtree and ${\rm Leb}\big(\lim(\name{T})\big)>0$ '', and $\delta=
\bigcup\limits_{n<\omega} i_n$, $i_n<i_{n+1}$, and without loss of
generality $\name{T}\cap {}^{n \ge}2$ is a hereditarily countable
$\bP_{i_n}$--name, 
\end{enumerate}
and to find an old $\eta\in {}^\omega 2$ such that $p_0
\nVdash_{\bP_\delta}$`` $\eta\notin\lim(\name{T})$ ''. Now, let $N$ be a
$\bar{\bQ}$--candidate to which $\{\name{T},p_0,\langle
i_n:n<\omega\rangle\}$ belongs (without loss of generality $N\prec({\mathcal
H}(\chi),\in)$, with $\chi$ large enough). Let $G_0\subseteq{\rm Levy}(
\aleph_0,(2^{|\bP_\delta|})^N)$ be generic over $N$. Let $p_1\in N[G_0]$ be
such that $p_1\in (\bP_\delta)^{N[G_0]}$, $N[G_0]\models$`` $p_0\le p_1$ and
$p_1$ is explicitly $(N,\bP_\delta)$--generic ''. Let $G_1=G_{1,\delta}
\subseteq\bP_\delta^{N[G_0]}$ be generic over $N[G_0]$ such that $p_1\in
G_{1,\delta}$ and let $G_{1,i_n}=G_1\cap\bP^{N[G_0]}_{i_n}$. Let $s$ be a
random real over $N[G_0][G_1]$ (if we replace random by other nep-explicitly
demand even more models), hence over $N$ too. Clearly we can choose such
$G_0,G_{1,\delta}$ in $\bV$. So $N[s]$ is a class (= definable) of
$N[G_0][G_1][s]$, and clearly $N[s]$ is a $\bar{\bQ}$--candidate. Also,
there is $s'\in\lim(\name{T}[G_1])$, $s'\equiv s$ (i.e., $s'\in {}^\omega 2$
and the set $\{\ell<\omega:s'(\ell)\ne s(\ell)\}$ is finite). Let us define
$\psi$ as 
\[\psi = p\ \&\ \bigwedge_{n<\omega}\bigl[s'\restriction n\in (\name{T}\cap
{}^n 2)[\name{G}_{\bP_{i_n}}]\bigr]\]
(note that, by the assumption, there is $\langle\psi_{n,\eta}:n<\omega,\
\eta\in {}^n 2\rangle$, $\psi_{n,\eta}$ an $L_{\omega_1,\omega}$--sentence
for $\bP_{i_n}$, i.e.~using (countably many) variables $q\in\bP_{i_n}$,
such that 
\[p_0\Vdash\mbox{`` }\eta\in\name{T}\cap {}^n 2\ \mbox{ iff }\
\psi_{n,\eta}\in \name{G}_{P_{i_n}}\mbox{ ''}\]
for $n<\omega$, $\eta\in {}^n2$, so, up to equivalence, 
\[\psi=p\ \& \ \bigwedge_n\psi_{n,s'\restriction n}\]
(recall that $\name{T}\cap {}^n 2$ is a $\bP_{i_n}$--name). The problem is
whether $\psi\in\bP_\delta$. Now, $\psi\in N[s]$, so by absoluteness to show
$\psi\in\bP_\delta$ it suffices to show that this holds in $N[s]$.  

So we need
\[N[s]\models\mbox{`` }\psi\ \&\ \bigwedge_n \Psi_{\bP_{i_n}}\mbox{ has a
Boolean-valued model '',}\]
where $\Psi_{\bP_{i_n}}=(\bigwedge\{\bigvee\limits_{q\in {\mathcal I}} q:
{\mathcal I}\subseteq\bP_{i_n}\mbox{ is predense }\})^{N[s]}$. But
$N[G_1][G_2][s]$ is a generic extension of $N[s]$, so it is enough to prove
it there; so there is no problem.  In fact for the case of ``non-null'' the
answer is yes (for CS iterations use \cite[Ch.XVIII,3.8, pp
912--916]{Sh:f}). 

On the other hand, in full generality the answer to \ref{4.17} is no; note
that life is harder if we want to preserve positiveness for $I_{(\bQ,
\name{\eta})}$, where $\bQ$ is (nep but) not c.c.c., on this and more see
\cite{Sh:669}. 

Another possible direction is

\begin{problem}
\label{4.9mv}
\begin{enumerate}
\item Is there an interesting theory of nicely definable forcing notions in
${\mathcal H}(\theta)$ or ${\mathcal H}_{<\theta}(\sigma)$?
\item Similarly for the theory of iterations (see later).
\item Generalize \cite{RoSh:470}, replacing $\aleph_0$ by $\kappa$. Probably
require $\kappa=\kappa^{<\kappa}>\aleph_0$ or even that $\kappa$ is strongly
inaccessible. 
\end{enumerate}
\end{problem}
E.g., let for simplicity $\kappa$ be strongly inaccessible, $D^*$ a normal
filter on $\kappa$. For a cardinal $\theta<\kappa$ let a $\theta$--creature
${\mathfrak c}$ consist of $(R^{\mathfrak c},{\rm pos}^{\mathfrak c},{\rm
val}^{\mathfrak c})$, where:  $R^{\mathfrak c}$ is a $\theta^+$--complete
forcing notion, ${\rm pos}^{\mathfrak c}$ is a non-empty set, ${\rm
val}^{\mathfrak c}$ is a function from $R^{\mathfrak c}$ to ${\mathcal
P}({\rm pos}^{\mathfrak c})\setminus\{\emptyset\}$ such that $R^{\mathfrak
c} \models$``$x\le y$''$\quad\Rightarrow\quad{\rm val}^{\mathfrak c}(y)
\subseteq{\rm val}^{\mathfrak c}(x)$.

A $\kappa$--normed tree is $(T^*,\bar{\mathfrak c},\bar{\theta})$, where
$T^* \subseteq {}^{\kappa >}{\rm Ord}$ is a subtree, usually closed under
increasing sequences of length $<\kappa$, $\bar{\mathfrak c}=\langle
{\mathfrak c}_\eta:\eta\in T^*\rangle$, ${\mathfrak c}_\eta$ is a
$\theta_\eta$--creature, $\bar{\theta}=\langle\theta_\eta:\eta\in T^*\rangle$
and ${\rm pos}^{{\mathfrak c}_\eta}={\rm Suc}_{T^*}(\eta)$.

We can consider
\[\begin{array}{rl}
\bQ^1_{(T^*,\bar{\mathfrak c}^*,\bar{\theta}),D^*}=\quad\ &\\
\bQ^1=\big\{(T,\bar{r}):\mbox{(a)}&T\mbox{ is a subtree of }T^*,\
(<\kappa)\mbox{-closed},\\ 
\mbox{(b)}&\bar{r}=\langle r_\eta:\eta\in {\rm sp}(T))\rangle,\ r_\eta\in
R^{{\mathfrak c}_\eta}\mbox{ and }{\rm val}^{{\mathfrak c}_\eta}(r_\eta)=
{\rm Suc}_T(\eta),\\
\mbox{(c)}&\mbox{for every }\eta\in\lim_\kappa(T), \mbox{ we have }
\{\zeta{<}\kappa\!:\eta\restriction\zeta\in {\rm sp}(T)\}\in D^*\big\}
\end{array}\]
with the natural order, where ${\rm sp}(T)=\{\eta\in T:(\exists^{\ge 2} x)
[\eta\conc\langle x \rangle \in T]\}$. Trivially, forcing with $\bQ^1$ adds
neither bounded subsets of $\kappa$ nor sequences of ordinals of length
$<\kappa$. For this forcing notion, if $T^*_\alpha=:\{\eta\in T^*:\ell
g(\eta)=\alpha\}$ has cardinality $<\kappa$ for each $\alpha<\kappa$, and
$(\forall \eta\in\lim_\kappa(T^*))[\{\zeta<\kappa:\theta_{\eta\restriction
\zeta}>\zeta\}\in D^*]$ (or there is $A\in D^*$ such that $\zeta\in A\ \&\
\eta\in T\ \&\ \ell g(\eta)=\zeta\quad\Rightarrow\quad\theta_\eta>\zeta$),
{\em then\/} forcing with $\bQ^1_{(T^*,\bar{\mathfrak c}^*,\bar{\theta}),
D^*}$ does not collapse $\kappa^+$. Also, this forcing is
$(<\kappa)$--complete ${}^\kappa \kappa$--bounding and $\bQ^1$--names of
$\name{\tau}:\kappa\longrightarrow {\rm Ord}$ can be read continuously on
$T$ for a dense set of $(T,\bar{r})\in\bQ^1$. Moreover, this property is
preserved by $(\le\kappa)$--support iterations (by \cite[\S1]{Sh:655}, or
directly). 

We can allow gluing (i.e., putting together with $\eta$ many nodes above,
creating a new forcing notion, i.e.~creatures, see \cite[3.3(2),\S
6.3]{RoSh:470}).  

We may consider $\bar{\mathfrak c}=\langle\bar{\mathfrak c}_\zeta:\zeta<
\kappa\rangle$, ${\mathfrak c}_\zeta$ a $\theta_\zeta$--creature and for
some club $E$ of $\kappa$, $\zeta\in E\ \&\ \xi\ge\zeta\quad\Rightarrow\quad
\theta_\xi>\zeta$, and consider
\[\bQ^2_{\bar{\mathfrak c}}=\big\{(w,\bar{r}):\mbox{for some }\zeta,\ w\in
\prod_{\varepsilon<\zeta}{\rm pos}^{{\mathfrak c}_\varepsilon},\ \bar{r}=
\langle r_\varepsilon:\varepsilon\in [\zeta,\kappa]\rangle,\ r_\varepsilon
\in R^{{\mathfrak c}_\varepsilon}\big\}\]
with the natural order. If $\zeta<\kappa\ \Rightarrow\ |R^{{\mathfrak
c}_\zeta}|+|{\rm pos}^{{\mathfrak c}_\zeta}|<\kappa$ and
\begin{enumerate} 
\item[$(\alpha)$] $\kappa=\sup\{\zeta:\mbox{for }\xi\in [\zeta,\kappa),\
R^{{\mathfrak c}_\xi}\mbox{ is }|\prod\limits_{\varepsilon<\zeta}{\rm
pos}^{{\mathfrak c}_\varepsilon}|^+\mbox{--complete}\}$, or 
\item[$(\beta)$]  $\diamondsuit_\kappa$
\end{enumerate}
then forcing does not collapse $\kappa^+$ and is $(< \kappa)$-complete but
the ``read continuously'' is problematic for case $(\beta)$.

Again, we may allow the forcing notion to be omittory (see
\cite[2.1.1]{RoSh:470}), i.e., allow 
\[\begin{array}{ll}
\bQ^3_{\bar{\mathfrak c}}=\{(w,\bar{r}):&\mbox{for some bounded }u\subseteq 
\kappa\mbox{ and } A\in D,\\
& w\in\prod\limits_{\varepsilon\in u}{\rm pos}^{{\mathfrak c}_\varepsilon},
\mbox{ and }\bar{r}=\langle r_\alpha:\alpha \in A\rangle\}
  \end{array}\]
and/or a combination of creatures (i.e., the function $\Sigma$) and/or we
may allow memory (=the object which ${\mathfrak c}_\alpha$ produces depends
on the earlier ${\mathfrak c}_\beta$'s) and/or gluing.   

Note that natural nice enough ${\mathfrak c}$'s make us regain ``read
continuously'' and its parallels.

All those generalize naturally. But if in $\bQ^1_{(T^*,\bar{\mathfrak c},
\theta),D^*}$ we allow $D^*$ to be the co-bounded filter, the proofs fail. 
There is more to be said and done.

Of course, we can now carry out generalizations of various independence
results on cardinal invariants, e.g.~between variants of $(f,g)$--bounding
and slaloms (= corsets), see \cite{Sh:326}, \cite[Ch.VI,\S2]{Sh:f},
Goldstern and Shelah \cite{GoSh:448}, Ros{\l}anowski and Shelah
\cite{RoSh:670}. We may consider in the tree version for every limit
$\delta\in S$ ($S\subseteq\kappa$ is a stationary set) to omit tops for one
or just $\le|\delta|$ branches in a condition $p=(T,\bar{r})$, provided for 
each $\eta\in T$ we have $\langle |\{\rho:\eta\vartriangleleft\rho\in T,\ 
\ell g(\rho)=\ell g(\eta)+\alpha\}|:\alpha<\kappa\rangle$ goes to $\kappa$. 

Recall that whereas Cohen forcing and many others have $\kappa$-parallels,
not so with random real forcing.

\begin{problem}
\label{4.9mu}
\begin{enumerate}
\item Prove that there is no reasonable parallel, say no
$\kappa^+$--c.c.~forcing notion such that any new member of ${}^\kappa
\kappa$ is bounded by some old one.
\item Similarly as far as the parallel of \cite{Sh:480} is concerned.
\end{enumerate}
\end{problem}

Concerning nep forcing notions, we may define when ``$\bQ$ is semi-nep'' as
follows:  for some nep forcing notion $\bQ^*$ we have 
\begin{enumerate}
\item[(a)] $\bQ\subseteq \bQ^*$ (that is the set of elements of $\bQ$ is a
subset of the set of elements of $\bQ^*$ and the order of $\bQ$ is the order of
$\bQ^*$ restricted to $\bQ$),
\item[(b)] if $N$ is a $\bQ$--candidate, $\bS\in N$, $\bS^N=\bQ^*\cap N$, $p
\in\bS$, $p\in N$ then for some $q$, $p\leq q$, $q$ is $(N,\bS)$--generic
(in an explicit way). 
\end{enumerate}

\section{To prove or to force, this is the question}
On many things we confidently ``know'' that they are independent, ``just'' a
proof is needed (for many others we know that forcing will not help by
absoluteness). The rest we may think are decidable, but actually we are not
sure. More fall in the middle; our intuitions do not give an answer or
worse: they give an answer which oscillates in time. In the mid-seventies, I
was interested in (see Abraham and Shelah \cite{AbSh:114}):

\begin{question}
\label{5.1}
Is there an Aronszajn tree $T$ and a function $c:T\longrightarrow\{{\rm
red,\ green}\}$ (or more colours) such that for any uncountable set
$A\subseteq T$, all colors appear on the set $\{\eta\cap\nu:\eta,\nu\in A
\mbox{ are incomparable}\}$, where $\eta\cap\nu$ is the maximal common lower
bound?  
\end{question}
Why?  Baumgartner \cite{B4} proved that among uncountable real linear orders
(e.g., with density $\aleph_0$) there may be a minimal one under
embeddability. This follows from
\[\mbox{ CON(if }A,B\in [\bR]^{\aleph_1}\mbox{ are $\aleph_1$--dense, then
they are isomorphic).}\]
So we may ask: among uncountable linear orders, can there be finitely many
such that any other embeds one of them? (call such a family a base).
 
A base should contain a real order and $\omega_1$ and $\omega^*_1$. Any
linear order into which none of them embeds is necessarily a Specker order
(= take an Aronszajn tree, order it lexicographically).  

You may ask: Can there be a ``minimal'' order among those? But there
cannot. It is known (\cite{Sh:50}, answering a question of Countryman) that
there is a Specker order $L$ such that the product $L\times L$ (with the
product order) is the union of countably many chains (comes from a very special
Aronszajn tree). Hence $L$ and $L^*$ (its inverse) embed no common
uncountable chain. 

Now, consistently any two candidates for $L$ are isomorphic or
anti-isomorphic. So, if we can also put the additional forcing together,
then we will have a candidate for a basis which seems extremely likely. We
just need to start with a Specker order, i.e., an Aronszajn tree with a
lexicographic order, under the circumstances (okay to force a little),
without loss of generality with any node having two immediate successors.
Look at it as an Aronszajn tree, and remember Abraham and Shelah
\cite{AbSh:114}. So  without loss of generality, on a club the tree is
isomorphic to the one for $L$, and let $\Dom(c)=T$,
\[c(\eta)=\left\{
\begin{array}{rl}
\mbox{ green}&\mbox{ if the two linear orders make the same}\\
&\quad\mbox{ decision about the two immediate successors of }\eta,\\
\mbox{ red}&\mbox{ otherwise,}
  \end{array}\right.\]
and so the problem \ref{5.1} arises.

\[* \qquad * \qquad *\]

A common property of some of the problems discussed below (\ref{5.2},
\ref{5.3A}, \ref{5.4}) is a difference between asking about
$S^\lambda_\theta = \{\delta<\lambda:\cf(\delta)=\theta\}$ and asking about
a stationary $S\subseteq S^\lambda_\theta$ such that $S^\lambda_\theta
\setminus S$ is stationary too; a difference of which I became aware in
\cite{Sh:64} (e.g. $\diamondsuit_{S_1}\ \&\ \neg\diamondsuit_{S_2}$ is
possible for disjoint stationary subsets of $\omega_1$) after much agony.

\begin{question}
[GCH]  
\label{5.2}
If $\mu$ is singular, do we have $\diamondsuit_{S^{\mu^+}_{\cf(\mu)}}$?  
(Those are the only cases left.)
\end{question}
Similarly for inaccessibles, see \cite{Sh:587}. 

If we try to force consistency of the negation, note that (for $\mu$ strong
limit singular) 
\[2^\mu=\mu^+ +\square_{\mu^+}\quad\Rightarrow\quad \diamondsuit_{S^{
\mu^+}_{\cf(\mu)}}\qquad\qquad \mbox{ when } \cf(\mu)>\aleph_0\]
(see \cite[3.2,p.1030]{Sh:186}). So we need large cardinals (hardly
surprising for successors of singular cardinals). See more in D\v{z}amonja and
Shelah \cite{DjSh:545}).

Probably it is wiser to try to force this for ``large'' $\mu$. Changing the
cofinality of a supercompact cardinal $\mu$ to $\aleph_0$, where ``$\mu$ is 
prepared'' is not helpful as after the forcing some old
$S^{\mu^+}_{\aleph_0}$ is added to the old $S^{\mu^+}_\mu$ to make the new
$S^{\mu^+}_{\aleph_0}$. 

If $\bV\models\neg\diamondsuit_{S^{\mu^+}_\mu}$ (which may be forced but
you need enough indestructibility of measure), you get $\neg\diamondsuit_S$ 
for some non-reflecting stationary $S\subseteq S^{\mu^+}_{\aleph_0}$, but we
have it ``cheaply'' by forcing (say, starting with $\bL$; see \cite{Sh:186},
better \cite{Sh:667}). So maybe it is wiser to start with $\mu$ of
cofinality $\aleph_0$. 

Let $\mu$ be a limit of large cardinals and try to add enough subsets to
$\mu$ to ``kill'' $\diamondsuit^{\mu^+}_{\cf(\mu)}$.  Our knowledge of such
forcing for such cases is limited at the present. But ZFC + GCH still give
an approximation (see more \cite{Sh:667}): 
\begin{enumerate}
\item[$(*)$] if $S\subseteq S^{\mu^+}_{\cf(\mu)}$ is stationary then for
some $\langle\langle\alpha_{\delta,i}:i<\cf(\mu)\rangle:\delta\in S\rangle$
we have:  
\begin{enumerate}
\item[(a)] $\alpha_{\delta,i}$ is increasing with limit $\delta$,
\item[(b)] if $\theta<\mu$, $f:\mu^+\longrightarrow\theta$,\quad {\em then}
\[(\exists\mbox{ stationarily many }\delta\in S)(\forall i<\cf(\mu))(f(
\alpha_{\delta,2i})=f(\alpha_{\delta,2i+1})).\] 
\end{enumerate}
\end{enumerate}
Note: Having two equal values inside a group calls for division (or
subtraction) so that we get a known value. So, if we are trying to guess
homomorphisms from a group $G$ with $|G|=\mu^+$ to $H$, $|H|=\theta<\mu$,
$G=\bigcup\limits_{i<\mu^+}G_i$, $G_i$ strictly increasing continuous,
$|G_i|<\mu^+$, and $S\subseteq S^{\mu^+}_{\cf(\mu)}$ is stationary, then we
can find $\bar{\eta}=\langle\langle g_{\delta,i}:i< \cf(\mu)\rangle:\delta
\in S\rangle$ such that $g_{\delta,i}\in G_{\alpha_{\delta,2i+2}}\setminus
G_{\alpha_{\delta,2i}}$ and for every homomorphism $h:G \longrightarrow H$
there are stationarily many $\delta$ such that $(\bigwedge\limits_i
h(g_{\delta,i})=e_H)$ (relevant to Whitehead groups). Without loss of
generality, $|G_{i+1}\setminus G_i|=\mu$, $|G_0|=\mu$, the universe of $G_i$
is $\mu\times (1+i)$ (or $g_i\in G_{i+1}\setminus G_i$ uses the question on
$f:f(i)=f(g_i)$). See \cite{Sh:667}. 

\noindent{\sc Question:}\quad Can we have something similar for any sequence
$\bar{\eta}$?
 
\noindent{\sc Answer:}\quad  No. We have quite a bit of freedom (e.g.,
demand $\alpha_{\gamma,i}\in A^*$, $A^*\in [\mu^+]^{\mu^+}$ fixed) but 
certainly not for any.

In fact, for any stationary non-reflecting set $S\subseteq
S^{\mu^+}_{\cf(\mu)}$ and any sequence $\bar{\eta}=\langle\langle
\alpha_{\delta,i}: i<\cf(\mu)\rangle:\delta\in S\rangle$ with
$\alpha_{\delta,i}$ increasing with limit $\delta$ we can force:
\begin{quotation}
for every $\langle h_\delta\in {}^{\cf(\mu)}\theta:\delta\in S\rangle$,
$\theta<\mu$, there is $h\in {}^{(\mu^+)}\theta$ such that 
\[(\forall\delta\in S)(\forall^\infty i)(h(\alpha_{\delta,i})=
h_\delta(i)).\] 
\end{quotation}
This is a strong negation of the earlier statements (see \cite{Sh:186},
\cite{Sh:667}). A related ZFC result is that for a singular cardinal $\mu$,
the restriction of the club filter ${\mathcal D}_{\mu^+}\restriction
S^{\mu^+}_{\cf(\mu)}$ is not $\mu^{++}$--saturated (see Gitik and Shelah
\cite{GiSh:577}). 

A well known problem is 

\begin{question}
\label{5.3A}
\begin{enumerate}
\item For a regular cardinal $\lambda>\aleph_2$, can ${\mathcal
D}_{\lambda^+} \restriction S^{\lambda^+}_\lambda$ be
$\lambda^{++}$--saturated? 
\item Similarly just adding the assumption GCH.
\end{enumerate}
\end{question}

\[* \qquad * \qquad *\]

Many ``club guessings'' are true (see \cite{Sh:E12}), but I have looked in
vain several times on: 

\begin{question}
\label{5.4}
Let $\lambda$ be a regular uncountable cardinal. Can we find a sequence
$\langle\langle\alpha_{\delta,i}:i<\lambda\rangle:\delta\in
S^{\lambda^+}_\lambda\rangle$ such that $\alpha_{\delta,i}$ are increasing
continuous in $i$ with limit $\delta$, and for every club $E$ of $\lambda^+$
\[(\exists\mbox{ stationarily many }\delta)(\exists\mbox{ stationarily many
}i)(\alpha_{\delta,i+1},\alpha_{\delta,i+2}\in E),\] 
or other variants (just $\alpha_{\delta,i+1}\in E$ is provable, $\alpha_{
\delta,i}\in E$ is trivial under the circumstances). 
\end{question}
This is interesting even under GCH, particularly as by Kojman and Shelah
\cite{KjSh:447} (essentially) we get from it that there is a
$\lambda^+$-Souslin tree. 

We may think instead of trying to prove for $S,S_1\subseteq
S^{\lambda^+}_\lambda$ being stationary disjoint, that we can force the
failure for $S$ (with GCH). This works (see \cite{Sh:587}) but $S=
S^{\lambda^+}_\lambda$ is harder. The present forcing proofs fail, but also
using ``first counterexample'' fails. We may consider proving:\quad ${\rm
GCH}\ \Rightarrow\ \neg{\rm GSH}$ (where GSH is Generalized Souslin
Hypothesis). Let us look at two successor cases $\lambda^+,\lambda^{++}$
($\lambda$ regular). How can this help?  Assume that there is no
$\lambda^{++}$--Souslin tree and GCH holds. It follows that every stationary
$S \subseteq S^{\lambda^{++}}_{\le\lambda}$ reflects in
$S^{\lambda^{++}}_{\lambda^+}$ (see Gregory \cite{Gre}), moreover it is
enough to assume just that there is no $(<\lambda^+$)--complete Souslin
tree, by Kojman and Shelah \cite{KjSh:447}. Hence
\begin{enumerate}
\item[$(\otimes)$] there is $S^*\subseteq\lambda^+$ such that:
\begin{enumerate}
\item[(a)]  (square on $S^*$)\\
$\bar{C}=\langle C_\alpha:\alpha\in S^*\rangle$, $C_\alpha\subseteq\alpha$
closed, $\otp(C_\alpha)\le\lambda$,\\ 
if $\alpha$ is limit then $\sup(C_\alpha)=\alpha$, and\\
$\beta\in C_\alpha\ \Rightarrow\ C_\beta=\beta\cap C_\alpha$,
\item[(b)] $S^*\cap S^{\lambda^+}_\lambda$ is stationary.
\end{enumerate}
\end{enumerate}
(Why?  This follows by \cite[\S4]{Sh:351} which says that, e.g.,
$S^{\lambda^{++}}_{< \lambda^+}$ is the union of $\lambda^+$ sets with
squares.) This looks as if it should help, but I did not find yet how.

D\v{z}amonja and Shelah \cite{DjSh:545} introduced

\begin{definition}
\label{5.3alpha} 
We say that $\lambda$ strongly reflects at $\theta$ if $\theta<\lambda$ are
regular uncountable cardinals and for some $F:\lambda\longrightarrow \theta$
for every $\delta\in S^\lambda_\theta$ for some club $C$ of $\delta$, $F
\restriction\delta$ is strictly increasing (equivalently, is one--to--one).
\end{definition}
This helps to prove variants of $\clubsuit$ on the critical stationary
subset of $\mu^+$ when $\mu$ is singular, i.e., on $S^{\mu^+}_{\cf(\mu)}$,
see D\v{z}amonja and Shelah \cite{DjSh:545}, \cite{DjSh:562}, and on
independence results Cummings, D\v{z}amonja and Shelah \cite{CDSh:571} and
D\v{z}amonja and Shelah \cite{DjSh:691}.

\begin{question}
\label{5.3beta}
Can we get something parallel when $\cf(\mu)=\aleph_0$?
\end{question}

\begin{question}
\label{5.3gamma}
Can we prove that for some strong limit singular cardinal $\mu$ and a
regular cardinal $\theta<\mu$ we have $\clubsuit_S$, where
$S=S^{\mu^+}_\theta$? 
\end{question}

\[* \qquad * \qquad *\]

On $I[\lambda]$ see \cite{Sh:108}, \cite{Sh:88a}, \cite{Sh:420}. We know
that e.g. 
\[\{\delta<\aleph_{\omega+1}:\cf(\delta)=\aleph_1\}\notin I[\aleph_{\omega
+1}]\] 
is consistent with GCH, but

\begin{problem} 
\label{5.4A}
\begin{enumerate}
\item Can $\{\delta<\aleph_{\omega +1}:\cf(\delta)=\aleph_2\}\notin
I[\aleph_{\omega +1}]$?  
\item Can $\{\delta<(2^{\aleph_0})^{+\omega +1}:\cf(\delta)=(2^{\aleph_0})^+
\}\notin I[(2^{\aleph_0})^{+\omega +1}]$?
\end{enumerate}
\end{problem}
Now \cite[\S4]{Sh:351}, dealing with successors of regulars, raises the
question 

\begin{question}
\label{5.4B}
\begin{enumerate}
\item Let $\lambda$ be inaccessible $>\aleph_0$. Is $I^{\rm sq}_\lambda$
(see Definition \ref{5.4C} below) non-trivial, i.e., does it include
stationary sets of cofinality $\sigma\in (\aleph_0,\lambda)\cap\Reg$? Does
it include such $S$ which is large in some sense (e.g., for every such
$\sigma$)?  
\item Similarly for successor of singular.
\end{enumerate}
\end{question}

\begin{definition}
\label{5.4C}
For a regular cardinal $\lambda>\aleph_0$ let
\[\begin{array}{ll}
I^{\rm sq}_\lambda=\big\{A\subseteq\lambda:&\mbox{for some partial square }
\bar{C}=\langle C_\delta:\delta\in S^1\rangle,\\ 
&S^1 \subseteq \lambda \mbox{ and the set } A\setminus S^1 \mbox{ is not
stationary in }\lambda\big\}.
  \end{array}\]
\end{definition}

\[* \qquad * \qquad *\]

\begin{definition}
\label{5.5A}
A linear order $I$ is $\mu$--entangled if: for any pairwise distinct
$t_{i,\ell}\in I$, $i<\mu$, $\ell<n$, for any $w\subseteq \{0,\ldots, n-1\}$
there are $i<j$ such that 
\[t_{i,\ell}<t_{j,\ell}\quad\Leftrightarrow\quad\ell\in w.\]
If $|I|=\mu$ then we omit $\mu$.
\end{definition}

\begin{question}
\label{5.5} 
Is there an entangled linear order of cardinality $\lambda^+$, where
$\lambda=\lambda^{\aleph_0}$? 
\end{question}
A ``yes'' answer will solve a problem of Monk \cite{M2} on the spread of
ultraproducts of Boolean Algebras; see \cite{Sh:462}.

With the help of pcf we can build entangled linear orders in $\lambda^+$ for
many $\lambda$ which means: provably for a proper class of $\lambda$'s.

The interesting phenomenon is: from instances of GCH, we can give a positive
answer, and if $\mu^{\aleph_0}=\mu$, $2^\mu>\aleph_{\mu^{+4}}$ we have: for
many $\delta<\mu^{+4}$ we get a positive answer to \ref{5.5} with $\lambda=
\aleph_\delta$.  On the other hand, from $\mu=\mu^{\aleph_0}$, $2^\mu<\mu^{+
\omega}$ we can also prove a positive answer.

In fact, in the remaining case there are quite heavy restrictions on pcf. A
typical universe with negative answer to \ref{5.5} (we think that it) will
satisfy: for a strong limit cardinal $\mu$, $2^\mu=\mu^{+\omega +1}$, and for
${\mathfrak a}\subseteq\Reg\cap\mu$, $\mu = \sup({\mathfrak a})$,
$\pcf({\mathfrak a})\setminus \mu$ essentially concentrates on
$\mu^{+\omega+1}$ (say $\langle \mu_i:i<\cf(\mu)\rangle$ increasing
continuous, if ${\mathfrak a}$ is disjoint from $\{\mu^{+n}_i:i<\cf(\mu),\
0<n<\omega\}$ then $\emptyset =\pcf({\mathfrak a})\cap (\mu,\mu^{+
\omega})$.) See \cite{Sh:462}. 

Maybe our knowledge of forcing will advance. Note that we need not only to
have pcf structure as indicated, but also to take care of the non-pcf
phenomena as well for constructing entangled linear order as in \ref{5.5}. 

Considering a ZFC proof of existence, it seems most reasonable to assume
toward a contradiction that the answer is no and consider strong limit
singular $\mu$ of uncountable cofinality. So we know $2^\mu<\aleph_{\mu^{
+4}}$ and being more careful even $2^\mu<\aleph_{\mu^+}$. Let $\gamma(*)=
\min\{\gamma:2^{(\mu^{+ \gamma})}>2^\mu\}$, so necessarily $\gamma(*)$ is a
successor ordinal, say $\gamma(*)=\beta(*)+1$. Let $\lambda=:\mu^{+\beta
(*)}$. We may consider trying to construct an entangled linear order of
cardinality $(2^\mu)^+$, using the weak diamond on $\lambda^+=\mu^{+\beta(*)
+1}$. Moreover, we know that there are trees $T$ with $\lambda^+$ levels and
$\le \lambda^+$ nodes and at least $(2^\mu)^+$ many $\lambda^+$--branches
(even ${}^{\lambda^+}2=\bigcup\limits_{\zeta<2^\mu}\lim_{\lambda^+}(
T_\zeta)$ for some subtrees $T_\zeta$ of ${}^{\lambda^+>}2$, $|T_\zeta|=
\lambda^+$ above).

Moreover, a relative of $\diamondsuit^*_{\lambda^+}$ holds. All this seems 
reasonably promising, but has failed so far to solve the problem.

I have also considered to repeat the proof of the weak diamond for
$\lambda^+$ to try to show that a tree with infinite splitting in the above 
representation is necessary.

\begin{problem}
\label{5.5B}
Can we prove that a stronger version of the weak diamond holds for some
$\lambda$?  E.g., a version with more than two colours and/or fixing the
cofinality. We shall be glad to get even just the definable weak
diamond. See \cite{Sh:208}, \cite[\S3]{Sh:576}, \cite{Sh:638} and
\cite{Sh:603}.  
\end{problem}

\[* \qquad * \qquad *\]

Our ignorance about such problems may well come from our gaps in forcing
theory.  

A major problem (more exactly a series of problems) is 

\begin{problem}
\label{5.11}
\begin{enumerate}
\item Can we have a reasonable theory of iterations (and/or forcing axioms)
for $(<\lambda)$-complete forcing notions ($\lambda=\lambda^{<\lambda}$)?
\item Similarly for forcing notions not changing cofinalities of cardinals
$<\lambda$? 
\item Similarly for forcing notions preserving $\mu^+$ and not adding
bounded subsets to $\mu$, $\mu$ a strong limit singular cardinal?
\end{enumerate}
See some recent information on the first in \cite{Sh:587}, \cite{Sh:655},
\cite{Sh:667}, and even much less on the second \cite[Ch.XIV]{Sh:f}, and on 
the third Mekler and Shelah \cite{MkSh:274}, D\v{z}amonja and Shelah
\cite{DjSh:659}. 
\end{problem}
Though much was done on forcing for the function $2^\lambda$ and some
specific problems, our flexibility is not as good as for $2^{\aleph_0}$ in
forcing theory. 

Particularly intriguing are solutions where we know some $\lambda$ exists
but do not know which.  The dual problem is iterated forcing of length Ord
(class forcing); now for such iteration it is particularly hard to control
in the  neighborhood of singulars.
 
\begin{problem}
\label{5.12}
Prove the consistency of:\quad for every $\lambda$ (or regular $\lambda$) a
suitable forcing axiom holds.
\end{problem}
Relevant is ``GCH fails everywhere'' (see Foreman and Woodin \cite{FW}). Now
Cummings and Shelah \cite{CuSh:530}, \cite{CuSh:541} is a modest try and
\ref{1.18A} is relevant.

Specific well known targets are 

\begin{problem}
\label{5.13}
Is GSH consistent?  (GSH is the generalized Souslin hypothesis: for every
regular uncountable $\lambda$ there is no $\lambda$-Souslin tree.)
\end{problem}

\begin{problem}
\label{5.14} 
Is it consistent that for no regular $\lambda>\aleph_1$ do we have a
$\lambda$-Aronszajn tree (see Abraham \cite{Ab83}, Cummings and Foreman
\cite{CuFo98}). 
\end{problem}

A relevant problem is \ref{6.8}.

\[* \qquad * \qquad *\]

I have found partition theorems on trees with $\omega$ levels very useful
and interesting (see Rubin and Shelah \cite{RuSh:117}, and \cite{Sh:136},
\cite{Sh:b}, \cite[X, XI, XV,2.6]{Sh:f}). In \cite[13,p.1453]{Sh:262} and
\cite[Ch.VIII,\S1]{Sh:a} trying to prove a theorem on the number of
non-isomorphic models of a pseudo elementary class we arrived at the
following problem [without loss of generality, try with $2^\lambda=
\lambda^+$ and see \cite{Sh:620}, by absoluteness]:

\begin{question}
\label{5.6}
Assume $m(*)<\omega$, $2^{\lambda_n}<\lambda_{n+1}$, $M$ is a model with
vocabulary of cardinality $\theta$, $\theta+\mu<\lambda_0$, $a^i_\eta\in M$
for $i<\mu$, $\eta\in T=\bigcup\limits_n\prod\limits_{\ell<n}\lambda_\ell$. 
Can we find a strictly increasing function $h:\omega\longrightarrow \omega$
and one-to-one functions $f^i_n:\prod\limits_{k<n}\lambda_k\longrightarrow
\prod\limits_{k<h(n)}\lambda_k$ such that
\begin{enumerate}
\item[(a)] for $n<m$, $\eta\in\prod\limits_{\ell<m}\lambda_\ell$ we have
$f^i_n(\eta\restriction\ell)=(f^i_m(\eta))\restriction h(\ell)$, 
\item[(b)] for $n<\omega$, $m(*)<\omega$, $i_0<\ldots<i_{m(*)-1}<\mu$ and
$\eta_\ell,\nu_\ell\in\prod\limits_{k<n}\lambda_k$ for $\ell<m(*)$, 
the tuples $\langle a^{i_0}_{\eta_0},a^{i_1}_{\eta_1},\ldots,a^{i_{m(*)
-1}}_{\eta_{m(*)-1}}\rangle, \langle a^{i_0}_{\nu_0},a^{i_1}_{\nu_1},\ldots,
a^{i_{m(*)-1}}_{\nu_{m(*)-1}}\rangle$ realize the same type in $M$ ?
\end{enumerate}
\end{question}
(If on $\lambda_n$ there is a ``large ideal'' (see \cite{Sh:542}) life is
easier, see \cite{Sh:262}.)

\section{Boolean Algebras and iterated forcing}
We turn to Boolean Algebras. Monk has made extensive lists of problems about
Boolean Algebras (which inspired more than few works of mine). His problems
mostly go systematically over all possible relations; our perspective is
somewhat different. 

Among my results on Boolean Algebras I like \ref{6.5A} stated below (see 
\cite{Sh:233}), but the result did not draw much attention though the paper
was noticed (see Bonnet and Monk \cite{BoMo89}, Juh\'asz \cite{Juh92}).

\begin{theorem}
\label{6.5A}
If $B$ is a Boolean Algebra of cardinality $\ge\beth_\omega$ and
$\lambda={\rm Id}(B)$ (the number of ideals of the Boolean Algebra) {\em
then\/} $\lambda=\lambda^{\beth_\omega}$.\\
(We can instead of ideals of Boolean Algebras speak about open subsets of a
compact Hausdorff topology, and we can replace $\beth_\omega$ by any
singular strong limit). 
\end{theorem}

So we are left with

\begin{question}
\label{6.6}
Is it true that for any large enough Boolean Algebra $B$ we have
${\rm Id}(B)={\rm Id}(B)^{< \theta}$ when, e.g., $\theta=\log_2(|B|)$, or at
least for some constant $n$, $\theta=\min\{\mu:\beth_n(\mu)\ge |B|\}$?

(Similarly for compact spaces).
\end{question}
By \cite{Sh:233}, for every $B$ there is such $n$. Of course, in
non-specially constructed universes the answer is yes. If you like to try
consistency, you have to use the phenomena proved consistent in Gitik and
Shelah \cite{GiSh:344}.

On the other hand, a ZFC proof may go in a different way than \cite{Sh:233}.
Related (see Juh\'asz \cite{Ju}) is

\begin{question}
\label{6.7}
What can be the number of open sets of a $T_2$ topology?  $T_3$ topology?
One with clopen basis? 
\end{question}
It seems interesting to consider the following 

\begin{problem}
\label{6.8}
Is there a class of cardinals $\lambda$ (or just two) such that there is a
$(\lambda^+,\lambda)$--thin tall superatomic Boolean Algebra $B$ (i.e.,
$|B|=\lambda^+$, $B$ is superatomic and for every $\alpha<\lambda^+$, $B$
has $\le\lambda$ atoms of order $\alpha$), provably in ZFC? 
\end{problem}

It is well known that if $\lambda=\lambda^{<\lambda}$ then there is a
$(\lambda^+,\lambda)$--thin tall superatomic Boolean Algebra, so for
$\lambda = \aleph_0$ there is one, so for negative consistency we need ``GCH
fails everywhere or at least for every large enough $\lambda$''. Also note
that trivially if there is a $\lambda^+$-tree (i.e., one with $\lambda^+$
levels each of cardinality $\le \lambda$), then there is such $(\lambda^+,
\lambda)$--thin tall superatomic Boolean Algebra. 

The point is that for several problems in Monk \cite{M2}: Problems 72, 74,
75 and ZFC versions of Problems 73, 77, 78, 79 (all solved in Ros{l}anowski
and Shelah \cite{RoSh:599} in the original version, i.e. showing
consistency) there is no point to try to get positive answers as long as we
do not know it for \ref{6.8}.

Also for several problems of \cite{M2} (Problems 49, 57, 58, 61, 63, 87)
there is no point to try to get consistency of non-existence as long as we
have not proved the consistency of the GSH (generalized Souslin hypothesis)
which says there are no $\lambda^+$-Souslin trees or there is no
$\lambda$-Souslin tree for $\lambda=\cf(\lambda)>\aleph_0$. For some others
this is not provable, but it still seems very advisable to wait for the
resolution of GSH.

\begin{problem}
\label{6.9}
Usually the question on cardinal invariants $\inv_1$, $\inv_2$ is ``do we
always have $\inv_1(B)\le\inv_2(B)$?'', or ``do we always have
$2^{\inv_1(B)} \ge\inv_2(B)$?'' But maybe there are relations like
$\inv_1(B) \le (\inv_2(B))^{+n}$ for some fixed $n<\omega$ (or for
$\omega$).  A particularly suspicious case is $|B|\le({\rm irr}_n(B))^+$,
where $n \in [2,\omega]$ and 
\[\begin{array}{ll}
{\rm irr}_n(B)=\sup\big\{|X|:&X\subseteq B\mbox{ is } n\mbox{-irredundant
which means:}\\
&\mbox{if }m\!<\!1\!+\!n\mbox{ and } a_0,\ldots,a_m\in X\mbox{ are pairwise 
distinct then}\\ 
&a_0\mbox{ is not in the subalgebra generated by }\{a_1,\ldots,a_m\}
\big\}. 
  \end{array}\]
(If $n=\omega$ we may omit it.)
\end{problem}
I think that for $n=\omega$ the case $|B|=\aleph_2$ had appeared in an old
list of Monk, but not in \cite{M2}. We may also ask $|B|\le({\rm irr}_n
(B))^{+m}$ for $n,m \in [2,\omega]$. 

Of course, the open case is when, say, $|B|=\lambda^{++}$, $2^\lambda\ge
\lambda^{++}$. Thinking about this problem, I was sure the answer is
consistently no (consistently yes is easy, even in the Easton
model). Moreover, I feel I know how to do it: let $\lambda=\lambda^{<
\lambda}$, $\lambda^+<\theta<\mu$, ${\mathfrak B}$ be a suitable algebra on
$\mu$ with $\le \lambda$ functions. A member $p$ of the forcing notion $\bP$
consists of $w^p\in [\mu]^{<\lambda}$ and a Boolean algebra $B^p$ generated
by $\{x_i:i\in w^p\}$ but such that
\begin{quotation}
if $B^p\models$`` $x_{i_0}=\sigma(x_{i_1},\ldots,x_{i_n})$ '' for ordinals
$i_0,\dotsc,i_n \in w^p$ and a Boolean term $\sigma$\\
then $i_0\in c\ell_{\mathfrak B}(\{i_1,\ldots,i_n\})$ and possibly more 
\end{quotation}
(the order of $\bP$ is the natural one). This is to reconcile the demand
``the Boolean algebra has cardinality $\mu$, so without loss of generality
we have to ask $i<j\ \Rightarrow\ x_i \ne x_j$'' and the $\lambda^+$-c.c. Of
course, $\name{B}=\bigcup\{B^p:p\in\name{G}_\bP\}$.

Note: if $\lambda=\aleph_0$ we have more freedom. (The expected proof goes:
if 
\[p\Vdash\mbox{`` }\name{X}\subseteq\name{B},\ \name{X}=\{\name{y}_i:i<
\theta\}\mbox{ with no repetition exemplifies }{\rm irr}(\name{B})^+\ge |B|
\mbox{ ''}\]
then we can find $p_i$, $p\le p_i\in\bP$, $p_i\Vdash$`` $y_i=\sigma_i(
x_{\alpha(i,0)},\ldots,x_{\alpha(i,n_i)})$ ''. Hence, if $(\forall\alpha<
\theta)(|\alpha|^{<\lambda}<\theta=\cf(\theta))$, without loss of generality,
$\sigma_i=\sigma$, $n_i=n^*$ and $\langle (p_i,\langle\alpha(i,\ell):\ell<
n^*\rangle):i<\theta\rangle$ forms a $\Delta$-system (moreover indiscernible
as a sequence of sequences of ordinals of length $<\lambda$). Let $\{\alpha(
i,\ell):\ell<n^*\}\subseteq w^p=w^*\cup\{\gamma_{i,\zeta}:\zeta<\zeta^*\}$,
$\zeta^*<\lambda$, etc, and we have to find $n<\omega_1$, $i_0<\ldots<i_n<
\theta$ and $q$ above $p_{i_0},\ldots,p_{i_n}$ such that $B^q\models$``
$y_{i_0}=\sigma^*(y_{i_1},\ldots,y_{i_n})$ ''. So it is natural to demand
$\zeta<\zeta^*\ \Rightarrow\ \gamma_{i_0,\zeta}\in c\ell_{\mathfrak B}
\{\gamma_{i_1,\zeta},\ldots,\gamma_{i_n,\zeta}\}$. So if $\lambda=\aleph_0$
we may use $n>|w^{p_i}|$. But this approach has not converged to a proof.)

So (see Monk \cite[Problem 28]{M2})

\begin{question}
\label{6.10}
\begin{enumerate}
\item Is there a class of (or just one) $\lambda$ such that for some Boolean
Algebra $B$ of cardinality $\lambda^+$ we have ${\rm irr}(B)=\lambda$?
\item Similarly for ${\rm irr}_n(B)$.
\end{enumerate}
\end{question}
Colouring theorems (e.g.~\cite{Sh:572}) are not enough for a construction.

\begin{question}
\label{6.11}
\begin{enumerate}
\item For which pairs $(\lambda,\theta)$ of cardinals $\lambda\ge\theta$ is
there a superatomic Boolean Algebra with $>\lambda$ elements, $\lambda$
atoms and every $f\in{\rm Aut}(B)$ moves $<\theta$ atoms?  

\noindent (That is $|\{x:B\models$`` $x$ an atom and $f(x)\ne x$ ''$\}|<
\theta$). 
\item In particular, is it true that for some $\theta$, for a proper class
of $\lambda$'s there is such Boolean Algebra? 
\item Replace ``automorphism'' by ``one--to-one endomorphism''.
\end{enumerate}
\end{question}
See some results in \cite[\S1,\S2]{Sh:641} for $\theta$ strong limit
singular. (It may be interesting to try: with $n$ depending on the arity of
the term as in \cite{Sh:641}.) 

Concerning attainment in ZFC: 

\begin{question}
\label{6.12}
\begin{enumerate}
\item Can we show the distinction made between the attainments of variants
of ${\rm hL}$ (and ${\rm hd}$), in a semi-ZFC way?  That is, in
Ros{\l}anowski and Shelah \cite{RoSh:651} such examples are forced. ``Semi
ZFC'' means can we prove such examples exist after adding to ZFC only
restrictions on cardinal arithmetic? 
\item Similarly for other consistency results.  (Well, preferably of low
consistency strength).
\end{enumerate}
\end{question}

In view of \cite{Sh:575}, \cite{Sh:620} it is reasonable to consider

\begin{problem}
\label{6.13}
In \ref{6.9}, \ref{6.10} replace ${\rm irr}_n(B)$ by ${\rm irr}_n(\bar{a},
B)$ (this is true for other cardinal invariants as well), see Definition
\ref{6.14} below.
\end{problem}

\begin{definition}
\label{6.14}
Let $B$ be a Boolean Algebra and $\bar{a}=\langle a_i:i<\lambda\rangle$ be a
sequence of elements of $B$. 
\begin{enumerate}
\item ${\rm irr}_n(\bar{a},B)=\sup\{|X|:X\subseteq\lambda\mbox{ and }\langle
a_i:i\in X\rangle\mbox{ is } n\mbox{--irredundant}\}$.
\item Similarly for other invariants of the ``universal family'' from
Ros{\l}anowski and Shelah \cite{RoSh:534} (see Definition 1.1 there).
\end{enumerate}
\end{definition}

\begin{question}
\label{6.15} 
Is there (at least consistently) a Boolean Algebra $B^*$, such that if $B$
is a Boolean Algebra extending $B^*$ {\em then\/} for some ultrafilters
$D_1,D_2$ on $B$ we have:\ \ $(B,D_1),(B,D_2)$ are not isomorphic, i.e., no
automorphism of $B$ maps $D_1$ onto $D_2$.\\
A close topological relative is: ``is there a homogeneous compact Hausdorff
space of cellularity $>2^{\aleph_0}$'' (van Douwen, see Kunen \cite{Ku90}).
\end{question}

\[* \qquad * \qquad *\]

There are some lemmas in \cite{Sh:233} which help to prove \ref{6.5A}, and I
would like to know whether the bounds used there are the best possible. 
Those lemmas also show that for some cardinal invariants (for Boolean
Algebras or topologies), defined by supremum, if the supremum is not
attained, then the value is ``almost" regular (the classical result of
Erd\"os and Tarski on the cellularity on Boolean Algebra (or topology) says
it is regular, whereas we get in \cite{Sh:233} that, e.g., the spread
satisfies $2^{\cf({\rm s}^+(B))}>{\rm s}^+(B)$, for a Boolean Algebra $B$,
$2^{2^{\cf({\rm s}^+(X))}}>{\rm s}^+(X)$ for a Hausdorff space, $2^{\cf({\rm
s}^+(X))}>{\rm s}^+(X)$ for a $T_3$ space $X$).

\begin{question}
\label{6.15a}
Can we find more applications of the theorems (and proofs) in \cite{Sh:233}
implying (or saying) that if the supremum in some cardinal invariants for a
space $\mathcal U$ (or a Boolean Algebra $B$, or whatever) is not attained,
then it has large cofinality?
\end{question}

A recent application is in Ros{\l}anowski and Shelah \cite[\S6]{RoSh:599}.
This is a converse to \ref{6.12}.

For the spread (and the hereditarily Lindelof degree and the hereditarily
density) the results are best possible (see Juh\'asz and Shelah
\cite{JuSh:231}), and for regular spaces we have better results
(\cite[5.1]{Sh:233} also best possible) but are the bounds in the claims
below best possible? 

\begin{definition}
\label{6.16} 
\begin{enumerate}
\item $\varphi$ is nice for $X$ if $\varphi$ is a function from the family
of subsets of the topological space $X$ to cardinals satisfying
\begin{quotation}
$\varphi(A)\le\varphi(A\cup B)\le\varphi(A)+\varphi(B)+ \aleph_0$\\ 
(i.e., monotonicity and subadditivity).
\end{quotation}
\item We say $\varphi$ is $(\chi,\mu)$--complete provided that if $A_i
\subseteq X,\varphi(A_i)<\chi$ for $i<\mu$ then $\varphi(\bigcup\limits_{i<
\mu} A_i)<\chi$.\\
Let $C(\varphi,\mu)=\{\chi:\varphi\mbox{ is }(\chi,\mu)\mbox{--complete}\}$.
\item We say $\varphi$ is $(<\lambda,\mu)$--complete, if for arbitrarily
large $\chi<\lambda$, $\varphi$ is $(\chi,\mu)$--complete.
\item Let ${\rm Ch}_\varphi$ be the following function from $X$ to
cardinals: 
\[{\rm Ch}_\varphi(y)=\min\{\varphi(u):y\in u\in\tau(X)\},\]
where $\tau(X)$ is the topology of $X$, that is the family of open sets.
\end{enumerate}
\end{definition}

\begin{remark}
\label{6.17}
\begin{enumerate}
\item We can replace $\mu$ by $<\mu$ and $i<\mu$ by $i<\alpha<\mu$, and make
suitable changes later.
\item In our applications we can restrict the domain of $\varphi$ to the
Boolean Algebra generated by $\tau(X)$ and even more, e.g., in \ref{6.18},
\ref{6.20} below to simple combinations of the $u_{i,\xi,\zeta}$.
\item We can change the definition of $(<\lambda,\mu)$--complete to
\begin{quotation}
if $A_i\subseteq X$ (for $i<\mu$) and $\sup\limits_{i<\mu} \varphi(A_i)<
\lambda$ then $\varphi(\bigcup\limits_{i<\mu} A_i)<\lambda$,
\end{quotation}
without changing our subsequent use. [We then will use: if $\varphi(A_\alpha
)<\chi_i$ for $\alpha<\mu$ then $\varphi(\bigcup\limits_{\alpha<\mu}
A_\alpha)<\chi_{i+1}$].
\end{enumerate}
\end{remark}

\begin{lemma}
\label{6.18}
Suppose $\lambda$ is a singular cardinal of cofinality $\theta$, $\lambda=
\sum\limits_{i<\theta}\chi_i$, $\chi_i<\lambda$, $\theta<\lambda$ and $\mu=
\beth_5(\theta)^+$, or even just $\mu=\beth_2(\beth_2(\theta)^+)^+$. Assume
that 
\begin{enumerate}
\item[(i)] $\varphi$ is nice for $X$,
\item[(ii)] $X_{\chi_i}=\{y\in X:{\rm Ch}_\varphi(y)\ge\chi_i\}$ has
cardinality $\ge\mu$ for $i<\theta$, 
\item[(iii)] $\varphi$ is $(<\lambda,\mu)$--complete.
\end{enumerate}
{\em Then\/} there are open sets $u_i\subseteq X$ (for $i<\theta$) such that
\[\varphi(u_i\setminus\bigcup_{j\ne i} u_j)\ge\chi_i.\]
\end{lemma}

\begin{remark}
\label{6.19} 
If $|\{y\in X:{\rm Ch}_\varphi(y)\ge\chi_i\}|<\mu$, it essentially follows
from $(\chi_i,\mu)$--completeness that $\varphi(X_{\chi_i})\ge\lambda$,
where $X_\chi=\bigcup\{u\in\tau(X): \varphi(u)<\chi\}$. Otherwise,
$\varphi(X\setminus X_{\chi_i})\ge\lambda$ by subadditivity, but $\varphi(X
\setminus X_{\chi_i})\le\prod\{\varphi(\{y\}):y\in X\setminus X_{\chi_i}\}$,
so by $(\chi_i,\mu)$--completeness for some $y\in X$, $\varphi(\{y\})\ge
\chi_i$, which is impossible for the instances which interest us. 
\end{remark}

\begin{lemma}
\label{6.20}
Suppose that $X$ is a Hausdorff space, $\lambda$ is a singular cardinal,
$\theta=\cf(\lambda)$, $\lambda=\sum\limits_{i<\theta}\chi_i$, $\chi_i<
\lambda$, $\mu<\lambda$ and clauses (i), (ii), (iii) of \ref{6.18} hold (for
$\varphi$). 
\begin{enumerate}
\item If $\mu=\beth_2(\theta)^+$ (or even $\sum\limits_{\sigma<\theta}
\beth_2 (\sigma)^+$), {\em then\/} there are open sets $u_i$ (for
$i<\theta$) such that $\varphi(u_i\setminus\bigcup\limits_{j>i}
u_j)\ge\chi_i$. 
\item If $X=\bigcup\{u:\varphi(u)<\lambda\}$, and $\mu$ is as in part (1),
{\em then\/} there are open sets $u_i$ (for $i<\theta$) such that $\varphi(
u_i\setminus\bigcup\limits_{j\ne i}u_j)\ge\chi_i$. 
\item If $\mu\ge\beth_3(2^{<\theta})^+$, $\varphi$ is $(<\chi_0,
\mu)$--complete, {\em then\/} there are $u_i$ (for $i<\theta$) such that
$\varphi(u_i\setminus\bigcup\limits_{j\ne i} u_j)\ge\chi_0$ (so $\lambda$,
$\chi_i$ ($0<i<\theta$) are irrelevant).
\end{enumerate}
\end{lemma}

\begin{remark}
\label{6.20A}
Part (1) of the lemma is suitable to deal with Boolean Algebras, part (2)
with the existence of $\{x_\alpha:\alpha<\lambda\}$ such that for every
$\alpha<\lambda$ for some $u$, $x_\alpha\in u\cap\{x_\beta:\beta<\lambda\}
\subseteq \{x_\beta:\beta\le\alpha\}$. 
\end{remark}

Now,

\begin{question}
\label{6.21}
Are the cardinal bounds in \ref{6.17} --- \ref{6.20} best possible?
\end{question}

\section{A taste of Algebra}
I have much interest in Abelian groups, but better see Eklof and Mekler
\cite{EM}. 

Thomas prefers to deal just with short elegant proofs of short elegant
problems (for me the second demand suffices). So he was rightly happy when
proving that for any infinite group $G$ with no center,
$\gamma(G)<(2^{|G|})^+$, where $G^{[0]}=G$, $G^{[1]}$ is the automorphism
group of $G$ considered as an extension of $G$, $G^{[i+1]}=(G^{[i]})^{[1]}$,
$G^{[\delta]}=\bigcup\limits_{i<\delta} G^{[i]}$, so $G^{[i]}$ is an
increasing sequence of groups with no center, and $\gamma(G)=\min\{\gamma:
G^{[\gamma]}=G^{[\gamma+1]}\}$. 

But is there a better cardinal bound?  No, for $|G|$ regular $>\aleph_0$,
see Just, Shelah and Thomas \cite{JsShTh:654}, but we are left with:

\begin{question}
\label{tal.1}
If $G$ is a countable group with a trivial center, then do we have
$\gamma(G)<\omega_1$? What about singular $|G|$?
\end{question}

I heard about the following problem (see Hamkins \cite{Hm98}).

\begin{problem}
\label{tal.2}
If $G$ is a group possibly with center, $G^{[i]}$ is defined as above but we
have just a homomorphism $h_{i+1,i}:G^{[i]}\longrightarrow G^{[i+1]}$ with
the center of $G^{[i]}$ being the kernel (and in limit stages take the
direct limit), is there a bound to $\gamma(G)$ really better than the first
strongly inaccessible $>|G|$ (gotten by Hamkins \cite{Hm98})?
\end{problem}

Thomas also had started investigating cofinalities of some natural groups,
(see Sharp and Thomas \cite{ST94}, \cite{ST95}, Thomas \cite{Th96}). He drew
me to it and I was particularly glad to see that pcf pops in naturally;
e.g., (see Shelah and Thomas \cite{ShTh:524}) if $\lambda_n\in {\rm
CF}({\rm Sym}(\omega))$ and $\lambda\in \pcf\{\lambda_n:n<\omega\}$ then
$\lambda\in {\rm CF}({\rm Sym}(\omega))$, where

\begin{definition}
\label{tal.3}
\begin{enumerate}
\item ${\rm CF}(G)=\{\theta:\theta=\cf(\theta)$ and there is an increasing
sequence of proper subgroups of $G$ of length $\theta$ with union $G\}$.
\item $\cf(G)=\min[{\rm CF}(G)\setminus\{\aleph_0\}]$.
\end{enumerate}
\end{definition}
Though we found some information about $\cf(\prod\limits_n{\rm alt}(n))$
(see Saxl, Shelah and Thomas \cite{ShST:584}, where ${\rm alt}(n)$ is the
group of even permutations of $\{0,\ldots,n-1\}$), we remained baffled by

\begin{question}
\label{tal.4}
Is it consistent that $\aleph_2\le\cf(\prod\limits_{n<\omega}{\rm alt}(n))$?
\end{question}
It is natural to try to use iterations of length $\omega_2$, where each
iterand consists of trees with norms (see Ros{\l}anowski and Shelah
\cite{RoSh:470}). Naturally, a norm on ${\mathcal P}({\rm alt}(n))$ will be
such that if ${\bf nor}(A)\ge m+1$ and $\sigma$ is a group term, then we can
have for ``many'' $f_1,\ldots,f_k,g \in{\rm alt}(n)$ that $A'=\{h\in A:
\sigma(h,f_1,\ldots,f_{k^*})=g\}$ has ${\bf nor}(A')\ge m$ toward destroying
a guess on an approximation to a lower subgroup exemplifying $\omega_1\in
{\rm CF}(\prod\limits_{n<\omega}{\rm alt}(n))$. This helps for $\aleph_0
\notin {\rm CF}(\prod\limits_{n<\omega}{\rm alt}(n))$, but fails for the
purpose of \ref{tal.4}.

\[* \qquad * \qquad *\]

My interest in lifting for the measure algebra started when Talagrand
promised me ``flowers on your grave from every measure theorist'' (a little
hard to check), and Fremlin said essentially the same, conventionally
(see \cite{Sh:185}, \cite[Ch.IV]{Sh:f}). But this does not solve some other
problems from Fremlin's list, from which I particularly like

\begin{question}
\label{tal.5}
Assume CH (or even GCH or just prove consistency).

Do we have lifting for every measure algebra? 

\noindent Which means: let ${\mathcal B}(I)$ be the algebra of subsets of
${}^I 2$ generated from clopen ones by countable unions and intersections,
$\mu_{\mathcal B}$ the Lebesgue measure on ${\mathcal B}$ (so we get the
so-called Maharam algebra), ${\bf I}=\{A\in {\mathcal B}:\mu_{\mathcal B}(A)
=0\}$ (so $\bf I$ is the ideal of null sets). A lifting is a homomorphism
from ${\mathcal B}/{\bf I}$ into ${\mathcal B}$ such that 
\[X\in {\mathcal B}\quad \Rightarrow\quad f(X/{\bf I})=X\mod{\bf I}.\] 
\end{question}
Naturally, we think $I_0\subseteq I_1\ \Rightarrow\ {\mathcal B}(I_0)
\subseteq {\mathcal B}(I_1)$ (by identifying) and for an increasing sequence
$\bar{I}=\langle I_\alpha:\alpha<\alpha^*\rangle$ we let ${\mathcal B}(
\bar{I})=\bigcup\limits_{\alpha<\alpha^*}{\mathcal B}(I_\alpha)$. In the
positive direction we may try to prove by induction on $\lambda$; then we
will be naturally drawn to proving: for any ${\mathcal P}^-(n)$-diagram
$\langle {\mathcal B}(\bar{I}_s):s\in {\mathcal P}^-(n)\rangle$, where
${\mathcal P}^-(A)=\{u:u\subseteq A,\ u\ne A\}$, and a sequence of liftings
$\bar{f}=\langle f_s:s\in {\mathcal P}^-(n)\rangle$ satisfying a reasonable
induction hypothesis, $\bigcup\limits_s f_s$ can be extended to a lifting of
$\bigcup\limits_s I_{s,\alpha}$ (as e.g.~in \cite{Sh:87b}, Sageev and Shelah
\cite{SgSh:217}). For the negative direction we may think of using a
partition theorem. 

\[* \qquad * \qquad *\]

For a long time I have been interested in compactness in singular cardinals;
i.e., whether if something occurs for ``many'' subsets of a singular
$\lambda$ of cardinality $<\lambda$, it occurs for $\lambda$. For the
positive side (on the filters see \ref{tal.7} below)

\begin{theorem}
\label{tal.6}
Let $\lambda$ be a singular cardinal, $\chi^*<\lambda$. Assume that ${\bf
F}$ is a set of pairs $(A,B)$ (written usually as $B/A$; ${\bf F}$ stands for
free) $A,B\subseteq\cU$ satisfying the axioms II, III, IV,VI,VII below. Let
$A^*,B^* \subseteq\cU$, $|A^*|=\lambda$, {\em then\/} $B^*/A^*\in{\bf F}$ if
$B^*/A^*$ is $\lambda$-free in a weak sense which means (see Definition
\ref{tal.7} below): 
\begin{enumerate}
\item[$(*)_0$] for the ${\mathcal D}_{\chi^*}(B^*)$--majority of $B\in
[B^*]^{< \lambda}$ we have $B/A^*\in{\bf F}$,\\ 
or just
\item[$(*)_1$] the set $\{\mu<\lambda:\{B\in [B^*]^\mu:B/A^*\in{\bf F}\}\in
{\mathcal E}^{\mu^+}_\mu(B^*)\}$ contains a club of $\lambda$,\\
or at least
\item[$(*)_2$] for some set $C$ of cardinals $<\lambda$, unbounded in
$\lambda$ and closed (meaningful only if $\cf(\lambda)>\aleph_0$), for every
$\mu\in C$, for an ${\mathcal E}^{\mu^+}_\mu(B^*)$--positive set of $B\in
[B^*]^\mu$ we have $B/A^*\in{\bf F}$. 
\end{enumerate}
\end{theorem}
The axioms are

\begin{enumerate}
\item[Ax II:]  $B/A\in{\bf F}\quad \Leftrightarrow\quad A\cup B/A\in {\bf
F}$, 
\item[Ax III:]  if $A\subseteq B\subseteq C$, $B/A\in{\bf F}$ and $C/B\in
{\bf F}$ {\em then\/} $C/A\in {\bf F}$,
\item[Ax IV:] if $\langle A_i:i\leq\theta\rangle$ is increasing continuous,
$\theta=\cf(\theta)$, $A_{i+1}/A_i\in {\bf F}$\\
then $A_\theta/A_0\in {\bf F}$,
\item[Ax VI:] if $A/B\in {\bf F}$ then for the ${\mathcal D}_{\chi^*
}$--majority of $A'\subseteq A$, $A'/B\in {\bf F}$ (see below), 
\item[Ax VII:] if $A/B\in {\bf F}$ then for the ${\mathcal D}_{\chi^*
}$--majority of $A'\subseteq A$, $A/B\cup A'\in {\bf F}$.
\end{enumerate}
(Of course we can get variants by putting more or less into the statement.)

\begin{definition}
\label{tal.7}
\begin{enumerate}
\item Let ${\mathcal D}$ be a function giving for any set $B^*$ a filter
${\mathcal D}(B^*)$ on ${\mathcal P}(B^*)$ (or on $[B^*]^\mu$). Then to say
``for the  ${\mathcal D}$--majority of $B\subseteq B^*$ (or $B\in
[B^*]^\mu$) we have $\varphi(B)$'' means $\{B\subseteq B^*:\varphi(B)\}\in
{\mathcal D}(B^*)$ (or $\{B\in [B^*]^\mu:\neg\varphi(B)\}=\emptyset\mod
{\mathcal D}$). 
\item Let ${\mathcal D}_\mu(B^*)$ be the family of $Y\subseteq {\mathcal P}(
B^*)$ such that for some algebra $M$ with universe $B^*$ and $\le\mu$
functions, 
\[Y\supseteq\{B\subseteq B^*:B\ne\emptyset \mbox{ is closed under the
functions of } M\}.\]  
\item ${\mathcal E}^{\mu^+}_\mu(B^*)$ is the collection of all $Y\subseteq
[B^*]^\mu$ such that:\quad for some $\chi,x$ such that $\{B^*,x\}\in
{\mathcal H}(\chi)$,\\
{\em if\/} $\bar{M}=\langle M_i:i<\mu^+\rangle$ is an increasing continuous
sequence of elementary submodels of $({\mathcal H}(\chi),\in)$ such that
$x\in M_0$ and $\bar{M}\restriction (i+1)\in M_{i+1}$,\\
{\em then} for some club $C$ of $\mu^+$, $i\in C\ \Rightarrow\ M_i\cap B^*
\in Y$. 
\end{enumerate}
\end{definition}

On the filters see Kueker \cite{Ku77}, and \cite[\S3]{Sh:52}. The theorem
was proved in \cite{Sh:52} but with two extra axioms, however it included
the full case of varieties (i.e., including the non-Schreier ones). Later,
the author eliminated those two extra axioms: Ax V in Ben David \cite{BD},
and Ax I in \cite{Sh:E18} (answering a question of Fleissner on providing a 
combinatorial proof of the compactness). Hodges \cite{Ho81} contains also
presentation of variants of this result. 

There are some cases of incompactness (see Fleissner and Shelah
\cite{FlSh:267}, and \cite{Sh:347}).

\begin{problem}
\label{tal.8} 
\begin{enumerate}
\item Are there general theorems covering the incompactness phenomena?
\item Are there significantly better compactness theorems (for uncountable
cofinality, of course)? 
\end{enumerate}
\end{problem}

Related is 

\begin{question}
\label{tal.9}
What can be 
\[\{\lambda:\mbox{there is a }\lambda\mbox{--free for }{\mathcal V}\mbox{
algebra $M$ of cardinality $\lambda$ which is not free}\},\]
for a variety ${\mathcal V}$ (at least with countable vocabulary) ?\\ 
(See Eklof and Mekler \cite{EM}, Mekler and Shelah \cite{MkSh:366}, Mekler,
Shelah and Spinas \cite{MShS:417}.) 
\end{question}

\[* \qquad * \qquad *\]

There are cases of strong dichotomy: if $\ge\lambda$ then $\ge 2^\lambda$,
related to groups (see \cite{Sh:273}, Grossberg and Shelah \cite{GrSh:302},
\cite{GrSh:302a}, and \cite{Sh:664}; on Abelian groups see Fuchs \cite{Fu}).

\begin{question}
\label{tal.10}
[$\bV=\bL$]  
If $\lambda>\cf(\lambda)>\aleph_0$, $G$ is a torsion free Abelian group of
cardinality $\lambda$, can $\lambda=\nu_p({\rm Ext}(G,{\mathbb Z}))$?
\end{question}
The cardinal $\nu_p({\rm Ext}(G,{\mathbb Z}))$ is the dimension of $\{x\in
{\rm Ext}(G,{\mathbb Z}):px =0\}$ as a vector space over ${\mathbb
Z}/p{\mathbb Z}$. To avoid ${\rm Ext}$ note that (see Fuchs \cite{Fu}) this
group can be represented as ${\rm Hom}(G,{\mathbb Z}/p{\mathbb Z})/\{h/p:
h\in {\rm Hom}(G,{\mathbb Z})\}$. If $G$ is torsion free, then the group
${\rm Ext}(G,{\mathbb Z})$ is divisible and hence the ranks $\nu_p({\rm
Ext}(G,{\mathbb Z}))$ for $p$ prime, and $\nu_0({\rm Ext}(G,{\mathbb Z}))$,
the rank of $\{x\in {\rm Ext}(G,{\mathbb Z}):x \mbox{ torsion}\}$, determine
${\rm Ext}(G,{\mathbb Z})$ up to an isomorphism. If we assume ($\bV=\bL$
and) there is no weakly compact cardinal, this question is the only piece
left for characterizing the possible such ${\rm Ext}(G,{\mathbb Z})$, see
Mekler, Ros{\l}anowski and Shelah \cite{MRSh:314}.  

\begin{question}
\label{tal.11} 
What is the first cardinal $\lambda=\lambda_\kappa$ such that: for every
ring $R$ of cardinality $\le \kappa$, if there is endorigid (or rigid, or
1-to-1 rigid) $R$--module of size $\ge\lambda$, {\em then\/} there are such
$R$-modules in arbitrarily large cardinals?  (I.e. Hanf numbers).
\end{question}

\section{Partitions and colourings}
Remember (see Erd\H{o}s, Hajnal, Mat\'e and Rado \cite{EHMR})

\begin{definition}
\label{col.1} 
\begin{enumerate}
\item $\lambda\rightarrow(\alpha)^n_\kappa$ means:\quad for every colouring
$c:[\lambda]^n\longrightarrow\kappa$ there is a set $X\subseteq\lambda$ of
order type $\alpha$ such that $c\restriction [X]^n$ is constant.
\item $\lambda\rightarrow [\alpha]^n_\kappa$ means:\quad for every colouring
$c:[\lambda]^n\longrightarrow\kappa$ there is a set $X\subseteq\lambda$ of
order type $\alpha$ such that $\Rang(c \restriction [X]^n)\ne\kappa$.
\item $\lambda\rightarrow [\alpha]^n_{\kappa,\sigma}$ means:\quad for every
colouring $c:[\lambda]^n\longrightarrow\kappa$ there is a set $X\subseteq
\lambda$ of order type $\alpha$ such that $\Rang(c\restriction [X]^n)$ has
cardinality $<\sigma$.
\end{enumerate}
\end{definition}

\begin{definition}
\label{col.2}
$M$ is a Jonsson algebra if it is an algebra with countably many functions
with no proper subalgebra of the same cardinality.\\ 
(See \cite{Sh:g}, \cite{Sh:E12}.)
\end{definition}

\begin{definition}
\label{col.3} 
${\rm Pr}_1(\lambda,\mu,\theta,\sigma)$ means:\quad there is $c:[\lambda]^2
\longrightarrow\theta$ such that if $u_i\in [\lambda]^{<\sigma}$ (for
$i<\mu$) are pairwise disjoint and $\gamma<\theta$ {\em then\/} for some
$i<j<\mu$ we have $c\restriction (u_i\times u_j)$ is constantly $\gamma$.\\
(See \cite{Sh:g}, \cite{Sh:E12}.)
\end{definition}
There are many more variants.

It irritates me that after many approximations, I still do not know (better
for me ``consistently no'', better for set theory ``yes'') the answer to the
following.

\begin{question}
\label{col.4} 
If $\mu$ is singular, is there a Jonsson algebra on $\mu^+$?  (and even
better ${\rm Pr}_1(\mu^+,\mu^+,\mu^+,\cf(\mu))$?)
\end{question}
Also the requirements on an inaccessible to get colouring theorems may well
be an artifact of our inability, so let us state minimal open cases.

\begin{question}
\label{col.6}
\begin{enumerate}
\item Let $\lambda$ be the first $\omega$--Mahlo cardinal. Does
$\lambda\nrightarrow [\lambda]^2_\lambda$ or at least $\lambda\nrightarrow
[\lambda]^2_\theta$ for $\theta<\lambda$? 
\item Let $\lambda$ be the first inaccessible cardinal which is
$(\lambda\cdot \omega)$--Mahlo. Is there a Jonsson algebra on $\lambda$?
(Even better $\lambda\nrightarrow [\lambda]^2_\lambda$?)
\end{enumerate}
In both parts it is better to have ${\rm Pr}_1(\lambda,\lambda,\lambda,
\aleph_0)$, etc; it is interesting even assuming GCH. 
\end{question}
Whereas under GCH the relation $\rightarrow$ for cardinals is essentially
understood (see Erd\H{o}s, Hajnal, Mat\'e and Rado \cite{EHMR}), the case of
ordinals is not.  As by \cite{Sh:26} (GCH for simplicity) $\alpha<\lambda=
\cf(\lambda)\ \Rightarrow\ \lambda^{++}\rightarrow(\lambda^+ +\alpha)^2_2$
(on the problem see \cite{EHMR}), and by Baumgartner, Hajnal and
Todor\v{c}evi\'{c} \cite{BHT93} for $k<\omega$,
\[\alpha<\lambda=\cf(\lambda)\quad \Rightarrow\quad \lambda^{++}\rightarrow
(\lambda^++\alpha)^2_k,\]
it remains open: 

\begin{question}
[GCH] 
\label{col.7}
When does $\lambda^{++}\rightarrow(\lambda^+ +\alpha)^2_{\aleph_0}$, where
$2 \le \alpha<\lambda$?\\
(The GCH assumption is for simplicity only.)
\end{question}
By \cite{Sh:424} we know that if $\lambda=\lambda^{<\lambda}$ (e.g.,
$\lambda=\aleph_1=2^{\aleph_0}$) then possibly $2^\lambda$ is very large (in
particular $>\lambda^{+\omega}$), but $2^\lambda \nrightarrow (\lambda\times
\omega)^2_2$. However, $\lambda^{+2k}\rightarrow (\lambda\times n)^2_k$ (by 
\cite[bottom of p.288]{Sh:424}), so 

\begin{question}
\label{col.8}
For $\lambda=\lambda^{<\lambda}$, $\lambda>\aleph_0$, and $k<\omega$ and
$n<\omega$, what is the minimal $m$ such that $\lambda^{+m}\rightarrow
(\lambda \times n)^2_k$ ?\\  
(Baumgartner, Hajnal and Todor\v{c}evi\'{c} \cite[p.2, end of \S0]{BHT93}
prefer to ask whether $\lambda^{++}\rightarrow (\lambda+\omega)^2_3$ for
$\lambda=\aleph_1$, so $\lambda=\lambda^{<\lambda}$ means CH as they choose
another extreme case of the unknown).
\end{question}
Now, for me a try at consistency of negative answers for \ref{col.7} calls
for using historic forcing (see Shelah and Stanley \cite{ShSt:258},
Ros{\l}anowski and Shelah \cite[\S 3]{RoSh:651}; it is explained below). On
the other hand, large cardinals may make some positive results easier. 

\begin{question}
\label{col.8a}
Assume that $\lambda>\kappa>|\zeta|+\sigma$, and $\kappa$ is a compact
cardinal, and $\lambda=\lambda^{<\lambda}$. Does it follow that $\lambda^+
\rightarrow (\lambda + \zeta)^2_\sigma$?
\end{question}
\noindent I tend to think the answer is yes.  If so, then we cannot expect,
when $\lambda =\lambda^{<\lambda}$, that a $\lambda$--complete
$\lambda^+$-c.c.~forcing notion will not add a counterexample $c$ to
$\lambda^+\rightarrow (\lambda+\alpha^*)^2_{\aleph_0}$. Let for simplicity
$\alpha^*=2$. Also we should expect that for such $c$, for every $u\in
[\lambda^+]^{<\lambda}$ we can find $f_n:u\times u\longrightarrow\omega$ for
$n<\omega$ such that $f_n(\alpha,\alpha)=n$, and for no distinct
$\alpha,\beta,\gamma\in u$ do we have:
\[c(\{\alpha,\beta\})=n=f_n(\alpha,\beta)\quad\mbox{ or }\quad f_n(\alpha,
\beta) =f_n(\alpha,\gamma)=c(\{\beta,\gamma\}).\]
The point of historic forcing is that we know what kind of object our
forcing notion has to add. In our case, we assume $\lambda=\lambda^{<
\lambda}$ and a condition $p$ has to give $u^p\in [\lambda^+]^{<\lambda}$
and $c^p:[u^p]^2\longrightarrow\omega$ and by the above considerations also
$f_n:u^p\times u^p\longrightarrow\omega$ (for $n<\omega$), and for having
some leeway we let $f_n:u^p\times u^p\longrightarrow [\omega]^{\aleph_0}$
such that $n\in f_n(\alpha,\alpha)$, and we demand
\begin{enumerate}
\item[$(*)$] for no distinct $\alpha,\beta\in u^p$ and $m$ do we have
\[m=c(\{\alpha,\beta\})\in f_n(\alpha,\beta)\cap f_n(\alpha,\alpha),\]
\item[$(**)$] for no distinct $\alpha,\beta,\gamma$ from $u^p$ do we
have
\[c(\{\beta,\gamma\})\in f_n(\alpha,\beta)\cap f_n(\alpha,\gamma).\]
\end{enumerate}
Moreover, we choose $\langle A_\eta:\eta\in {}^{\omega>}\omega\rangle$,
$A_\eta \in [\omega]^{\aleph_0}$, $\langle A_{\eta\conc\langle\ell\rangle}:
\ell<\omega\rangle$ are pairwise disjoint subsets of $A_\eta$ and demand
$f_n(\alpha,\beta)\in \{A_\eta:\eta\in {}^{\omega>}\omega\}$.

Considering that our forcing will be strategically $(<\lambda)$--complete
(as $\lambda$--complete seems too much both if the answer to \ref{col.8a} is
yes and because of the properties of historic forcing in general), in order
that there will be no $A\subseteq\lambda^+$, $\otp(A)=\lambda+\alpha^*$ such
that $\name{c}\restriction [A]^2$ is constant, we have a set ${\mathcal T}$
of cardinality $<\lambda$ such that (for simplicity $\lambda=\mu^+$, maybe
also $\mu$ regular): 
\begin{enumerate}
\item[(a)] each $x\in {\mathcal T}^p$ has the form $x=(a,\delta,b)=(a^x,
\delta^x,b^x)$, where 
\begin{enumerate}
\item[(i)] $\delta<\lambda^+$, $\cf(\delta)=\lambda$,
\item[(ii)] $a\subseteq u^p\setminus\delta$, $\otp(a)=\zeta$,
\item[(iii)] $b\subseteq u^p\cap\delta$, $\otp(b)=\mu$,
\end{enumerate}
\item[(b)] if $(a',\delta',b'),(a'',\delta'',b'')\in {\mathcal T}^p$ then
$|b' \cap b''|<\mu$ (probably not necessary), 
\item[(c)] $c^p\restriction [a\cup b]^2$ is constantly $n(x)$,
\item[(d)]  for no $\alpha\in u^p\cap\delta\setminus\sup(b)$ do we have
\[(\forall\beta\in a\cup b)[c^p(\{\alpha,\beta\})=n(x)].\]
\end{enumerate}
A $p$ like above will be a precondition. Now the ``history'' enters and the
proof should be clear just as the definition of the set of the conditions
should roll itself.

So atomic conditions will have $u^p$ a singleton, each condition will have a
history telling how it was created: each step in the history corresponds to
one of the reasons for creating a condition in the proof. Naturally, a major
reason is the proof of the $\lambda^+$-c.c.~by the $\Delta$--system
lemma. So assume for $\ell= 1,2$ that $p_\ell\in\bP$, $\alpha_1<\alpha_2$,
$u^{p_1}\cap\alpha_1=u^{p_2}\cap\alpha_2$, $u^{p_1}\subseteq\alpha_2$,
$\otp(u^{p_1})=\otp(u^{p_2})$ and the order preserving mapping ${\rm
OP}_{u^{p_2},u^{p_1}}$ from $u^{p_1}$ onto $u^{p_2}$ maps $p_1$ to $p_2$. We
have to amalgamate $p_1$ and $p_2$ getting $q$, so we have to determine
$f^q_n$ on $(u^{p_1}\setminus u^{p_2})\times (u^{p_2}\setminus u^{p_1})$,
and $c^q(\{\alpha,\beta\})$ for $\alpha\in u^{p_1}\setminus u^{p_2}$,
$\beta\in u^{p_2}\setminus u^{p_1}$.

But our vision is that there is one line of history, so should the history
of $q$ continue the history of $p_1$ with $p_2$ joining or should the
history of $q$ continue the history of $p_2$ with $p_1$ joining? Both are
O.K., but we get two distinct conditions $q',q''$ which, however, are
equivalent; i.e.~$q'\le q''\le q'$.  

Generally going back and changing the history as above we get an equivalent
condition if this is done finitely many times (this also explains why we get
strategical $(<\lambda)$--completeness and not $\lambda$--completeness). 
That is, we define $p\le_{\rm pr} q$ iff $p$ appears in the history of $q$,
$p \sim q$ if $p$ is gotten from $q$ by finitely many changes as above in
the history, and lastly $p\le q$ iff $(\exists p')(p'\sim p\ \&\ p'\le_{\rm
pr} q)$. 

What is left?  ``Only'' carrying out the amalgamation (using and
guaranteeing the conditions, and probably changing them so that they will
fit). 

\begin{question}
\label{col.8A}
What are the best cardinals needed for the canonization theorems in
\cite{Sh:95}?
\end{question}

An old well-known problem is 

\begin{question}
\label{Col.8B}
Is $\aleph_1\rightarrow [\aleph_1;\aleph_1]^2_2$ consistent?  (And variants,
connected to the $L$-space problem; this seems related to \ref{5.1}).
\end{question}
Another problem of Erd\"os is (the answer is consistently yes, see
\cite{Sh:289}, even colouring also no edges, but provability in ZFC is not
clear): 

\begin{question}
\label{col.9}
Is there a graph $G$ with no $K_4$ (complete graph on 4 vertices) such that
$G \rightarrow (K_3)^2_{\aleph_0}$, that is for any colouring of the edges
by $\aleph_0$ colours there is a monochromatic triangle?
\end{question}
We can ask it for $K_k,K_{k+1}$ instead of $k=3$ and colouring $r$--tuples
instead of pairs; the answer still is consistently yes (see \cite{Sh:289}),
so the problem is in ZFC. See \cite[Ch.III,\S1]{Sh:e} on a connection to
model theory. 

\[* \qquad * \qquad *\]

A very nice theorem of Hajnal \cite{Ha93} says that, e.g., for any finite
graph $G$ and $\kappa$ for some graph $H$, $H\rightarrow (G)^2_\kappa$, but
leaves as a mystery: 

\begin{question}
\label{col.10}
Let $G$ be a countable graph, is there a graph $H$ such that 
\[H\rightarrow(G)^2_{\aleph_0}\ ?\] 
\end{question}
Starting from a problem of Erd\"os and Hajnal \cite{EH}, I have been very
interested in consistency results, e.g., of the form $\lambda\rightarrow
[\mu]^2_3$, $\kappa<\mu<\lambda\le 2^\kappa$ (see \cite{Sh:546}). Usually
those are really canonization theorems, for a fixed natural coloring, any
other is, on a large set, computable from it. Those results help sometimes
in consistency results (just as, e.g., the Erd\"os-Rado theorem helps in ZFC
results). Still it seems to me worthwhile to know.

\begin{question}
\label{col.11}
\begin{enumerate}
\item Can we put together the results of, e.g., \cite{Sh:546}, Shelah and
Stanley \cite{ShSt:608} and \cite{Sh:288}?  Assume that $\kappa=\kappa^{<
\kappa}<\lambda$, $\lambda$ is, e.g., strongly inaccessible large
enough. Can we find a $(<\kappa)$--complete, $\kappa^+$--c.c.~forcing notion
$\bP$ such that in $\bV^\bP$: 
\begin{enumerate}
\item[(a)] for $\sigma<\kappa$ and $\mu<\lambda$ we can find $\mu',\lambda'$
such that $\mu<\mu'<\lambda'<\lambda$ and $\lambda'\rightarrow
[\mu']^2_{\sigma,2}$;
\item[(b)] if $\kappa$ is a measurable indestructibly by adding many Cohens,
then also the parallel results for colouring $n$-tuples (see \cite{Sh:288});
\item[(c)] if $\kappa=\aleph_0$, we also have results on colouring
$n$--tuples simultaneously for all $n$ ?
\end{enumerate}
\item Add the hopeful consistency answer for \ref{col.12}, \ref{col.13}.
\end{enumerate}
\end{question}
Probably easier are, e.g.

\begin{question}
\label{col.12}
Is it consistent that for some $n$, $2^{\aleph_0}=\aleph_n=\lambda
\rightarrow [\aleph_1]^2_3$ ?\\
(The exact $n$ is less exciting for me, the main division line seems to me
$\aleph_\omega$, of course best to know the exact $\lambda$.) 
\end{question}

\begin{question}
\label{col.13}
Is it consistent that $2^{\aleph_0}>\lambda\rightarrow [\aleph_1]^2_3$ ?
\end{question}

\[* \qquad * \qquad *\]

{\bf On Finite Combinatorics.}\qquad Spencer, Szemeredi and Alon told me
that finding $\lim(\log(r^2_2(n)/n))$ is a major problem (see Definition
\ref{col.21} below), but the difference between lower and upper bounds seems
to me negligible.   

\begin{definition}
\label{col.21}
$r^m_k(n)$ is the minimal $r$ such that $r \rightarrow (n)^m_k$.
\end{definition}
Erd\"os and Hajnal ask, and I find more convincing, the following.

\begin{question}
\label{col.22}
What is the order of magnitude of $r^3_2(n)$?
\end{question}
We expect it should be $2^{2^n}$, or e.g.~$2^{2^{(n^\varepsilon)}}$ for some
$\varepsilon>0$. But we cannot rule out its being $2^n$ or
e.g.~$2^{n^{1/\varepsilon}}$ for some $\varepsilon>0$.\\
Here the difference is large.

Note that for four colours the problem (what is $r^3_4(n)$) is settled; but
I think the true question is:

\begin{question}
\label{col.23}
Determine (order of magnitude is OK) $f_{k}(n,r),f^+_k(n,c)$ where:
\begin{enumerate}
\item[(a)] $f_3(n,c)$ is the minimal $m$ such that\\
for every ${\bf d}:[m]^3\longrightarrow\{0,\ldots,c-1\}$, there are $A\in
[m]^n$ and a strictly increasing function $h:A\longrightarrow\{0,\ldots,
2^n-1\}$ such that for $\ell_0<\ell_1<\ell_2$ in $A$, the value of ${\bf
d}(\{\ell_0,\ell_1,\ell_2\})$ is determined by the quantifier--free type
$\langle h(\ell_0),h(\ell_1),h(\ell_2)\rangle$ in $B^3_n$,\\
where $B^3_n$ has the universe $\{0,\ldots,2^n-1\}$ and two relations:
(viewed as ${}^n 2$) the lexicographic order and 
\[\{(\eta_0,\eta_1,\eta_2):\eta_0<_{\ell x}\eta_1<_{\ell x}\eta_2\mbox{ and
} \ell g(\eta_0 \cap \eta_1)<\ell g(\eta_1 \cap \eta_2)\}.\] 
\item[(b)] $f^+_3(n,c)$ is defined similarly but for every pregiven
$\Rang(h)$ we can find such $h$. 
\item[(c)] $B^k_n$ is defined below by induction of $k$, and then
$f_k(n,c),f^+_k(n,c)$ are defined analogously to $f_3,f^+_3$.
\item[(d)] Define canonization numbers $g_k(n), g^+_k(n)$, which is the
first $m$ such that: if $d:[n]^k\rightarrow C$, with no restriction on the
cardinality of $C$, then we can find $A,h$ as above, and a quantifier free
formula $\phi$ in the vocabulary of $B^k_n$ such that for any $u_1,u_2 \in
[A]^k$ we have
\[d(u_1)=d(u_2)\quad\mbox{iff}\quad \phi (\ldots,h(\ell_1),\ldots;\ldots,
h(\ell_2),\ldots)_{\ell_1\in u_1,\ell_2 \in u_2}\mbox{ is satisfied in
}B^k_n.\] 
\end{enumerate}
\end{question}
The explicit way to describe $B^k_n$, by induction on $k$ is: it has a
linear order $<_k$ :\\
$B^2_n$ is the structure $(n,<)$,\\
$B^{k+1}_n$ has universe ${}^{B^k_n}\{0,1\}$ and the relation:

$\eta <_{k+1} \nu$ if and only if 

for some $y=y(\eta,\nu)\in B^k_n$ we have $\eta(y)=0$, $\nu(y)=1$ and
\[\eta\restriction\{x\in B^k_n:x <_k y\}=\nu\restriction\{x\in B^k_n:x<_k
y\}.\]

\noindent For an $m$-place relation $R^{B^k_n}$ of $B^k_n$, $R^{B^{k+1}_n}$
is a $2m$-place relation on $B^{k+1}_n$, namely 
\[\{\langle \eta_0,\ldots,\eta_{2m-1}\rangle:\eta_\ell\in B^{k+1}_n,\
\eta_\ell<_{k+1}\eta_{\ell+1}\mbox{ and }\langle y(\eta_0,\eta_1),
y(\eta_2,\eta_3),\ldots\rangle\in R^{B^k_n}\}.\]
Now it is not clear how fast the number in \ref{col.23} grows, e.g., we
cannot exclude $2^{2^{n+c+k}}$. The main question is whether it grows like 
$h(k)$--iterated exponentiation in $n$ (say $c$ fixed) with $h$ going to
infinity, or with $h$ constant. Of course, enriching somewhat the structure
is not a great loss to me.

\begin{question}
\label{col.24}
Let $f^*(n,c)$ be the first $m$ such that if $\langle A_\ell:\ell<m\rangle$
are pairwise disjoint, $|A_\ell|=m$ for $\ell<n$ and 
\[F:\{w\subseteq\bigcup_\ell A_\ell:|w\cap A_\ell|\in \{1,2\},\ (\exists !
\ell)(|w|=1)\}\longrightarrow C,\]
where $|C|=c$ {\em then\/} for some $x_\ell\ne y_\ell$ from $A_\ell$ for
$\ell<n$ we have 
\[\ell^*<n\ \Rightarrow\ F(\{x_\ell,y_\ell:\ell\ne\ell^*\}\cup\{x_{\ell^*}
\})= F(\{x_\ell,y_\ell:\ell\ne \ell^*\}\cup\{y_{\ell^*}\}).\]
Again the main question for me is: Does $f^*(n,c)$ grow as a fixed iterated
exponentiation? \\
(This is connected to the van der Waerden theorem, see \cite{Sh:329}).
\end{question}

On the Ramsey Theory see Graham, Rothschild and Spencer \cite{GrRoSp}.

\begin{definition}
\label{col.25} 
\begin{enumerate}
\item For a group $G$ and a subset $A$ of $G$, and a group $H$ let
$H\rightarrow (G)^A_\sigma$ mean:
\begin{quotation}
{\em if\/} ${\bf d}$ is a function with domain $H$ and range of cardinality
$\le\sigma$,\\ 
{\em then\/} for some embedding $h$ of $G$ into $H$ the function ${\bf d}$
restricted to $h(A)$ is constant.
\end{quotation}
\item  For a group $G$ and a subset $A$ of $G$, and group $H$ let
$H\rightarrow [G]^A_{\sigma, \tau}$ mean:
\begin{quotation}
{\em if\/} ${\bf d}$ is a function with domain $H$ and range of cardinality
$\le\sigma$,\\
{\em then\/} for some embedding $h$ of $G$ into $H$ the range of the function 
${\bf d}$ restricted to $h(A)$ has cardinality $<\tau$.
\end{quotation}
If we omit $\tau$, we mean just that the range is not equal to $\sigma$.
\item For a group $G$ and an equivalence relation $E$ on $G$, and group $H$
let $H \rightarrow (G)^E_\sigma$ mean:
\begin{quotation}
{\em if\/} ${\bf d}$ is a function with domain $H$ and range of cardinality
$\le \sigma$,\\
{\em then\/} for some embedding $h$ of $G$ into $H$, for any $x,y\in G$
which are $E$-equivalent we have ${\bf d}(h(x))={\bf d}(h(y))$.
\end{quotation}
\item  In part (1) (or (2), or (3)) we can replace $A$ (or the domain of
$E$) by the family of subgroups of $G$ isomorphic to a fixed group $K$, and 
then ${\bf d}$ is a function with domain being the set of subgroups of $G$
isomorphic to $K$.
\item Like part (4), but we replace ``subgroups isomorphic to $K$'' by
``embedding of $K$'', and then replace ``$\rightarrow$'' by
``$\rightarrow^*$''. 
\end{enumerate}
\end{definition}

\begin{discussion}
\label{col.26}
There is a connection between the two last definitions: the first implies a
special case of the second one, when we restrict ourselves to permutation
groups of some finite set, and $A$ is the set of conjugates of the
permutation just interchanging two elements.
\end{discussion}

\begin{problem}
\label{8.22}
\begin{enumerate}
\item Investigate the arrows from Definition \ref{col.25}.
\item In particular, consider the case when $A$ is a set of pairwise
conjugate members of $G$ each of order two. 
\end{enumerate}
\end{problem}

\bigskip

\noindent{\sc Recent advances:}\\
During the winter of 1999, Gitik told me that he can start with 
\[\begin{array}{ll}
\bV\models&\mbox{`` $\kappa_n$ hypermeasurable of order $\lambda_n$},\\
& \lambda_n\mbox{ first (strongly) inaccessible }>\kappa_n,\ \lambda_n<
\kappa_{n+1},\ \lambda>\kappa=\sum\limits_n \kappa_n\mbox{ '',}
  \end{array}\]
and find a forcing notion $\bP$, not adding bounded subsets of $\kappa=
\sum\limits_{n<\omega}\kappa_n$, satisfying the $\kappa^{++}$--c.c., and
making $2^\kappa\ge\lambda$. I have conjectured that combining this proof
with earlier proof, you can demand $\bP$ makes $\lambda_n$ be the $n$--th
inaccessible cardinal, $\kappa=\sum\limits_n \lambda_n$, GCH holds below
$\kappa$ and $2^\kappa\ge\lambda$, so $\pp(\kappa)\ge\lambda$. Gitik has
confirmed this conjecture with $\kappa_n$ Mahlo. This proves that though
\ref{1.17} is open, other theorem which holds for $\lambda$ of cofinality
$\aleph_1$ cannot be generalized to cofinality $\aleph_0$. 

%\bibliographystyle{hplain}
%\bibliography{listb,lista,listx,liste,listf}

\begin{thebibliography}{100}

\bibitem{Ab83}
Uri Abraham.
\newblock {Aronszajn trees on $\aleph_2$ and $\aleph_3$}.
\newblock {\em Ann. Pure Appl. Logic}, 24:213--230, 1983.

\bibitem{ARSh:153}
Uri Abraham, Matatyahu Rubin, and Saharon Shelah.
\newblock {On the consistency of some partition theorems for continuous
  colorings, and the structure of $\aleph_ 1$-dense real order types}.
\newblock {\em {Annals of Pure and Applied Logic}}, 29:123--206, 1985.

\bibitem{AbSh:114}
Uri Abraham and Saharon Shelah.
\newblock {Isomorphism types of Aronszajn trees}.
\newblock {\em {Israel Journal of Mathematics}}, 50:75--113, 1985.

\bibitem{AbTo97}
Uri Abraham and Stevo Todor\v{c}evi\'{c}.
\newblock {Partition properties of $\omega_1$ compatible with CH}.
\newblock {\em {Fundamenta Mathematicae}}, 152:165--181, 1997.

\bibitem{ADSh:81}
Uri Avraham~(Abraham), Keith~J. Devlin, and Saharon Shelah.
\newblock {The consistency with CH of some consequences of Martin's axiom plus
  $2^{\aleph _{0}}>\aleph _{1}$}.
\newblock {\em {Israel Journal of Mathematics}}, 31:19--33, 1978.

\bibitem{BaJu95}
Tomek Bartoszy\'nski and Haim Judah.
\newblock {\em {Set Theory: On the Structure of the Real Line}}.
\newblock A K Peters, Wellesley, Massachusetts, 1995.

\bibitem{BRSh:616}
Tomek Bartoszy\'nski, Andrzej Ros{\l}anowski, and Saharon Shelah.
\newblock {After all, there are some inequalities which are provable in ZFC}.
\newblock {\em {Journal of Symbolic Logic}}, to appear,
math.LO/9711222\footnote{Papers with numbers of the form math.XX/xxxxxxx are
electronically available from the XXX Mathematics Archive at {\tt
http://front.math.ucdavis.edu/}}.  

\bibitem{BRSh:490}
Tomek Bartoszy\'nski, Andrzej Ros{\l}anowski, and Saharon Shelah.
\newblock {Adding one random real}.
\newblock {\em {Journal of Symbolic Logic}}, 61:80--90, 1996, math.LO/9406229.

\bibitem{BHT93}
James Baumgartner, Andras Hajnal, and Stevo Todor\v{c}evi\'{c}.
\newblock {Extensions of the Erdos--Rado Theorems}.
\newblock In {\em {Finite and Infinite Combinatorics in Set Theory and Logic}},
  pages 1--18. Kluwer Academic Publishers, 1993.
\newblock {N.W. Sauer et. al. eds}.

\bibitem{B}
James~E. Baumgartner.
\newblock {Decomposition of embedding of trees}.
\newblock {\em Notices Amer. Math. Soc.}, 17:967, 1970.

\bibitem{B4}
James~E. Baumgartner.
\newblock {All $\aleph_1 $-dense sets of reals can be isomorphic}.
\newblock {\em Fund. Math.}, 79:101--106, 1973.

\bibitem{B6}
James~E. Baumgartner.
\newblock {Ultrafilters on $\omega$}.
\newblock {\em The Journal of Symbolic Logic}, 60:624--639, 1995.

\bibitem{BD}
Shai Ben~David.
\newblock {On Shelah's compactness of cardinals}.
\newblock {\em Israel J. of Math.}, 31:34--56 and 394, 1978.

\bibitem{BsSh:242}
Andreas Blass and Saharon Shelah.
\newblock {There may be simple $P_ {\aleph_ 1}$- and $P_ {\aleph_ 2}$-points
  and the Rudin-Keisler ordering may be downward directed}.
\newblock {\em {Annals of Pure and Applied Logic}}, 33:213--243, 1987.

\bibitem{BsSh:257}
Andreas Blass and Saharon Shelah.
\newblock {Ultrafilters with small generating sets}.
\newblock {\em {Israel Journal of Mathematics}}, 65:259--271, 1989.

\bibitem{BoMo89}
Robert Bonnet and Donald Monk.
\newblock {\em {Handbook of Boolean Algebras, vol.~1--3 }}.
\newblock North Holland, 1989.

\bibitem{BnSh:642}
Joerg Brendle and Saharon Shelah.
\newblock {Ultrafilters on $\omega$ --- their ideals and their cardinal
  characteristics}.
\newblock {\em {Transactions of the American Mathematical Society}}, accepted,
  math.LO/9710217.

\bibitem{Ca93}
Timothy~J. Carlson.
\newblock {Strong measure zero and strongly meager sets}.
\newblock {\em Proceedings of the American Mathematical Society}, 118:577--586,
  1993.

\bibitem{Ci97}
Krzysztof Ciesielski.
\newblock Set theoretic real analysis.
\newblock {\em J. of Applied Analysis}, 3:143--190, 1997.

\bibitem{CiSh:695}
Krzysztof Ciesielski and Saharon Shelah.
\newblock {Category analog of sub-measurability problem}.
\newblock {\em Fundamenta Mathematicae}, submitted, math.LO/9905147.

\bibitem{CDSh:571}
James Cummings, Mirna D\v{z}amonja, and Saharon Shelah.
\newblock {A consistency result on weak reflection}.
\newblock {\em {Fundamenta Mathematicae}}, 148:91--100, 1995, math.LO/9504221.

\bibitem{CuFo98}
James Cummings and Matthew Foreman.
\newblock The tree property.
\newblock {\em Adv. Math.}, 133:1--32, 1998.

\bibitem{CuSh:530}
James Cummings and Saharon Shelah.
\newblock {A model in which every infinite Boolean algebra has many
  subalgebras}.
\newblock {\em {Journal of Symbolic Logic}}, 60:992--1004, 1995,
  math.LO/9509227.

\bibitem{CuSh:541}
James Cummings and Saharon Shelah.
\newblock {Cardinal invariants above the continuum}.
\newblock {\em {Annals of Pure and Applied Logic}}, 75:251--268, 1995,
  math.LO/9509228.

\bibitem{DvSh:65}
Keith~J. Devlin and Saharon Shelah.
\newblock {A weak version of $\diamondsuit $ which follows from $2^{\aleph
  _{0}}<2^{\aleph _{1}}$}.
\newblock {\em {Israel Journal of Mathematics}}, 29:239--247, 1978.

\bibitem{DJ2}
Tony Dodd and Ronald~B. Jensen.
\newblock {The covering lemma for $K$}.
\newblock {\em Ann. of Math Logic}, 22:1--30, 1982.

\bibitem{DjSh:659}
Mirna D\v{z}amonja and Saharon Shelah.
\newblock {Universal graphs at successors of singular strong limits}.
\newblock {\em {preprint}}.

\bibitem{DjSh:691}
Mirna D\v{z}amonja and Saharon Shelah.
\newblock {Weak reflection can first happen at the successor of singular}.
\newblock {\em Preprint}.

\bibitem{DjSh:562}
Mirna D\v{z}amonja and Saharon Shelah.
\newblock {On squares, outside guessing of clubs and $I_{<f}[\lambda]$}.
\newblock {\em {Fundamenta Mathematicae}}, 148:165--198, 1995, math.LO/9510216.

\bibitem{DjSh:545}
Mirna D\v{z}amonja and Saharon Shelah.
\newblock {Saturated filters at successors of singulars, weak reflection and
  yet another weak club principle}.
\newblock {\em {Annals of Pure and Applied Logic}}, 79:289--316, 1996,
  math.LO/9601219.

\bibitem{EM}
Paul~C. Eklof and Alan Mekler.
\newblock {\em {Almost free modules; Set theoretic methods}}.
\newblock North Holland Library, 1990.

\bibitem{EH}
Paul Erd\H{o}s and Andras Hajnal.
\newblock {Unsolved Problems in set theory}.
\newblock In {\em Axiomatic Set Theory}, volume XIII Part I of {\em Proc. of
  Symp. in Pure Math.}, pages 17--48, Providence, R.I., 1971. AMS.

\bibitem{EHMR}
Paul Erd\H{o}s, Andras Hajnal, A.~Mat\'e, and Richard Rado.
\newblock {\em {Combinatorial set theory: Partition Relations for Cardinals}},
  volume 106 of {\em Studies in Logic and the Foundation of Math.}
\newblock North Holland Publ. Co, Amsterdam, 1984.

\bibitem{FlSh:267}
William~G. Fleissner and Saharon Shelah.
\newblock {Collectionwise Hausdorff: incompactness at singulars}.
\newblock {\em {Topology and its Applications}}, 31:101--107, 1989.

\bibitem{FW}
Matthew Foreman and Hugh Woodin.
\newblock {The generalized continuum hypothesis can fail everywhere}.
\newblock {\em Annals Math.}, 133:1--36, 1991.

\bibitem{Fe94}
David Fremlin.
\newblock Problem list.
\newblock circulated notes (1994).

\bibitem{Fe93}
David Fremlin.
\newblock {Real--valued--measurable cardinals}.
\newblock In {\em Set Theory of the Reals}, volume~6 of {\em Israel
  Mathematical Conference Proceedings}, pages 151--304, 1993.

\bibitem{Fu}
Laszlo Fuchs.
\newblock {\em {Infinite Abelian Groups}}, volume I, II.
\newblock Academic Press, New York, 1970, 1973.

\bibitem{GFJ97}
Salvador Garcia-Ferreira and Winfried Just.
\newblock {Two examples of relatively pseudocompact spaces}.
\newblock {\em Questions and Answers in General Topology}, 17:35--45, 1999.

\bibitem{Gi80}
Moti Gitik.
\newblock {All uncountable cardinals can be singular}.
\newblock {\em Israel Journal of Mathematics}, 35:61--88, 1980.

\bibitem{GiSh:582}
Moti Gitik and Saharon Shelah.
\newblock {More on real-valued measurable cardinals and forcing with ideals}.
\newblock {\em {Israel Journal of Mathematics}}, submitted, math.LO/9507208.

\bibitem{GiSh:357}
Moti Gitik and Saharon Shelah.
\newblock {Forcings with ideals and simple forcing notions}.
\newblock {\em {Israel Journal of Mathematics}}, 68:129--160, 1989.

\bibitem{GiSh:344}
Moti Gitik and Saharon Shelah.
\newblock {On certain indestructibility of strong cardinals and a question of
  Hajnal}.
\newblock {\em {Archive for Mathematical Logic}}, 28:35--42, 1989.

\bibitem{GiSh:412}
Moti Gitik and Saharon Shelah.
\newblock {More on simple forcing notions and forcings with ideals}.
\newblock {\em {Annals of Pure and Applied Logic}}, 59:219--238, 1993.

\bibitem{GiSh:577}
Moti Gitik and Saharon Shelah.
\newblock {Less saturated ideals}.
\newblock {\em {Proceedings of the American Mathematical Society}},
  125:1523--1520, 1997, math.LO/9503203.

\bibitem{GoSh:448}
Martin Goldstern and Saharon Shelah.
\newblock {Many simple cardinal invariants}.
\newblock {\em {Archive for Mathematical Logic}}, 32:203--221, 1993,
  math.LO/9205208.

\bibitem{GrRoSp}
Ronald Graham, Bruce~L. Rothschild, and Joel Spencer.
\newblock {\em {Ramsey Theory}}.
\newblock {Willey -- Interscience Series in Discrete Mathematics}. Willey, New
  York, 1980.

\bibitem{Gre}
John Gregory.
\newblock {Higher Souslin trees and the generalized continuum hypothesis}.
\newblock {\em Journal of Symbolic Logic}, 41(3):663--671, 1976.

\bibitem{GrSh:302}
Rami Grossberg and Saharon Shelah.
\newblock {On the structure of ${\rm Ext}_ p(G,{\bf Z})$}.
\newblock {\em {Journal of Algebra}}, 121:117--128, 1989.
\newblock See also [GrSh:302a] below.

\bibitem{GrSh:302a}
Rami Grossberg and Saharon Shelah.
\newblock {On cardinalities in quotients of inverse limits of groups}.
\newblock {\em {Mathematica Japonica}}, 47(2), 1998.

\bibitem{Ha93}
Andras Hajnal.
\newblock {True embedding partition relations}.
\newblock In {\em {Finite and Infinite Combinatorics in Sets and Logic}}, pages
  135--152. Kluwer Academic Publishers, 1993.

\bibitem{HH}
Andras Hajnal and Peter Hamburger.
\newblock {\em {Halmazelm\'elet (Set Theory)}}.
\newblock Tank\"onyvkia\'o V\'allalat, Budapest, Hungary.

\bibitem{HaJuSz}
Andras Hajnal, Istvan Juh\'asz, and Z.~Szentmiklossy.
\newblock {On the structure of CCC partial orders}.
\newblock {\em Algebra Universalis}, to appear.

\bibitem{Hm98}
Joel~D. Hamkins.
\newblock {Every group has a terminating transfinite automorphism tower}.
\newblock {\em Proc. Amer. Math. Soc.}, 126:3223--3226, 1998.

\bibitem{Ho81}
Wilfrid Hodges.
\newblock {For singular $\lambda$, $\lambda$-free, implies free.}
\newblock {\em Algebra Universalis}, 12:205--220, 1981.

\bibitem{JdSh:292}
Jaime Ihoda (Haim~Judah) and Saharon Shelah.
\newblock {Souslin forcing}.
\newblock {\em {The Journal of Symbolic Logic}}, 53:1188--1207, 1988.

\bibitem{JeSh:476}
Thomas Jech and Saharon Shelah.
\newblock {Possible pcf algebras}.
\newblock {\em {Journal of Symbolic Logic}}, 61:313--317, 1996,
  math.LO/9412208.

\bibitem{JdSh:335}
Haim Judah and Saharon Shelah.
\newblock {${\rm MA}(\sigma$-centered): Cohen reals, strong measure zero sets
  and strongly meager sets}.
\newblock {\em {Israel Journal of Mathematics}}, 68:1--17, 1989.

\bibitem{JdSh:446}
Haim Judah and Saharon Shelah.
\newblock {Baire Property and Axiom of Choice}.
\newblock {\em {Israel Journal of Mathematics}}, 84:435--450, 1993,
  math.LO/9211213.

\bibitem{Ju}
Istvan Juh\'asz.
\newblock {Cardinal functions. II.}
\newblock In {\em {Handbook of set-theoretic topology}}, pages 63--109. North
  Holland Publ. Co, 1984.

\bibitem{Juh92}
Istvan Juh\'asz.
\newblock {Cardinal functions}.
\newblock In Miroslav~Hu\v sek and Jan van Mill, editors, {\em {Recent progress
  in general topology (Prague, 1991)}}, pages 417--441. North-Holland
  Publishing Co., Amsterdam, 1992.

\bibitem{JuSh:231}
Istvan Juh\'asz and Saharon Shelah.
\newblock {How large can a hereditarily separable or hereditarily Lindelof
  space be?}
\newblock {\em {Israel Journal of Mathematics}}, 53:355--364, 1986.

\bibitem{JMPS}
Winfried Just, A.~R.~D. Mathias, Karel Prikry, and Petr Simon.
\newblock {On the existence of large p--ideals}.
\newblock {\em Journal of Symbolic Logic}, 55(2):457--465, 1990.

\bibitem{JsShTh:654}
Winfried Just, Saharon Shelah, and Simon Thomas.
\newblock {The Automorphism Tower Problem III: Closed Groups of Uncountable
  Degree}.
\newblock {\em {Advances in Mathematics}}, accepted.

\bibitem{KeSo95}
Alexander~S. Kechris and S{\l}awomir Solecki.
\newblock {Approximation of analytic by Borel sets and definable countable
  chain conditions}.
\newblock {\em Israel Journal of Mathematics}, 89:343--356, 1995.

\bibitem{Kt76}
Jussi Ketonen.
\newblock {On the existence of $P$-points in the Stone--\v Cech
  compactification of integers}.
\newblock {\em Fund. Math.}, 92:91--94, 1976.

\bibitem{KjSh:409}
Menachem Kojman and Saharon Shelah.
\newblock {Non-existence of Universal Orders in Many Cardinals}.
\newblock {\em {Journal of Symbolic Logic}}, 57:875--891, 1992,
  math.LO/9209201.

\bibitem{KjSh:447}
Menachem Kojman and Saharon Shelah.
\newblock {The universality spectrum of stable unsuperstable theories}.
\newblock {\em {Annals of Pure and Applied Logic}}, 58:57--72, 1992,
  math.LO/9201253.

\bibitem{Ko}
Peter Komjath.
\newblock On second-category sets.
\newblock {\em Proc. Amer. Math. Soc.}, 107:653--654, 1989.

\bibitem{Ku77}
David~W. Kueker.
\newblock {Countable approximations and Lowenheim-Skolem theorems}.
\newblock {\em Ann. Math. Logic}, 11:57--103, 1977.

\bibitem{Ku84}
Kenneth Kunen.
\newblock {Random and Cohen Reals}.
\newblock In K.~Kunen and J.~E. Vaughan, editors, {\em {Handbook of
  Set--Theoretic Topology}}, pages 887--911. Elsevier Science Publishers B.V.,
  1984.

\bibitem{Ku90}
Kenneth Kunen.
\newblock {Large homogeneous compact spaces}.
\newblock In J.~van Mill and G.M. Reed, editors, {\em {Open Problems in
  Topology}}. Elsevier Science Publishers B.V, 1990.

\bibitem{L1}
Richard Laver.
\newblock {On the consistency of Borel's conjecture}.
\newblock {\em Acta Math.}, 137:151--169, 1976.

\bibitem{Mt78}
A.R.D. Mathias.
\newblock {$0^\#$ and the p-point problem}.
\newblock In {\em {Higher Set Theory (Proc. Conf., Math. Forschungsinst.,
  Oberwolfach, 1977)}}, volume 669 of {\em Lecture Notes in Math.}, pages
  375--383. Springer, Berlin, 1978.

\bibitem{MRSh:314}
Alan~H. Mekler, Andrzej Ros{\l}anowski, and Saharon Shelah.
\newblock {On the $p$-rank of Ext}.
\newblock {\em {Israel Journal of Mathematics}}, to appear, math.LO/9806165.

\bibitem{MkSh:274}
Alan~H. Mekler and Saharon Shelah.
\newblock {Uniformization principles}.
\newblock {\em {The Journal of Symbolic Logic}}, 54:441--459, 1989.

\bibitem{MkSh:366}
Alan~H. Mekler and Saharon Shelah.
\newblock {Almost free algebras }.
\newblock {\em {Israel Journal of Mathematics}}, 89:237--259, 1995,
  math.LO/9408213.

\bibitem{MShS:417}
Alan~H. Mekler, Saharon Shelah, and Otmar Spinas.
\newblock {The essentially free spectrum of a variety}.
\newblock {\em Israel Journal of Mathematics}, 93:1--8, 1996, math.LO/9411234.

\bibitem{Mi91}
Arnold~W. Miller.
\newblock {Arnie Miller's problem list}.
\newblock In Haim Judah, editor, {\em {Set Theory of the Reals}}, volume~6 of
  {\em Israel Mathematical Conference Proceedings}, pages 645--654.
\newblock {Proceedings of the Winter Institute held at Bar--Ilan University,
  Ramat Gan, January 1991}.

\bibitem{M2}
Donald Monk.
\newblock {\em {Cardinal Invariants of Boolean Algebras}}, volume 142 of {\em
  {Progress in Mathematics}}.
\newblock Birkh\"auser Verlag, Basel--Boston--Berlin, 1996.

\bibitem{Rt79}
Judy Roitman.
\newblock {Adding a random or Cohen real: topological consequences and the
  effect on Martin's axiom}.
\newblock {\em Fundamenta Mathematicae}, 103:47--60, 1979.

\bibitem{RoSh:534}
Andrzej Ros{\l}anowski and Saharon Shelah.
\newblock {Cardinal invariants of ultrapoducts of Boolean algebras}.
\newblock {\em {Fundamenta Mathematicae}}, 155:101--151, 1998, math.LO/9703218.

\bibitem{RoSh:651}
Andrzej Ros{\l}anowski and Saharon Shelah.
\newblock {Forcing for hL, hd and Depth}.
\newblock {\em {Colloquium Mathematicum}}, accepted, math.LO/9808104.

\bibitem{RoSh:599}
Andrzej Ros{\l}anowski and Saharon Shelah.
\newblock {More on cardinal functions on Boolean algebras}.
\newblock {\em {Annals of Pure and Applied Logic}}, accepted, math.LO/9808056.

\bibitem{RoSh:470}
Andrzej Ros{\l}anowski and Saharon Shelah.
\newblock {Norms on possibilities I: forcing with trees and creatures}.
\newblock {\em {Memoirs of the AMS}}, to appear, math.LO/9807172.

\bibitem{RoSh:628}
Andrzej Ros{\l}anowski and Saharon Shelah.
\newblock {Norms on possibilities II: more ccc ideals on
  $2^{\textstyle\omega}$}.
\newblock {\em {Journal of Applied Analysis}}, 3:103--127, 1997,
  math.LO/9703222.

\bibitem{RoSh:670}
Andrzej Ros{\l}anowski and Saharon Shelah.
\newblock {Norms on possibilities III: strange subsets of the real line}.
\newblock {\em {in preparation}}.

\bibitem{RoSh:672}
Andrzej Ros{\l}anowski and Saharon Shelah.
\newblock {Norms on possibilities IV: ccc forcing notions}.
\newblock {\em {Preprint}}.

\bibitem{RuSh:117}
Matatyahu Rubin and Saharon Shelah.
\newblock {Combinatorial problems on trees: partitions, $\Delta$-systems and
  large free subtrees}.
\newblock {\em {Annals of Pure and Applied Logic}}, 33:43--81, 1987.

\bibitem{SgSh:217}
Gershon Sageev and Saharon Shelah.
\newblock {Noetherian ring with free additive groups}.
\newblock {\em {Abstracts of the American Mathematical Society}}, 7:369, 1986.

\bibitem{ST94}
James~D. Sharp and Simon Thomas.
\newblock {Uniformization problems and the cofinality of the infinite symmetric
  group}.
\newblock {\em Notre Dame J. Formal Logic}, 35:328--345, 1994.

\bibitem{ST95}
James~D. Sharp and Simon Thomas.
\newblock {Unbounded families and the cofinality of the infinite symmetric
  group}.
\newblock {\em Arch. Math. Logic}, 34:33--45, 1995.

\bibitem{Sh:E9}
Saharon Shelah.
\newblock {Remarks on $\aleph_1$--CWH not CWH first countable spaces}.
\newblock In {\em Set Theory, Boise ID, 1992--1994}, volume 192 of {\em
  Contemporary Mathematics}, pages 103--145.

\bibitem{Sh:E11}
Saharon Shelah.
\newblock {Also quite large ${\frak b}\subseteq \pcf({\frak a})$ behave
  nicely}.
\newblock math.LO/9906018.

\bibitem{Sh:E12}
Saharon Shelah.
\newblock {Analytical Guide and Corrections to \cite{Sh:g}.}
\newblock math.LO/9906022.

\bibitem{Sh:589}
Saharon Shelah.
\newblock {Applications of PCF theory}.
\newblock {\em {Journal of Symbolic Logic}}, submitted, math.LO/9804155.

\bibitem{Sh:576}
Saharon Shelah.
\newblock {Categoricity of an abstract elementary class in two successive
  cardinals}.
\newblock {\em {Israel Journal of Mathematics}}, submitted, math.LO/9805146.

\bibitem{Sh:575}
Saharon Shelah.
\newblock {Cellularity of free products of Boolean algebras (or topologies)}.
\newblock {\em {Fundamenta Mathematica}}, accepted, math.LO/9508221.

\bibitem{Sh:641}
Saharon Shelah.
\newblock {Constructing Boolean algebras for cardinal invariants}.
\newblock {\em {Algebra Universalis}}, {submitted}, math.LO/9712286.

\bibitem{Sh:592}
Saharon Shelah.
\newblock {Covering of the null ideal may have countable cofinality}.
\newblock {\em {Fundamenta Mathematicae}}, submitted, math.LO/9810181.

\bibitem{Sh:603}
Saharon Shelah.
\newblock {Few non-minimal types and non-structure}.
\newblock {\em {in preparation}}, math.LO/9906023.

\bibitem{Sh:E16}
Saharon Shelah.
\newblock The future of set theory.
\newblock In Haim Judah, editor, {\em {Set Theory of the Reals}}, volume~6 of
  {\em Israel Mathematical Conference Proceedings}.
\newblock {Proceedings of the Winter Institute held at Bar--Ilan University,
  Ramat Gan, January 1991}.

\bibitem{Sh:538}
Saharon Shelah.
\newblock {Historic iteration with $\aleph_\varepsilon$-support}.
\newblock {\em Archive for Mathematical Logic}, accepted, math.LO/9607227.

\bibitem{Sh:655}
Saharon Shelah.
\newblock {Iteration of $\lambda$-complete forcing notions not collapsing
  $\lambda^+$.}
\newblock {\em {in preparation}}, math.LO/9906024.

\bibitem{Sh:638}
Saharon Shelah.
\newblock {More on Weak Diamond}.
\newblock {\em {Preprint}}, math.LO/9807180.

\bibitem{Sh:e}
Saharon Shelah.
\newblock {\em {Non--structure theory}}, volume accepted.
\newblock {Oxford University Press}.

\bibitem{Sh:630}
Saharon Shelah.
\newblock {Non-elementary proper forcing notions}.
\newblock {\em {Journal of Applied Analysis}}, submitted, math.LO/9712283.

\bibitem{Sh:669}
Saharon Shelah.
\newblock {Non-elementary proper forcing notions II}.
\newblock {\em {In preparation}}.

\bibitem{Sh:622}
Saharon Shelah.
\newblock {Non-existence of universal members in classes of Abelian groups}.
\newblock {\em {Journal of Group Theory}}, submitted, math.LO/9808139.

\bibitem{Sh:587}
Saharon Shelah.
\newblock {Not collapsing cardinals $\leq\kappa$ in $(<\kappa)$--support
  iterations}.
\newblock {\em {Israel Journal of Mathematics}}, accepted, math.LO/9707225.

\bibitem{Sh:667}
Saharon Shelah.
\newblock {Not collapsing cardinals $\leq\kappa$ in $(<\kappa)$--support
  iterations II}.
\newblock {\em {Israel Journal of Mathematics}}, {submitted}, math.LO/9808140.

\bibitem{Sh:668}
Saharon Shelah.
\newblock {On Arhangelskii's Problem}.
\newblock {\em {in preparation}}, math.LO/9906025.

\bibitem{Sh:656}
Saharon Shelah.
\newblock {On CS iterations not adding reals}.
\newblock {\em {in preparation}}.

\bibitem{Sh:513}
Saharon Shelah.
\newblock {PCF and infinite free subsets}.
\newblock {\em {Archive for Mathematical Logic}}, accepted, math.LO/9807177.

\bibitem{Sh:295}
Saharon Shelah.
\newblock {Projective measurability does not imply projective Baire}.
\newblock {\em {in preparation}}.

\bibitem{Sh:620}
Saharon Shelah.
\newblock {Special Subsets of ${}^{{\rm cf}(\mu)}\mu$, Boolean Algebras and
  Maharam measure Algebras}.
\newblock {\em {General Topology and its Applications - Proc. of Prague
  Topological Symposium 1996}}, accepted, math.LO/9804156.

\bibitem{Sh:580}
Saharon Shelah.
\newblock {Strong covering without squares}.
\newblock {\em {Fundamenta Mathematicae}}, accepted, math.LO/9604243.

\bibitem{Sh:664}
Saharon Shelah.
\newblock {Strong dichotomy of cardinality}.
\newblock {\em {Results in Mathematics}}, {submitted}, math.LO/9807183.

\bibitem{Sh:460}
Saharon Shelah.
\newblock {The Generalized Continuum Hypothesis revisited}.
\newblock {\em {Israel Journal of Mathematics}}, accepted, math.LO/9809200.

\bibitem{Sh:619}
Saharon Shelah.
\newblock {The null ideal restricted to a non-null set may be saturated}.
\newblock {\em {Asian Journal of Mathematics}}, submitted, math.LO/9705213.

\bibitem{Sh:546}
Saharon Shelah.
\newblock {Was Sierpi\'nski right? IV}.
\newblock {\em {Journal of Symbolic Logic}}, accepted, math.LO/9712282.

\bibitem{Sh:26}
Saharon Shelah.
\newblock {Notes on combinatorial set theory}.
\newblock {\em {Israel Journal of Mathematics}}, 14:262--277, 1973.

\bibitem{Sh:52}
Saharon Shelah.
\newblock {A compactness theorem for singular cardinals, free algebras,
  Whitehead problem and transversals}.
\newblock {\em {Israel Journal of Mathematics}}, 21:319--349, 1975.

\bibitem{Sh:50}
Saharon Shelah.
\newblock {Decomposing uncountable squares to countably many chains}.
\newblock {\em {Journal of Combinatorial Theory. Ser. A}}, 21:110--114, 1976.

\bibitem{Sh:E18}
Saharon Shelah.
\newblock {A combinatorial proof of the singular compactness theorem}.
\newblock {\em Mineograph notes and lecture in a mini-conference, Berlin,
  August'78}, 1977.

\bibitem{Sh:64}
Saharon Shelah.
\newblock {Whitehead groups may be not free, even assuming CH. I}.
\newblock {\em {Israel Journal of Mathematics}}, 28:193--204, 1977.

\bibitem{Sh:a}
Saharon Shelah.
\newblock {\em {Classification theory and the number of nonisomorphic models}},
  volume~92 of {\em {Studies in Logic and the Foundations of Mathematics}}.
\newblock {North-Holland Publishing Co., Amsterdam-New York, xvi+544 pp,
  \$62.25}, 1978.

\bibitem{Sh:108}
Saharon Shelah.
\newblock {On successors of singular cardinals}.
\newblock In {\em {Logic Colloquium '78 (Mons, 1978)}}, volume~97 of {\em
  {Stud. Logic Foundations Math}}, pages 357--380. {North-Holland,
  Amsterdam-New York}, 1979.

\bibitem{Sh:92}
Saharon Shelah.
\newblock {Remarks on Boolean algebras}.
\newblock {\em {Algebra Universalis}}, 11:77--89, 1980.

\bibitem{Sh:95}
Saharon Shelah.
\newblock {Canonization theorems and applications}.
\newblock {\em {The Journal of Symbolic Logic}}, 46:345--353, 1981.

\bibitem{Sh:126}
Saharon Shelah.
\newblock {On saturation for a predicate}.
\newblock {\em {Notre Dame Journal of Formal Logic}}, 22:239--248, 1981.

\bibitem{Sh:b}
Saharon Shelah.
\newblock {\em {Proper forcing}}, volume 940 of {\em {Lecture Notes in
  Mathematics}}.
\newblock {Springer-Verlag, Berlin-New York, xxix+496 pp}, 1982.

\bibitem{Sh:87b}
Saharon Shelah.
\newblock {Classification theory for nonelementary classes, I. The number of
  uncountable models of $\psi \in L_{\omega _{1},\omega }$. Part B}.
\newblock {\em {Israel Journal of Mathematics}}, 46:241--273, 1983.

\bibitem{Sh:136}
Saharon Shelah.
\newblock {Constructions of many complicated uncountable structures and Boolean
  algebras}.
\newblock {\em {Israel Journal of Mathematics}}, 45:100--146, 1983.

\bibitem{Sh:185}
Saharon Shelah.
\newblock {Lifting problem of the measure algebra}.
\newblock {\em {Israel Journal of Mathematics}}, 45:90--96, 1983.

\bibitem{Sh:176}
Saharon Shelah.
\newblock {Can you take Solovay's inaccessible away?}
\newblock {\em {Israel Journal of Mathematics}}, 48:1--47, 1984.

\bibitem{Sh:186}
Saharon Shelah.
\newblock {Diamonds, uniformization}.
\newblock {\em {The Journal of Symbolic Logic}}, 49:1022--1033, 1984.

\bibitem{Sh:177}
Saharon Shelah.
\newblock {More on proper forcing}.
\newblock {\em {The Journal of Symbolic Logic}}, 49:1034--1038, 1984.

\bibitem{Sh:208}
Saharon Shelah.
\newblock {More on the weak diamond}.
\newblock {\em {Annals of Pure and Applied Logic}}, 28:315--318, 1985.

\bibitem{Sh:233}
Saharon Shelah.
\newblock {Remarks on the numbers of ideals of Boolean algebra and open 
  sets of a topology}.
\newblock In {\em {Around classification theory of models}}, volume 1182 of
  {\em {Lecture Notes in Mathematics}}, pages 151--187. {Springer, Berlin},
  1986.

\bibitem{Sh:88a}
Saharon Shelah.
\newblock {Appendix: on stationary sets (in ``Classification of nonelementary
  classes. II. Abstract elementary classes'')}.
\newblock In {\em {Classification theory (Chicago, IL, 1985)}}, volume 1292 of
  {\em {Lecture Notes in Mathematics}}, pages 483--495. {Springer, Berlin},
  1987.
\newblock {Proceedings of the USA--Israel Conference on Classification Theory,
  Chicago, December 1985; ed. Baldwin, J.T.}

\bibitem{Sh:273}
Saharon Shelah.
\newblock {Can the fundamental (homotopy) group of a space be the rationals?}
\newblock {\em {Proceedings of the American Mathematical Society}},
  103:627--632, 1988.

\bibitem{Sh:329}
Saharon Shelah.
\newblock {Primitive recursive bounds for van der Waerden numbers}.
\newblock {\em {Journal of the American Mathematical Society}}, 1:683--697,
  1988.

\bibitem{Sh:289}
Saharon Shelah.
\newblock {Consistency of positive partition theorems for graphs and models}.
\newblock In {\em {Set theory and its applications (Toronto, ON, 1987)}},
  volume 1401 of {\em {Lecture Notes in Mathematics}}, pages 167--193.
  {Springer, Berlin-New York}, 1989.
\newblock {ed. Steprans, J. and Watson, S.}

\bibitem{Sh:262}
Saharon Shelah.
\newblock {The number of pairwise non-elementarily-embeddable models}.
\newblock {\em {The Journal of Symbolic Logic}}, 54:1431--1455, 1989.

\bibitem{Sh:347}
Saharon Shelah.
\newblock {Incompactness for chromatic numbers of graphs}.
\newblock In {\em {A tribute to Paul Erd\H{o}s}}, pages 361--371. {Cambridge
  Univ. Press, Cambridge}, 1990.

\bibitem{Sh:280}
Saharon Shelah.
\newblock {Strong negative partition above the continuum}.
\newblock {\em {The Journal of Symbolic Logic}}, 55:21--31, 1990.

\bibitem{Sh:351}
Saharon Shelah.
\newblock {Reflecting stationary sets and successors of singular cardinals}.
\newblock {\em {Archive for Mathematical Logic}}, 31:25--53, 1991.

\bibitem{Sh:288}
Saharon Shelah.
\newblock {Strong Partition Relations Below the Power Set: Consistency, Was
  Sierpi\'nski Right, II?}
\newblock In {\em {Proceedings of the Conference on Set Theory and its
  Applications in honor of A.Hajnal and V.T.Sos, Budapest, 1/91}}, volume~60 of
  {\em Colloquia Mathematica Societatis Janos Bolyai. Sets, Graphs, and
  Numbers}, pages 637--638. 1991, math.LO/9201244.

\bibitem{Sh:400a}
Saharon Shelah.
\newblock {Cardinal arithmetic for skeptics}.
\newblock {\em {American Mathematical Society. Bulletin. New Series}},
  26:197--210, 1992, math.LO/9201251.

\bibitem{Sh:326}
Saharon Shelah.
\newblock {Viva la difference I: Nonisomorphism of ultrapowers of countable
  models}.
\newblock In {\em {Set Theory of the Continuum}}, volume~26 of {\em
  {Mathematical Sciences Research Institute Publications}}, pages 357--405.
  {Springer Verlag}, 1992, math.LO/9201245.

\bibitem{Sh:420}
Saharon Shelah.
\newblock {Advances in Cardinal Arithmetic}.
\newblock In {\em {Finite and Infinite Combinatorics in Sets and Logic}}, pages
  355--383. Kluwer Academic Publishers, 1993.
\newblock {N.W. Sauer et al (eds.)}.

\bibitem{Sh:424}
Saharon Shelah.
\newblock {On $CH + 2^{\aleph_1}\rightarrow(\alpha)^2_2$ for
  $\alpha<\omega_2$}.
\newblock In {\em {Logic Colloquium'90. ASL Summer Meeting in Helsinki}},
  volume~2 of {\em Lecture Notes in Logic}, pages 281--289. Springer Verlag,
  1993, math.LO/9308212.

\bibitem{Sh:g}
Saharon Shelah.
\newblock {\em {Cardinal Arithmetic}}, volume~29 of {\em {Oxford Logic
  Guides}}.
\newblock {Oxford University Press}, 1994.

\bibitem{Sh:454a}
Saharon Shelah.
\newblock {Cardinalities of topologies with small base}.
\newblock {\em {Annals of Pure and Applied Logic}}, 68:95--113, 1994,
  math.LO/9403219.

\bibitem{Sh:480}
Saharon Shelah.
\newblock {How special are Cohen and random forcings i.e. Boolean algebras of
  the family of subsets of reals modulo meagre or null}.
\newblock {\em {Israel Journal of Mathematics}}, 88:159--174, 1994,
  math.LO/9303208.

\bibitem{Sh:473}
Saharon Shelah.
\newblock {Possibly every real function is continuous on a non--meagre set}.
\newblock {\em {Publications de L'Institute Math\'ematique - Beograd, Nouvelle
  S\'erie}}, 57(71):47--60, 1995, math.LO/9511220.

\bibitem{Sh:430}
Saharon Shelah.
\newblock {Further cardinal arithmetic}.
\newblock {\em {Israel Journal of Mathematics}}, 95:61--114, 1996,
  math.LO/9610226.

\bibitem{Sh:542}
Saharon Shelah.
\newblock {Large Normal Ideals Concentrating on a Fixed Small Cardinality}.
\newblock {\em {Archive for Mathematical Logic}}, 35:341--347, 1996,
  math.LO/9406219.

\bibitem{Sh:572}
Saharon Shelah.
\newblock {Colouring and non-productivity of $\aleph_2$-cc}.
\newblock {\em {Annals of Pure and Applied Logic}}, 84:153--174, 1997,
  math.LO/9609218.

\bibitem{Sh:552}
Saharon Shelah.
\newblock {Non existence of universals for classes like reduced torsion free
  abelian groups under non neccessarily pure embeddings}.
\newblock In {\em Advances in Algebra and Model Theory. Editors: Manfred Droste
  und Ruediger Goebel}, volume~9, pages 229--286. Gordon and Breach, 1997,
  math.LO/9609217.

\bibitem{Sh:462}
Saharon Shelah.
\newblock {On $\sigma $-entangled linear orders}.
\newblock {\em {Fundamenta Mathematicae}}, 153:199--275, 1997, math.LO/9609216.

\bibitem{Sh:497}
Saharon Shelah.
\newblock {Set Theory without choice: not everything on cofinality is
  possible}.
\newblock {\em Archive for Mathematical Logic}, 36:81--125, 1997,
  math.LO/9512227.
\newblock {A special volume dedicated to Prof. Azriel Levy}.

\bibitem{Sh:f}
Saharon Shelah.
\newblock {\em {Proper and improper forcing}}.
\newblock {Perspectives in Mathematical Logic}. {Springer}, 1998.

\bibitem{Sh:594}
Saharon Shelah.
\newblock {There may be no nowhere dense ultrafilter}.
\newblock In {Makowsky, J.A. and Ravve, E.V.}, editor, {\em Proceedings of the
  Logic Colloquium Haifa'95}, volume~11 of {\em {Lecture Notes in Logic}},
  pages {305--324}. {Springer}, 1998, math.LO/9611221.

\bibitem{ShST:584}
Saharon Shelah, Jan Saxl, and Simon Thomas.
\newblock {Infinite products of finite simple groups}.
\newblock {\em {Transactions of the American Mathematical Society}},
  348:4611--4641, 1996, math.IG/9605202.

\bibitem{ShSt:608}
Saharon Shelah and Lee Stanley.
\newblock {Consistency of partition relation for cardinal in $(\lambda
  ,2^\lambda )$}.
\newblock {\em {Journal of Symbolic Logic}}, submitted, math.LO/9710216.

\bibitem{ShSt:258}
Saharon Shelah and Lee Stanley.
\newblock {A theorem and some consistency results in partition calculus}.
\newblock {\em {Annals of Pure and Applied Logic}}, 36:119--152, 1987.

\bibitem{ShSr:427}
Saharon Shelah and Juris Steprans.
\newblock {Somewhere trivial automorphisms}.
\newblock {\em {Journal of the London Mathematical Society}}, 49:569--580,
  1994, math.LO/9308214.

\bibitem{ShTh:524}
Saharon Shelah and Simon Thomas.
\newblock {The Cofinality Spectrum of The Infinite Symmetric Group}.
\newblock {\em {Journal of Symbolic Logic}}, 62:902--916, 1997,
  math.LO/9412230.

\bibitem{ShZa:610}
Saharon Shelah and Jindra Zapletal.
\newblock {Canonical models for $\aleph_1$ combinatorics}.
\newblock {\em {Transactions of the American Mathematical Society}}, accepted,
  math.LO/9806166.

\bibitem{Sh:F257}
{Shelah, Saharon}.
\newblock {On ``$I$ is a $\kappa$-complete ideal on $\kappa$ such that 
${\mathcal P}(\kappa)/I$ is the universal meager forcing''}. 

\bibitem{So98}
S{\l}awomir Solecki.
\newblock {Analytic ideals and their applications}.
\newblock {\em Annals of Pure and Applied Logic}, to appear.

\bibitem{So96}
S{\l}awomir Solecki.
\newblock {Analytic ideals}.
\newblock {\em Bulletin of Symbolic Logic}, 2:339--348, 1996.

\bibitem{So2}
Robert~M. Solovay.
\newblock {Real valued measurable cardinals}.
\newblock In {\em {Axiomatic set theory}}, volume XIII, Part 1 of {\em Proc. of
  Symposia in Pure Math.}, pages 397--428, Providence R.I, 1971. A.M.S.

\bibitem{SW87}
Juris Stepr\=ans and Stephen~W. Watson.
\newblock {Homeomorphisms of manifolds with prescribed behaviour on large dense
  sets}.
\newblock {\em Bull. London Math. Soc.}, 19:305--310, 1987.

\bibitem{Th96}
Simon Thomas.
\newblock {The cofinalities of the infinite-dimensional classical groups}.
\newblock {\em J. Algebra}, 179:704--719, 1996.

\bibitem{To89}
Stevo Todor\v{c}evi\'{c}.
\newblock {\em {Partition problems in Topology}}, volume~84 of {\em
  Contemporary Mathematics}.
\newblock American Mathematical Society, Providence, RI, 1989.

\bibitem{vM84}
Jan van Mill.
\newblock {An Introduction to $\beta\omega$}.
\newblock In K.~Kunen and J.~E. Vaughan, editors, {\em {Handbook of
  Set--Theoretic Topology}}, pages 503--567. Elsevier Science Publishers,
  Amsterdam, 1984.

\bibitem{Ve92}
Boban Veli{\v c}kovi\'c.
\newblock {Applications of the open coloring axiom}.
\newblock In H.~Judah, W.~Just, and H.~Woodin, editors, {\em {Set Theory of the
  Continuum (Berkeley, CA, 1989)}}, volume~26 of {\em Mathematical Sciences
  Research Institute Publications}, pages 137--154. Springer Verlag, 1992.

\bibitem{Wl84}
Scott~W. Williams.
\newblock {Box Products}.
\newblock In K.~Kunen and J.~E. Vaughan, editors, {\em {Handbook of
  Set--Theoretic Topology}}, pages 169--200. Elsevier Science Publishers,
  Amsterdam, 1984.

\bibitem{Woxx}
Hugh Woodin.
\newblock {\em {The Axiom of Determinacy, Forcing Axioms and the Nonstationary
  Ideal}}, volume~1 of {\em de Gruyter Series in Logic and its Applications}.
\newblock De Gruyter, in press.

\end{thebibliography}

\shlhetal
\end{document}